\documentclass[a4paper,10pt]{article}
\pdfoutput=1 
\usepackage[top=1in,bottom=1in,left=0.75in,right=0.75in]{geometry}
\usepackage{amsmath}
\usepackage{amssymb}
\usepackage{amsthm}
\usepackage{mathtools}
\usepackage[only,llbracket,rrbracket]{stmaryrd}
\usepackage[scr=boondoxo]{mathalfa} 
\usepackage[english]{babel}
\usepackage{graphicx}
\usepackage[font=small,figurewithin=section]{caption}
\usepackage{subcaption}
\usepackage{float}
\usepackage[numbers]{natbib} 
\usepackage{hyperref}
\usepackage[all]{xy}
\usepackage{multirow}
\usepackage{array}
\usepackage{enumitem}
\usepackage{tikz}
\usetikzlibrary{matrix,arrows,decorations.pathmorphing}
\usepackage{color}
\usepackage{authblk}
\usepackage[mathscr]{euscript}

\usepackage{relsize}

\numberwithin{equation}{section}
\makeatletter
\g@addto@macro\bfseries{\boldmath}
\makeatother

\hypersetup{colorlinks=true,citecolor=black,filecolor=black,linkcolor=black,urlcolor=black}
\hypersetup{pdfstartview=FitB,pdfpagemode=UseNone}

\setlist{nolistsep}

\newcolumntype{L}[1]{>{\raggedright\let\newline\\\arraybackslash\hspace{0pt}}m{#1}}
\newcolumntype{C}[1]{>{\centering\let\newline\\\arraybackslash\hspace{0pt}}m{#1}}
\newcolumntype{R}[1]{>{\raggedleft\let\newline\\\arraybackslash\hspace{0pt}}m{#1}}
\newcolumntype{N}{@{}m{0pt}@{}}
\SetSymbolFont{stmry}{bold}{U}{stmry}{m}{n}
\usepackage{derivative}

\DeclareMathOperator{\Sym}{Sym}

\newcommand{\MS}{\mathbb{S}}

\newcommand{\R}{\mathbb{R}} 
\newcommand{\HSo}{H} 
\newcommand{\Leb}{L} 

\newcommand{\idmatrix}{\textup{\uppercase\expandafter{\romannumeral 1}}}
\newcommand{\mesh}{\mathcal{T}} 
\newcommand{\T}{\mathsf{T}} 
\newcommand{\bdry}{\partial} 
\newcommand{\Tr}{\mathrm{Tr}} 

\newcommand{\bveps}{\boldsymbol{\varepsilon}}

\newcommand{\bsigma}{{\mathlarger{\boldsymbol{\sigma}}}}

\DeclareMathOperator{\curl}{curl}
\let\div\relax 
\DeclareMathOperator{\div}{div}

\newcommand{\idtens}{\mathbf{I}} 

\newtheoremstyle{boldremark}
    {\dimexpr\topsep/2\relax} 
    {\dimexpr\topsep/2\relax} 
    {}          
    {}          
    {\bfseries} 
    {.}         
    {.5em}      
    {}          
\theoremstyle{definition}
\newtheorem*{definition*}{Definition}
\theoremstyle{plain}

\theoremstyle{boldremark}
\newtheorem{remark}{Remark}

\makeatletter
\newcommand{\dashto}[1][2pt]{
  \settowidth{\@tempdima}{${}\rightarrow{}$}
  \makebox[\@tempdima]{${}\rightarrow{}$}
  \makebox[-\@tempdima]{\hspace{-0.1\@tempdima}\color{white}\rule[0.5ex]{#1}{1pt}}
  \makebox[\@tempdima]{}
  }
\makeatother
\let\tilde\widetilde

\newcommand{\mcF}{\mathcal{F}}

\newcommand{\mcO}{\mathcal{O}}

\newcommand{\sfC}{\mathsf{C}}

\newcommand{\bsb}{\boldsymbol{b}}

\newcommand{\bsi}{\boldsymbol{i}}
\newcommand{\bsj}{\boldsymbol{j}}

\newcommand{\bsq}{\boldsymbol{q}}

\newcommand{\bsu}{\boldsymbol{u}}
\newcommand{\bsv}{\boldsymbol{v}}
\newcommand{\bsw}{\boldsymbol{w}}

\newcommand{\scV}{\mathscr{V}}

\newcommand{\fkv}{\mathfrak{v}}

\newcommand{\nml}{\mathbf{n}}

\newcommand{\bfx}{\mathbf{x}}

\usepackage{mathtools}

\makeatletter
\DeclareRobustCommand\widecheck[1]{{\mathpalette\@widecheck{#1}}}
\def\@widecheck#1#2{%
    \setbox\z@\hbox{\m@th$#1#2$}%
    \setbox\tw@\hbox{\m@th$#1%
       \widehat{%
          \vrule\@width\z@\@height\ht\z@
          \vrule\@height\z@\@width\wd\z@}$}%
    \dp\tw@-\ht\z@
    \@tempdima\ht\z@ \advance\@tempdima2\ht\tw@ \divide\@tempdima\thr@@
    \setbox\tw@\hbox{%
       \raise\@tempdima\hbox{\scalebox{1}[-1]{\lower\@tempdima\box
\tw@}}}%
    {\ooalign{\box\tw@ \cr \box\z@}}}
\makeatother

\newcommand{\bfb}{\mathbf{b}}
\newcommand{\bfs}{\mathbf{s}}
\newcommand{\bff}{\mathbf{f}}
\newcommand{\bfd}{\mathbf{d}}

\newcommand{\bfM}{\mathbf{M}}

\newcommand{\bfK}{\mathbf{K}}

\newcommand{\bsfC}{\boldsymbol{\sfC}}

\newcommand{\bzero}{\mathbf{0}}

\newcommand{\tend}{t_\text{end}}
\newcommand{\interface}{\Gamma_{\textrm{int}}}
\newcommand{\nmlint}{\nml_{\textrm{int}}}
\newcommand{\iapp}{I_{\textrm{app}}}

\newcommand{\ibv}{I_{\textrm{BV}}}
\newcommand{\kbv}{k_{\textrm{BV}}}
\newcommand{\csmax}{c_{s,\text{max}}}
\newcommand{\ocp}{\Phi_{s,\textrm{open}}}
\newcommand{\ocpsa}{\Phi_{sa,\textrm{open}}}
\newcommand{\ocpsc}{\Phi_{sc,\textrm{open}}}
\newcommand{\csamax}{c_{sa,\text{max}}}
\newcommand{\cscmax}{c_{sc,\text{max}}}
\newcommand{\Gpos}{\Gamma_{\text{cc}+}}
\newcommand{\Gneg}{\Gamma_{\text{cc}-}}
\newcommand{\Gtop}{\Gamma_{\text{top}}}
\newcommand{\Gbot}{\Gamma_{\text{bottom}}}
\newcommand{\soca}{\check{c}_{sa}}
\newcommand{\socc}{\check{c}_{sc}}

\usepackage{tikz}
\usetikzlibrary{shapes.geometric, shapes.misc, arrows}
\tikzstyle{startstop} = [rectangle, rounded corners, minimum width=1.8cm, minimum height=0.5cm,text centered, draw=black, fill=red!30]
\tikzstyle{io} = [trapezium, trapezium left angle=70, trapezium right angle=110, minimum width=1.8cm, minimum height=0.5cm, text centered, draw=black, fill=blue!30]
\tikzstyle{process} = [rectangle, minimum width=1.8cm, minimum height=0.5cm, text centered, draw=black, fill=orange!30]
\tikzstyle{decision} = [diamond, minimum width=1.8cm, minimum height=0.5cm, text centered, draw=black, fill=green!30]
\tikzstyle{forloop} = [chamfered rectangle, chamfered rectangle xsep = 0.5cm, text centered, draw=black, fill=teal!30]
\tikzstyle{cont} = [circle, minimum width=0.3cm,draw=black,fill=teal!30]
\tikzstyle{arrow} = [thick,->,>=stealth]

\usepackage{listings}
\usepackage{lscape}
\usepackage{xfrac}
\usepackage{xcolor}
\usepackage{cancel}

\def\be{\begin{equation}}
\def\ee{\end{equation}}
\def\ba{\begin{array}}
\def\ea{\end{array}}
\def\bea{\begin{eqnarray}}
\def\eea{\end{eqnarray}}
\def\beas{\begin{eqnarray*}}
\def\eeas{\end{eqnarray*}}

\begin{document}
%
\title{High-order transient multidimensional simulation of a thermo-electro-chemo-mechanical model for Lithium-ion batteries}
\author[1]{Jaime Mora-Paz\thanks{e-mail: jaimed.morap@konradlorenz.edu.co}}
\affil[1]{Fundación Universitaria Konrad Lorenz, Bogota DC, Colombia}
\date{\vspace{-10mm}}
%
\maketitle
\renewcommand{\abstractname}{\large Abstract}
\begin{abstract}
\small 
\noindent
We build a transient multidimensional multiphysical model based on continuum theories, involving the coupled mechanical, thermal and electrochemical phenomena occurring simultaneously in the discharge or charge of lithium-ion batteries. The process delivers a system of coupled nonlinear partial differential equations. Besides initial and boundary conditions, we highlight the treatment of the electrode-electrolyte interface condition, which corresponds to a Butler-Volmer reaction kinetics equation. We present the derivation of the strong and weak forms of the model, as well as the discretization procedure in space and in time. The discretized model is computationally solved in two dimensions by means of a finite element method that employs $hp$ layered meshes, along with staggered second order semi-implicit time integration. The expected error estimate is of higher order than any other similar work, both in space and in time. A representative battery cell geometry, under distinct operating scenarios, is simulated. The numerical results show that the full model allows for important additional insights to be drawn than when caring only for the electrochemical coupling. Considering the multiphysics becomes more important as the applied current is increased, whether for discharge or for charge. Our full model provides battery design professionals with a valuable tool to optimize designs and advance the energy storage industry.

\noindent\textbf{Keywords: lithium-ion batteries, multiphysics, coupled PDEs, finite elements, semi-implicit time integration} 
\end{abstract}
\section{Introduction}
\label{sec:Introduction}
\subsection{Motivation}
Rechargeable batteries have become ubiquitous in the past decades, and their research, development and market size will keep growing in the future \cite{fichtner2022rechargeable,fleischmann2023battery,precedence2023}. The need of storing energy from alternate sources and/or for powering electric vehicles and many consumer goods pushes the research in battery technology and, consequently, the demand for modeling and simulation tools that deliver additional insight and even predict the functioning of new battery systems, structures, materials and applications \cite{fichtner2022rechargeable,granholm2021national}. The present paper aims at a new effort to obtain a high-fidelity multiphysics simulation for Lithium-ion batteries. Such an effort must begin by understanding and appropriately modeling the coupled phenomena that unfold inside a battery.

Modern-day Lithium-ion batteries comprise a single or multiple cells, which is the basic unit inside a battery where the energy storage occurs \cite{linden2002handbook}. Each battery cell mainly consist of two solid electrodes, which can be porous, that are named anode (or negative electrode) and cathode (or positive electrode), with a separator between them. In the macroscale, conventional cell designs comprise a cathode built with a Lithium compound; a carbon-based anode, made, for instance, of graphite, which is a material that offer stability and sufficient space to allocate Lithium ions; and a porous membrane immersed in an electrolytic solution that is able to conduct ions but not electrons. On the opposite side of the electrodes lie metallic foils known as current collectors. Typically, the negative current collector (in contact with the anode) is of copper and the positive one (in contact with the cathode) is of aluminum \cite[pp4]{bai_chemo-mechanical_2021}.
\begin{figure}
    \centering
    \includegraphics[width=0.48\textwidth]{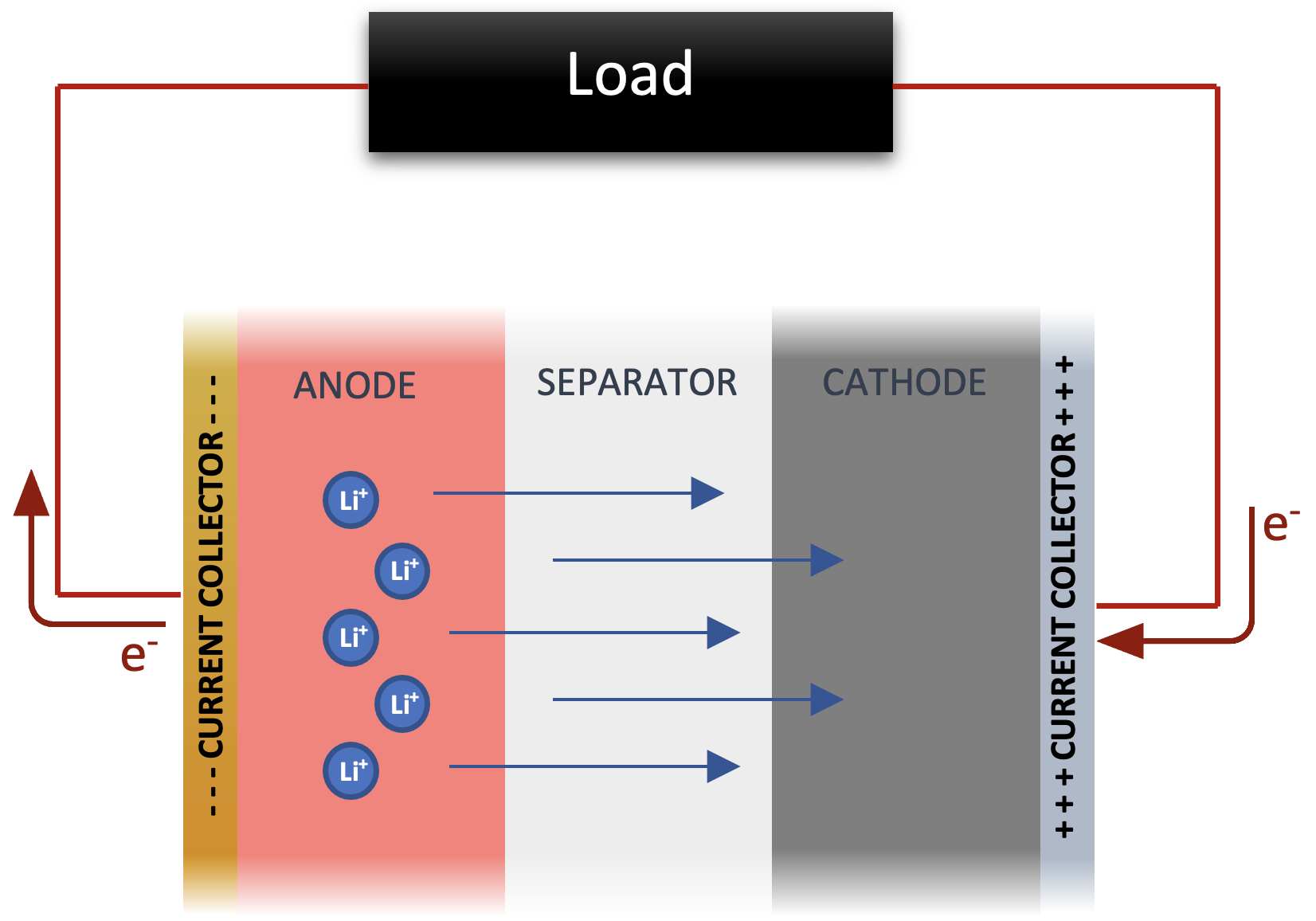}
    \includegraphics[width=0.48\textwidth]{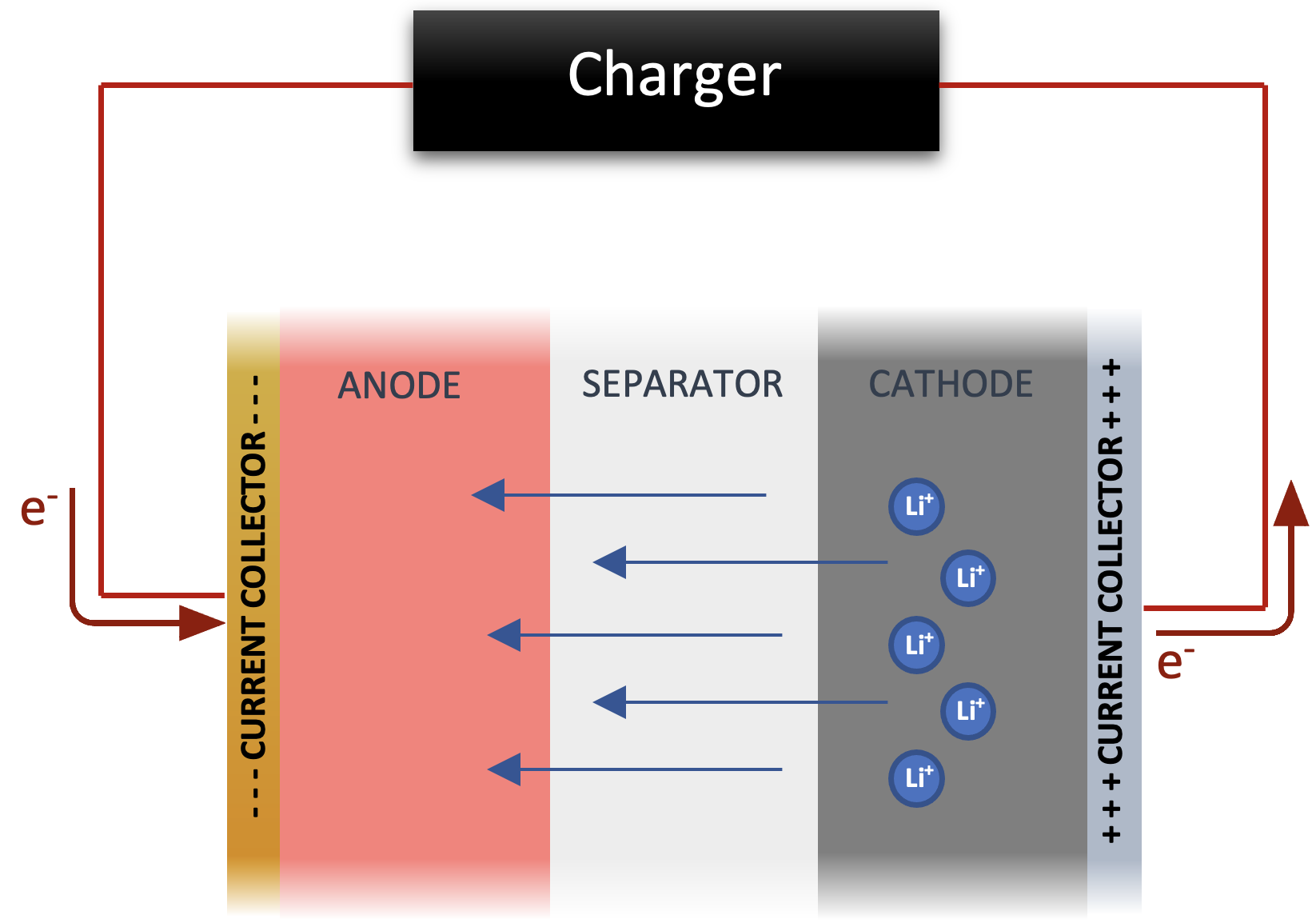}
    \caption{(Left) Diagram of a Lithium-ion battery discharge. (Right) Diagram of a Lithium-ion battery charge.}
    \label{fig:charge-discharge}
\end{figure}

We now explain in a few words the working principles of this type of batteries, with the help of the diagrams included in Figure \ref{fig:charge-discharge}. Consider the process of discharge (left). Given certain electric load, it demands an outflow of current at the positive current collector, so that electrons go from the load onto the positive side of the battery. The abundance of negative charge (of the electrons) attract positively charged particles, which are the Lithium ions produced by the electrochemical reaction between the electrodes and the electrolyte upon a certain difference of electric potential across their interface. A high concentration of Lithium in the anode (i.e., the battery is charged) maintains this process for a long time.
On the other hand, during the charging process (right),  an inflow of current is imposed, and the induced electrochemical reaction generates motion of Lithium particles from the cathode to the anode. If there is high concentration of ions at the cathode (i.e., the battery is discharged), more time or more current is required to fully charge the battery. The flow and allocation of ions into the electrodes is a mass transport phenomenon that is also known as insertion or \emph{intercalation}, a reversible process \cite[\S 35.2.1]{linden2002handbook}.

\begin{figure}
    \centering
    \includegraphics[width=0.9\textwidth]{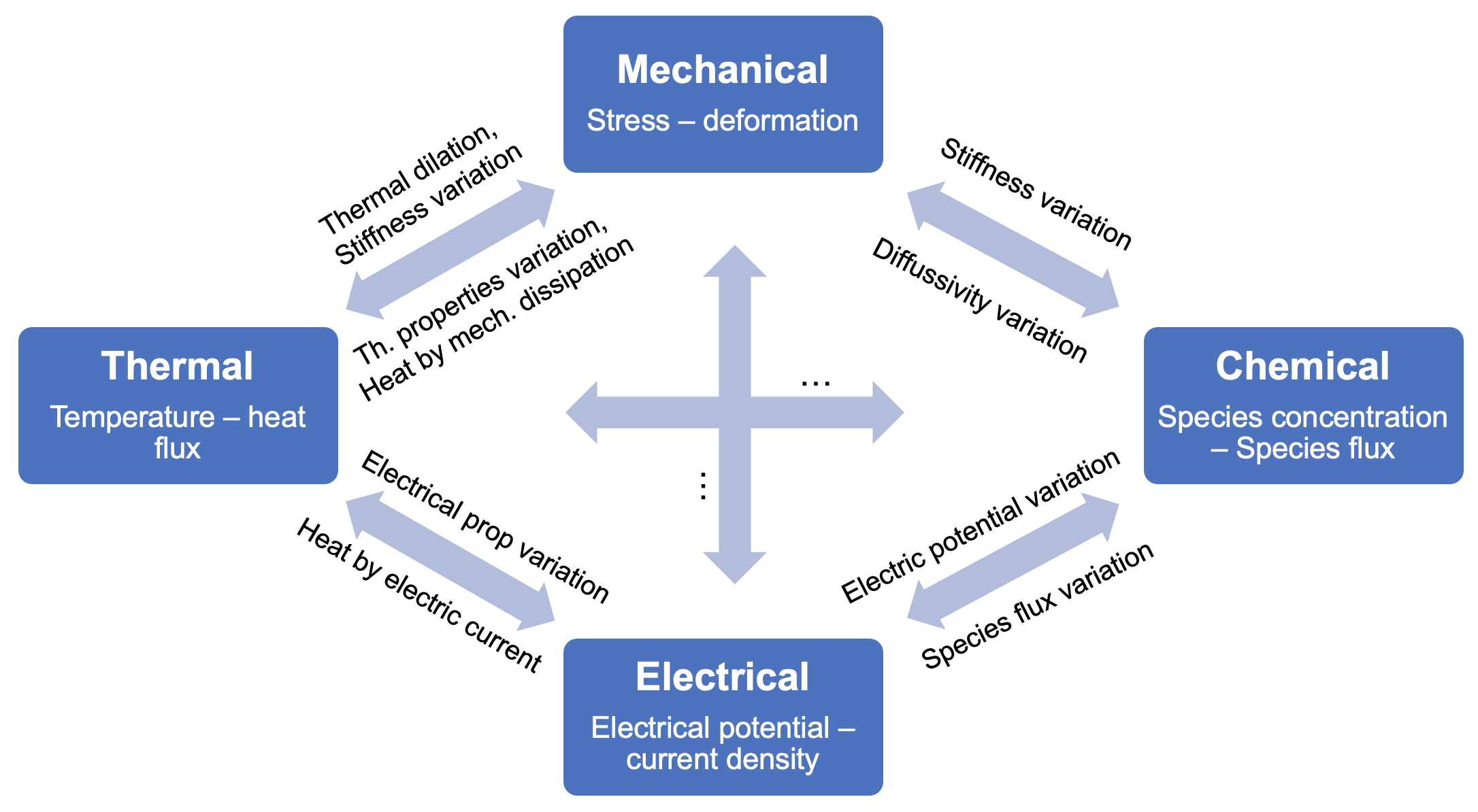}
    \caption{Multiphysical interactions occuring within a Lithium-ion battery during charge/discharge.}
    \label{fig:interactions}
\end{figure}

During charging/discharging of Lithium-ion batteries, apart from the obvious electrical and chemical phenomena that drive the process (electric charge transport and mass diffusion) there are additional physical processes simultaneously going on. For instance, the mentioned electrochemical processes induce thermal effects (heating by Ohm's law) as well as mechanical effects (dilation due to temperature and concentration change). Mechanical stresses, in turn, can perturb the mass diffusion through the change of diffusivity, thus altering too the ion transport behavior \cite{gatica_analysis_2018,gritton_using_2017,anand_continuum_2020,miranda_effect_2019,carlstedt_coupled_2022}. Of course, this is not an exhaustive list of the interactions. Battery charge/discharge is consequently a strongly coupled multiphysical system. Figure \ref{fig:interactions} illustrates in short the variables of interest and the distinct mechanisms of interaction among them.

In what follows, we quickly review some lithium-ion-battery coupled multi-physics models based on continuum theories, which help giving an idea of where the state of the art lies. With this context, we shortly declare the goals and scope and of this research paper.

\subsection{Background briefing}

Seminal works by Newman, Doyle and Fuller in the 90s \cite{doyle1993modeling,fuller_simulation_1994,doyle_comparison_1996} on the electrochemical model equations for a Lithium-ion battery have been widely adopted by researchers since then. Of key importance to us is the fact that these works provide experimental models for the open-circuit potential of two common electrode materials in this kind of batteries, which are evaluated in the electrolysis kinetics model and that we must utilize in our numerical experiments below. Because of the distinct scales in play when studying battery physics, most work available (including Newman's) involves the interaction of two scales: a microscale, where the electrolytical reaction unfolds across the interfaces of the electrodes and the electrolytic medium; and a macroscale, where the electrodes are porous and a homogenization approach is adopted, that is, the battery cell is regarded as a mixture of three materials, their volume fractions are taken into account for computing the effective material properties and the resulting fluxes of the microscale reaction are brought in as a local source function \cite{trembacki_fully_2015,bai2019two,miranda_effect_2019}. This is nevertheless a modeling path different from ours, as will be seen below.

Multiple models which account for the interactions of two or three of the physical aspects in question exist and have been applied within numerical simulations that can be utilized to analyze various battery designs, regarding geometry and materials. 
For instance, Anand and Govindjee \cite{anand_continuum_2020} illustrates three coupled theories that are partly applicable to our problem: thermoelasticity, species diffusion coupled elasticity, and chemoelasticity, considering small deformations in all cases. 
Yang Bai dedicates his doctoral research to several features of the chemo-mechanical modeling of Lithium-ion batteries, including large deformations \cite{zhao2019review,bai2019two,bai2020chemo,bai_chemo-mechanical_2021}. 
Another case with three physical attributes, namely, a thermo-electrochemical coupling, is that by Miranda et al. \cite{miranda_effect_2019}, who carry out a numerical comparison of the thermal effects from different cathode materials and cell geometries. A second example is found in \cite{northrop2015efficient}, focusing on simulation efficiency. Review papers on the coupled modeling of these three aspects are also available \cite{park2010review,zhao2019review}.

Regarding the interaction of all four physical attributes, we have found two references by Carlstedt and Asp \cite{carlstedt2019thermal,carlstedt_coupled_2022}, about studying the coupled performance of structural batteries with composite materials. The 2019 paper studies thermal and diffusion-induced cyclic stresses, simplified in space to one-dimensional (axisymmetric) equations. However, the 2022 publication is considerably more ambitious and thorough, to a level that is comparable to the scope of the present work. That paper has been devoted to rigorously derive and simulate galvanostatic cycles with a model that involves our four aspects of interest: thermal, mechanical, electrical, chemical. This model is specialized to a structural-battery-electrode design, and solved in two-dimensions through Galerkin methods and implicit Euler time-stepping.

Finally, at least two major commercial finite element packages, namely, ANSYS and COMSOL, are offering battery simulation modules that incorporate Newman-Doyle-Fuller's equations, with the possibility of combining those with additional physics attributes and 3D capabilities, besides all the robust geometric and numerical utilities equipped with those programs. Because of all this, they are advertised as advanced solutions for design and analysis of this type of systems \cite{ansys_batteries,comsol_batteries}. 






\subsection{Goals and scope}
We intend to derive and use a multi-way coupled PDE-based multiphysics model, built upon classical continuum theories, by means of high order Galerkin discretizations and time-stepping methods, in a staggered, semi-implicit fashion, without homogenization or multi-scale approaches.

Besides demonstrating that the numerical methods of choice are a good fit for the system of non-linear equations at hand, we take a glimpse on results' difference when accounting for all thermo-electrochemo-mechanical effects in contrast with the more usual electrochemical models available in the literature.

All the material characterization and battery configurations are drawn from the literature, and suitable assumptions are made as well in order to complete the model. The model herein developed is nonetheless original and at the highest degree of complexity available, to the best of our knowledge.

\subsection{Organization of the document}
Section \ref{sec:model} presents the battery multi-physics model derivation, both in its strong and weak forms. Next, the discretization approach is detailed in Section \ref{sec:method}. Numerical experiments with real parameter and geometries are included in Section \ref{sec:results}. We conclude in Section \ref{sec:conclusions} with final remarks.
\section{Multiphysics Li-ion battery model} 
\label{sec:model}
The present section begins with a description of the geometry, and the basic assumptions of our modeling effort. As explained in the motivation subsection, there are four aspects of the multi-physical phenomenon of battery charging/discharging that we attempt to mathematically couple. The electrical and chemical aspects are jointly introduced in the first place, followed by the thermal (energy balance) and mechanical (momentum balance) aspects. All those equations are gathered in a system of PDEs. Equally important is the introduction of the reaction kinetics of the electrochemical processes at the electrode-electrolyte interface. With that, we bring together the remaining pieces for completing the mathematical model: initial, boundary and interface conditions. As the present model is solved through finite element methods, the detailed derivation of a suitable weak problem is also given.

\subsection{Geometric representation of the battery}

We begin by stating important hypotheses and assumptions of our work:
\begin{itemize}
    \item the present model is purely deterministic;
    \item the model domain is circumscribed to a small portion of the battery cell geometry, namely, a representative unit of the micro-structure, which is supposed to repeat numerously inside a real-size battery cell;
    \item the separator is modeled simply as a volume filled with electrolyte, with no other material present, and the contact with the electrodes is direct and occurs on the electrode-separator interface only;
    \item nor gravitational effects neither fluid dynamics phenomena (e.g., in the electrolyte) are accounted for;
    \item the electrodes are not modeled as porous or composite materials, but voidless isotropic solids.
\end{itemize}

\begin{figure}
    \centering
    \includegraphics[width=0.35\textwidth]{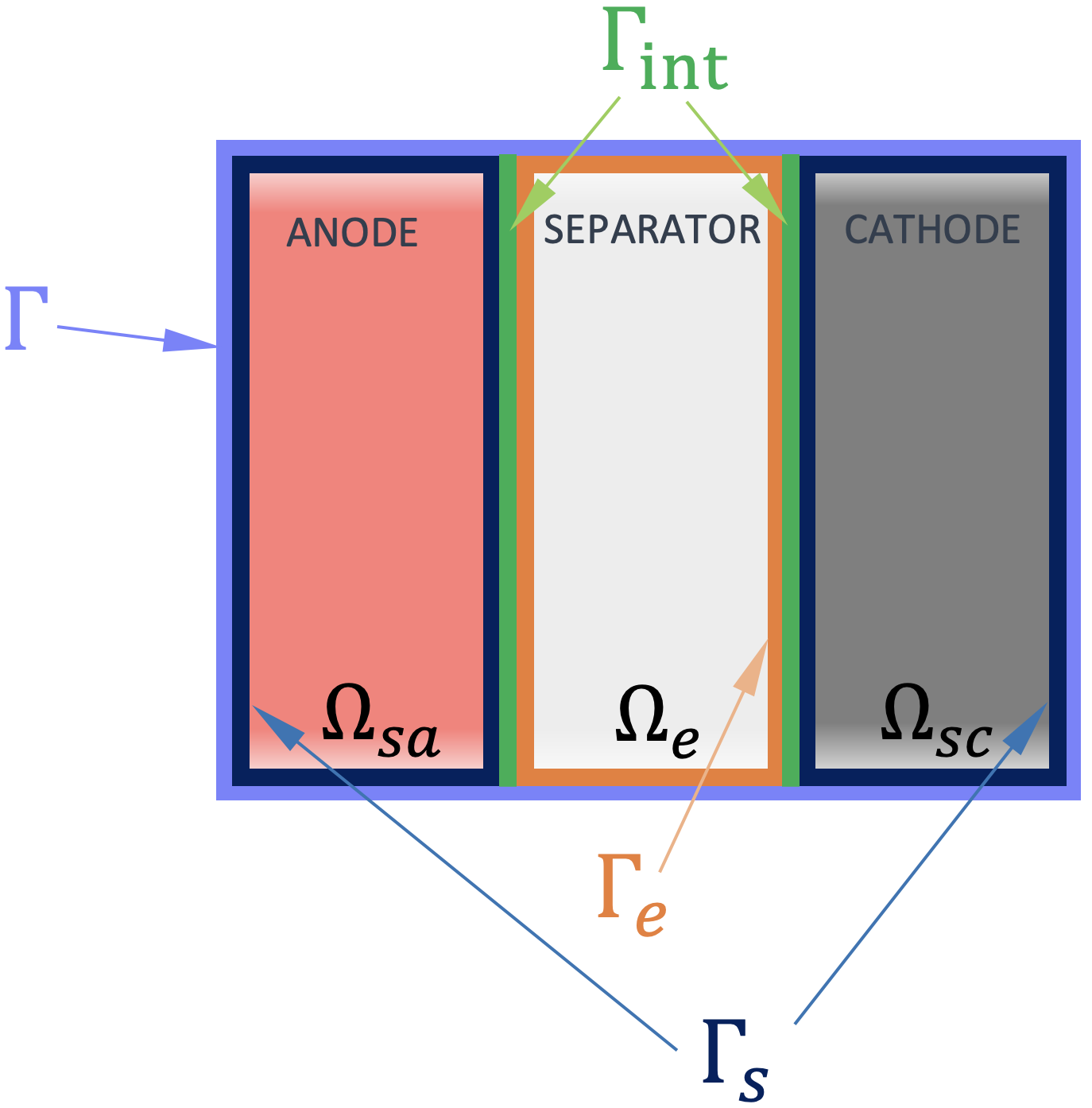}
    \label{fig:subdomains-diagram}
    \caption{Battery diagram showing subdomains, boundaries and the interface.}
\end{figure}

The representative geometry of the battery internal structure herein modeled lies in $\Omega\subset\R^2$, which is a domain (open, connected, bounded set) with Lipschitz boundary, $\Gamma = \bdry\Omega$. The domain has three disjoint subdomains or parts, the (solid) anode subdomain 
$\Omega_{sa}$, the (solid) cathode subdomain $\Omega_{sc}$, and the electrolyte subdomain $\Omega_e$. Each of these is a Lipschitz domain as well.
Subset $\Omega_s$ is the union of the two electrodes' subdomains, but it is not a domain per se, since it is not connected (not even its closure is connected), as the  two electrodes cannot be in contact. We have then
\begin{align*}
    \Omega_s    &= \Omega_{sa} \cup \Omega_{sc}, \vspace{0.2cm}\\
    \overline{\Omega}      &= \overline{\Omega_{s}} \cup \overline{\Omega_{e}}
\end{align*}
In order to deal with three components, we use an index $\iota$ as subscript in many symbols below, relating to geometry, material properties and unknown fields. When $\iota=s$, we are dealing with either of the solid electrodes. If more precision is required, then we set $\iota=sa$ or $\iota=sc$, for the anode or the cathode, respectively. The advantage of using the $s$ subscript is that it helps avoiding duplicated expressions for the two electrodes. Finally, if $\iota=e$, we are dealing with the electrolyte. This is next summarized:
\begin{equation*}
\iota = 
\left\{
    \begin{array}{rll}
        s   & \text{solid} 
            & \left\{\begin{array}{rl} sa & \text{solid anode} \\ sc & \text{solid cathode} \end{array}\right. 
            \vspace{0.1cm}\\
        e   & \text{electrolyte}
    \end{array}
\right.            
\end{equation*}

Notice that the union of the boundaries yields
$\bdry\Omega_s \cup \bdry\Omega_e = \Gamma \cup \interface$,
where $\interface$ is the interface between the electrodes and the electrolyte ($\interface$ has two disjoint components). 
Boundaries $\bdry\Omega_s$ and $\bdry\Omega_e$ can be decomposed similarly:
$\bdry\Omega_s=(\Gamma\cap\bdry\Omega_s)\cup \interface$ and $\bdry\Omega_e=(\Gamma\cap\bdry\Omega_e)\cup \interface$. 
If we write $\nml$ for the normal unit vector, it represents the vector pointing out of the domain in which we are working. However, as we assume that $\interface$ is assigned a unique orientation, its unit normal vector $\nmlint$, always points outward from either electrode ($\Omega_s$) into the electrolyte ($\Omega_e$). This is taken into consideration in the interface integrals.

\subsection{Governing equations}

\subsubsection{Charge/Current density conservation.}
Within each medium, Gauss' law will govern the conservation of current density $\bsi_\iota$ throughout medium $\iota$ ($\iota=s$ for solid electrodes, $\iota=e$ for the electrolyte): 
\begin{equation}
\nabla\cdot \bsi_\iota = 0.
\label{eq:gauss}
\end{equation}
Symbol $\nabla\cdot$ is the divergence operator. Here, the current density constitutive models are based on Ohm's law, adding a concentration modification for the $e$ case:
\begin{equation}
\bsi_\iota =
\left\{
    \begin{array}{rlrl}
        -\gamma_s\nabla \phi_s,
            & \iota = s;
            \vspace{0.2cm}\\
        -\kappa_e\nabla\phi_e - \kappa_D\nabla \ln c_e,
            & \iota = e;
            \vspace{0.2cm}\\
    \end{array}
\right.
\label{eq:current-density}
\end{equation}
where $\nabla$ is the gradient operator. The material properties in \eqref{eq:current-density} are: electrical conductivity in the solid $\gamma_s$, ionic conductivity in the electrolyte $\kappa_e$, diffusional conductivity in the electrolyte $\kappa_D$. The latter is given by
\begin{equation}
    \kappa_D = -\frac{2R\theta \kappa_e}{F}( 1 - t_+ ),
    \label{eq:diffusional_cond}
\end{equation}
where $t_+$ is known as the ion transference number, $\theta$ is the temperature (another primal variable, whose evolution is described below), and $R$ and $F$ are physical constants: the universal gas constant and Faraday's constant, respectively. As $0<t_+<1$, notice that $\kappa_D$ is a negative number. Moreover, notice that $\nabla \ln c_e$ is equal to $\frac{1}{c_e}\nabla c_e$.

In summary, the primal variables for this physical aspect are the electric potential in the solid electrodes $\phi_s(t,\bfx)$ and the electric potential in the electrolyte $\phi_e(t,\bfx)$, and the PDEs that govern them are
\begin{align}
    -\nabla\cdot \left( \gamma_s\nabla \phi_s \right)
     &=0
    \label{eq:solid-potential-pde} \\
    -\nabla\cdot \left( \kappa_e\nabla\phi_e + \kappa_D\nabla \ln c_e\right)
     &=0
    \label{eq:elyte-potential-pde}
\end{align}
\subsubsection{Mass conservation / Lithium ions balance.}
The species (Li+ ions) concentration balance in each medium is represented by $c_\iota(t,\bfx)$ ($\iota=s$ for solid, $\iota=e$ for the electrolyte) and follows the diffusion equation
\begin{equation*}
    \dot{c}_\iota\,+\,\nabla\cdot\bsj_\iota\,=\,0\, ,
\end{equation*}
where $\bsj_\iota$ is the species flux in medium $\iota$ \cite{trembacki_fully_2015}. The dot on $c_\iota$ means partial derivative with respect to time $t$.
Notice that we are not taking into account any advective transport that could occur in the electrolytic medium. 
Indeed, in the present work, the disregard for the mechanical behavior in this part of the battery facilitates as well the constitutive law for $\bsj_e$, yielding
\begin{equation}
    \bsj_e = -D_e\nabla c_e,
\end{equation}
which is Fick's first law \cite{park2010review}, assuming a constant mass diffusivity of the electrolyte $D_e$.

The mechanical aspect in the solid electrodes, to be described below, has effects over the electrochemical phenomena, through a so called stress-assisted  diffusion process. 
Although there are multiple ways to model this coupling, see \cite{anand_continuum_2020} and \cite{gatica_analysis_2018}, here we follow \cite{gritton_using_2017}, \S3.1, by picking a convenient constitutive model that incorporates the stress into the picture in the following form:
\begin{align}
\bsj_s &= -D_s(c_s,\pi)\nabla c_\iota \,,
\label{eq:species-flux-law} \\
\pi &=-\frac{1}{3}\Tr(\bsigma) \,,
\label{eq:pressure} \\
D_s (c_s,\pi) &=
\left\{ \begin{array}{cc}
D_{s,0}\,e^{\alpha_{D}c_s/c_{s,\max}}                         & \text{if}\ \pi\leq 0, \\
D_{s,0}\,e^{\alpha_{D}c_s/c_{s,\max}-\beta_{D}\pi/\pi_{\max}} & \text{if}\ 0 < \pi < \pi_{\max}, \\
D_{s,0}\,e^{\alpha_{D}c_s/c_{s,\max}-\beta_{D}}               & \text{if}\ \pi\geq \pi_{\max},
\end{array}
\right.
    \label{eq:stress-diffusivity}
\end{align}
where $\pi$ is the hydrostatic pressure, $\bsigma$ is the Cauchy stress tensor (detailed below), $D_{s,0}$ is a reference level diffusivity for each electrode, and $\alpha_D,\beta_D,\pi_{\max},$ are suitable positive real constants. 
As a disclaimer, in \cite{gritton_using_2017} the model \eqref{eq:stress-diffusivity} is applied to a silicon electrode, which differs from the carbon-based electrodes of the numerical experiments below, but we have chosen it as a mathematically convenient option. 
Unlike other constitutive laws in the literature, this choice has as an advantage, in the author's opinion, that the hydrostatic pressure (hence the stress field) needs no differentiation in the final weak formulation, thus being consistent with the functional spaces used below\footnote{At a fixed time $t$, the displacement field $\bsu\in \HSo^1(\Omega_s;\R^3)$, hence the stress field $\bsigma\in \Leb^2(\Omega_s;\MS^{3 \times 3})$, and consequently doesn't have well-defined derivatives in space.}.

Finally, the mass conservation PDEs for the primal variables $c_s(t,\bfx)$ and $c_e(t,\bfx)$ are
\begin{align}
    \dot{c}_s-\nabla\cdot (D_s(c_s,\pi)\nabla c_s) = 0,
    \label{eq:solid-conc-pde} \\
    \dot{c}_e-\nabla\cdot (D_e\nabla c_e) = 0.
    \label{eq:elyte-conc-pde}
\end{align}

\subsubsection{Energy balance.}
By taking the first law of thermodynamics, neglecting any work by a mechanical source, and representing the internal energy of each medium as proportional to the temperature $\theta$, the energy balance equation has the form \cite{miranda_effect_2019}:
\begin{equation*}
    \rho_\iota C_{v,\iota} \dot{\theta}
         + \nabla\cdot \bsq_\iota 
            = Q_\iota;
    \label{eq:energy-balance}
\end{equation*}
here, $\rho_\iota$ is the density of medium $\iota$, while $C_{v,\iota}$ is its specific heat capacity at constant volume; in the second term, $\bsq_\iota$ is the heat flux in medium $\iota$, related to the temperature spatial variation through the first Fourier's law:
\begin{equation}
\bsq_\iota=-\lambda_\iota\nabla\theta,
\label{eq:Fourierlaw}
\end{equation}
with $\lambda_\iota$ being the corresponding (constant) thermal conductivity coefficient. 
On the right-hand side (RHS) of this energy balance, we have the volumetric heat source $Q_\iota$ \cite{miranda_effect_2019,carlstedt_coupled_2022}, which hereinafter consists of the Ohmic heat generation only,
\begin{equation}
        Q_\iota = \bsi_\iota \cdot \nabla \phi_\iota.
        \label{eq:bulk_heat_source}
\end{equation}
The PDE that governs the primal variable $\theta(t,\bfx)$ across medium $\iota$ yields
\begin{equation}
    \rho_\iota C_{v,\iota} \dot{\theta} - \nabla\cdot (\lambda_\iota\nabla\theta)
    = \bsi_\iota \cdot \nabla \phi_\iota
    \label{eq:temperature-pde}
\end{equation}

\subsubsection{Linear momentum balance} 
For the mechanical characterization of the solid electrodes, we neglect accelerations and use the quasi-static equation for linear momentum balance \cite{anand_continuum_2020}:
\begin{equation*}
    \nabla\cdot\bsigma\,+\,\bsb\, =\, \bzero ,
\end{equation*}
where $\bsigma$ is the Cauchy stress tensor, and $\bsb$ is a vector of volumetric body forces. Hereinafter, we assume that $\bsb\equiv\bzero$ \cite{zhao2019review,carlstedt_coupled_2022}. We shall incorporate generalized Hooke's law as the stress constitutive law, namely,
\begin{equation*}
    \bsigma\,= \bsfC:\bveps^{\mathrm{me}},
\end{equation*}
with $\bsfC$ representing the fourth-order elasticity tensor. Hereinafter, this tensor corresponds to that of a compressible, isotropic and homogeneous linearly elastic material, so it will be fully defined with just two constants, e.g., the shear modulus $G$ and the bulk modulus $K$. 

The mechanical strain $\bveps^{\mathrm{me}}$ is defined by \cite{anand_continuum_2020,carlstedt_coupled_2022}:
\begin{equation}
    \bveps^{\mathrm{me}} = \bveps - \bveps^{\mathrm{th}} - \bveps^{\mathrm{ch}} ,
    \label{eq:me-strain}
\end{equation}
with the full strain tensor being
\begin{equation}
    \bveps = \nabla_{_{\Sym}}\bsu = \frac{1}{2}\left(\nabla\bsu+\nabla^\T\bsu\right),
    \label{eq:full-strain}
\end{equation}
that is, the symmetric part of the gradient of the displacement vector field $\bsu$; moreover, the thermal and chemical strains are defined by
\begin{align}
    \bveps^{\mathrm{th}} &= \alpha (\theta - \theta_{0})\idtens, \label{eq:th-strain} \\
    \bveps^{\mathrm{ch}} &= \omega (c_s - c_{s,0}) \idtens, \label{eq:ch-strain}
\end{align}
where $\alpha$ and $\omega$ are the thermal dilation coefficient and the concentration dilation coefficient, respectively; $\theta_0$ is the reference (uniform) temperature, and $c_{s,0}$ is the reference species concentration in the solid, which is taken as a constant in each electrode. 

For completeness, we show the application of the elasticity tensor $\bsfC$ on the mechanical strain tensor $\bveps^{\mathrm{me}}$, for the current case of a homogeneous, isotropic, linear elastic solid \cite{anand_continuum_2020}:
\begin{equation}
    \bsfC:\bveps^{\mathrm{me}} =  2 G \bveps^{\mathrm{me}} + (K-2G/3) \Tr(\bveps^{\mathrm{me}}) \idtens
    \label{eq:generalized_Hooke}
\end{equation}
where $G$ and $K$ can be computed from the well-known pair $(E,\nu)$, that is, Young's modulus and Poisson's ratio, as follows:
\begin{equation}
    G = \frac{E}{2(1+\nu)} ,\qquad
    K = \frac{E}{3(1-2\nu)}.
    \label{eq:Lame_parameters}
\end{equation}

By combining the linear momentum balance and definitions \eqref{eq:me-strain}-\eqref{eq:Lame_parameters}, we get
\begin{equation*}
\bsigma\,
= \bsfC:\nabla_{_{\Sym}}\bsu - 3 K \alpha (\theta - \theta_{0})\idtens - 3 K \omega (c_s - c_{s,0})\idtens
\end{equation*}
The primal variable for the mechanical aspect is the displacement $\bsu(t,\bfx)$. The corresponding PDE is reorganized as follows
\begin{equation}
    -\nabla\cdot (\bsfC:\nabla_{_{\Sym}}\bsu) = - 3 K \alpha \nabla(\theta - \theta_{0}) - 3 K \omega \nabla (c_s - c_{s,0}).
    \label{eq:displacement-pde}
\end{equation}

For a two-dimensional application, we utilize the plane strain model, as we are considering a very small plane section of a three-dimensional object many times larger in every direction (including the normal direction), thereby the scenario for plane strain holds. 

\begin{remark}
    We point out that, even though several material properties in the governing equations derived above are assumed to be constant, they could be variable as well. Those include $\rho$, $C_{v,\iota}$, $\lambda_\iota$, $ D_e$, $\gamma_s$, $\kappa_e$, $K$ and $G$. As long as they range between positive bounds, the present framework should still hold.
\end{remark}
\subsubsection{Full system of partial differential equations}
\begin{figure}
    \centering
    \includegraphics[width=0.25\textwidth]{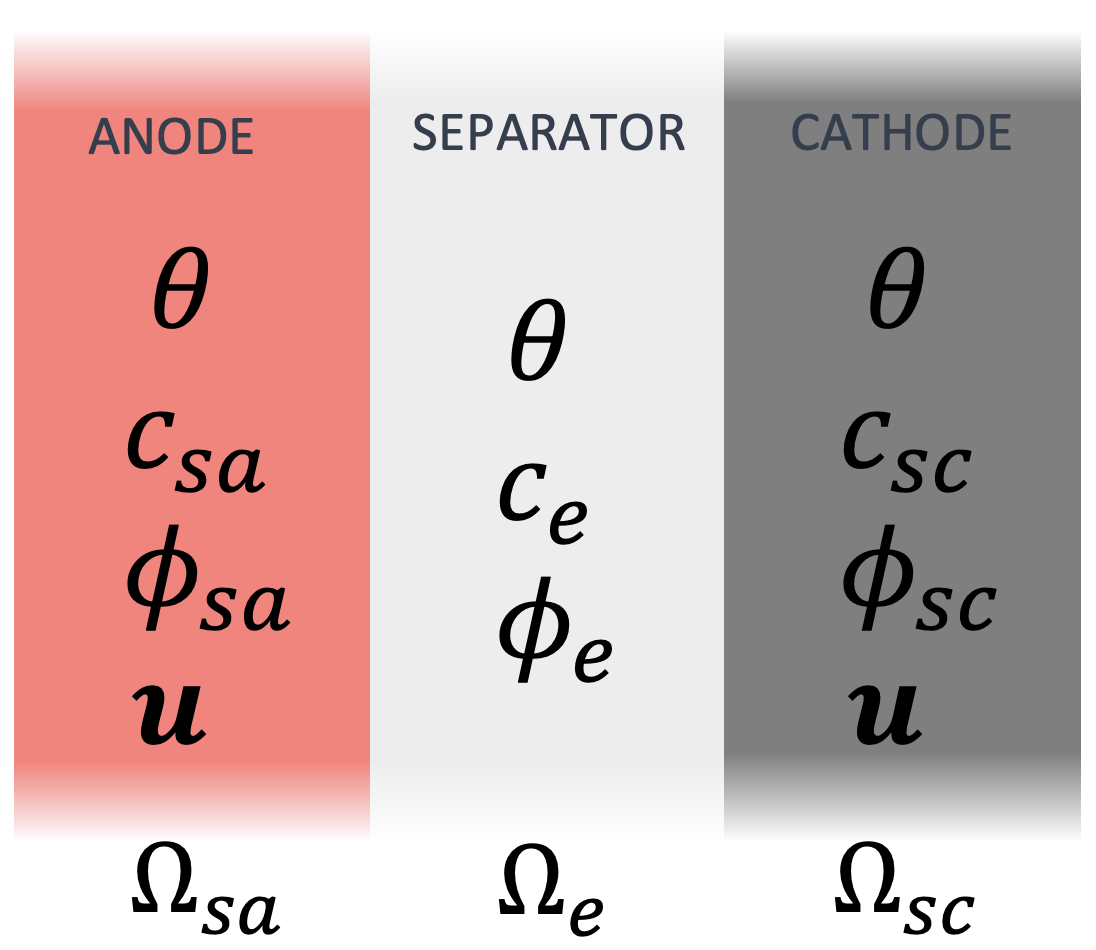}
    \caption{Battery diagram showing the primal unknowns supported in each subdomain.}
    \label{fig:variables-support-diagram}
\end{figure}

The collection of individual equations presented above are gathered in a system of six partial differential equations with six primal unknowns, namely,
\begin{equation}
\left\{
    \begin{array}{rrlrl}
        \rho_\iota C_{v,\iota} \dot{\theta} &
        - \nabla\cdot\left(\lambda_\iota\nabla\theta\right)
            &= \bsi_\iota \cdot \nabla \phi_\iota
            & \text{ in}
            & \Omega\times(0,\tend],  \vspace{0.2cm}\\
        \dot{c_s} &
        - \nabla\cdot \left(D_s(c_s,\pi)\nabla c_s\right)
            & = 0
            & \text{ in}
            & \Omega_s \times(0,\tend] ,  \vspace{0.2cm} \\
        \dot{c_e} &
        - \nabla\cdot \left(D_e\nabla c_e\right)
            & = 0
            & \text{ in}
            & \Omega_e \times(0,\tend],  \vspace{0.2cm}\\
        & -\nabla\cdot \left(\sigma_s\nabla \phi_s\right) 
            &= 0
            & \text{ in}
            & \Omega_s \times(0,\tend],  \vspace{0.2cm}\\
        & -\nabla\cdot\left(\kappa_e\nabla\phi_e\right) 
            & = \nabla\cdot\left(\kappa_D\nabla \ln c_e\right)
            & \text{ in}
            & \Omega_e \times(0,\tend], \\
        & - \nabla\cdot \left( \bsfC : \nabla_{_{\Sym}} \bsu \right) 
            &= -3K \alpha \nabla (\theta - \theta_{0}) 
                 - K \omega \nabla (c_s - c_{s,0})  
            & \text{ in}
            & \Omega_s \times(0,\tend],  \vspace{0.2cm}
    \end{array}
\right.      
\label{eq:pde_system}
\end{equation}
with the unknowns
\begin{enumerate}
    \item Temperature: scalar field $\theta(t,\bfx)$, defined in $\overline{\Omega}\times[0,\tend]$.
    \item Concentration in the solid electrodes: scalar field $c_s(t,\bfx)$, defined in $\overline{\Omega_s}\times[0,\tend]$.
    \item Concentration in the electrolyte: scalar field $c_e(t,\bfx)$, defined in $\overline{\Omega_e}\times[0,\tend]$.
    \item Displacement: vector field $\bsu(t,\bfx)$, defined in $\overline{\Omega_s}\times[0,\tend]$.
    \item Electric potential in the solid electrodes: scalar field $\phi_s(t,\bfx)$, defined in $\overline{\Omega_s}\times[0,\tend]$.
    \item Electric potential in the electrolyte: scalar field $\phi_e(t,\bfx)$, defined in $\overline{\Omega_e}\times[0,\tend]$.
\end{enumerate}
Figure \ref{fig:variables-support-diagram} illustrate how these variables are distributed over the battery subdomains. This system has been reorganized to view the time-dependent equations first. Notice that this is a transient, non-linear, coupled, second-order PDE system. The first three equations in \eqref{eq:pde_system} are parabolic, while the last three are elliptic.

To finish the picture, we must close the model with appropriate boundary conditions (BC) and interface conditions, as well as initial conditions. Due to its central role in setting BC and interface conditions, we first proceed with the introduction of the electrochemical reaction kinetics model.
\subsection{Butler-Volmer Reaction Kinetics}
\label{subsec:butler-volmer}
What happens on the electrode-electrolyte interface, $\interface$, is a key piece in the battery modeling machinery. Indeed, the electrochemical reaction that arises on this interface is the main mechanism that couples the physical processes occurring in the electrodes to those in the electrolyte.
Upon contact between either solid electrode and the electrolytic solution, an electrochemical reaction may arise given certain potential difference across their interface.
Such a reaction produces transport of Lithium ions from one side to another, hence, an electric current normal to the interface is obtained, whose value $\ibv$ is evaluated through the Butler-Volmer equation:
\begin{equation}
    \begin{aligned}
    \ibv & = 2 I_c \sinh\left(\frac{F\eta}{2R\theta}\right)\ , \\
    I_c & = \kbv F\, c_e^{\sfrac{1}{2}}\left(\csmax-c_s\right)^{\sfrac{1} {2}}c_s^{\sfrac{1}{2}}\ ,
\end{aligned}
\label{eq:butler-volmer} 
\end{equation}
where $I_c$ is known as the exchange current density, $\kbv$ is a specific reaction constant, $\csmax$ is the saturation or maximum concentration possible at each electrode, and $\eta$ represents the overpotential:
\begin{equation}
    \eta=\phi_s-\phi_e-\ocp(\check{c}_s) \,.
    \label{eq:overpotential}
\end{equation}
In \eqref{eq:overpotential}, $\ocp$ denotes the so called open circuit potential, an electrode's material property that varies depending on the local \emph{state of charge}, 
\begin{equation}
    \check{c}_s=c_s/\csmax.
    \label{eq:stateofcharge}
\end{equation}
There exist experimentally fitted functions relating the open-circuit potential to the state of charge depending on the material of the electrode, thereby we are required to specify two different functions for the anode and the cathode.

The computation of $\ibv$, $I_c$, $\eta$, $\ocp$ and $\check{c}$ is carried out using the traces of $c_s$, $c_e$, $\phi_s$, $\phi_e$ on $\interface$. Notice that, according to \eqref{eq:butler-volmer}, whenever $\eta$ vanishes, so does $\ibv$ (because $\sinh(0)=0$), meaning that there is no ion transfer at the point. Per \eqref{eq:overpotential}, such local equilibrium holds if $\phi_s-\phi_e=\ocp(\check{c}_s)$.

In addition to quantifying the current density passing from and to the electrodes, Butler-Volmer equation is useful in determining the species flux across the interface, which is essentially proportional to $\ibv$ via Faraday's constant $F$. Furthermore, $\ibv$ together with the overpotential $\eta$ induce an interfacial Ohmic heat source. In the next subsection, these effects are detailed and listed along with the BC.

\subsection{Specific BC and interface conditions}
\begin{figure}
    \centering
    \includegraphics[width=0.35\textwidth]{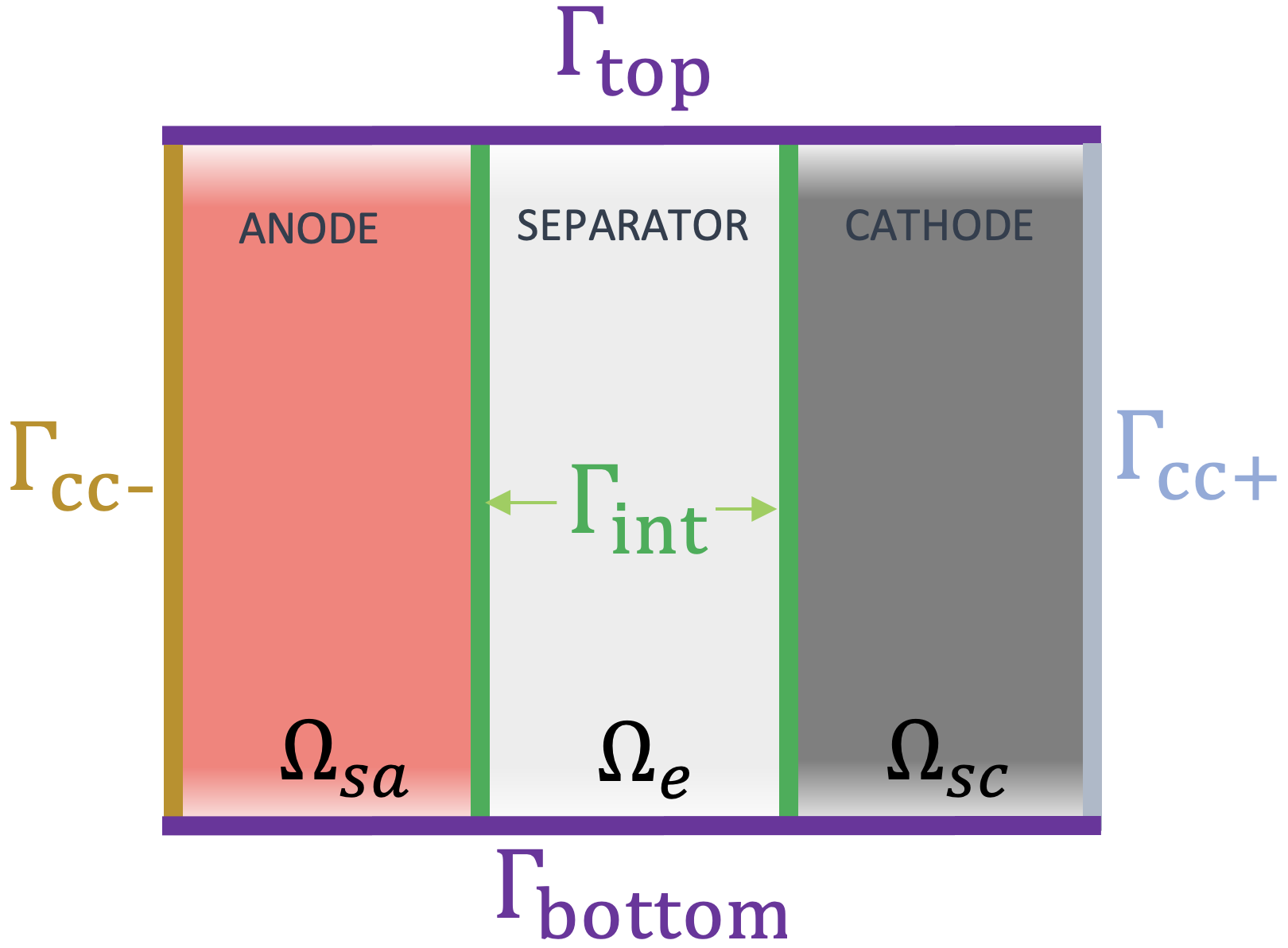}
    \caption{Battery diagram showing boundary portions that are referenced in the specification of boundary conditions.}
    \label{fig:boundaries-diagram}
\end{figure}

We have assumed that we are studying a small representative part of the battery microstructure. For that reason, it is possible to prescribe symmetry or periodicity boundary conditions on some portions of $\Gamma$. We opt for the former applied to the top and bottom faces of our domain. Figure \ref{fig:boundaries-diagram} shows the parts of the boundary where the BC and the interface conditions are to be imposed.
\subsubsection{Current density and electric potential}
We begin with the electrodes. Their boundary is split as follows
\begin{itemize}
    \item $\Gneg$ denotes a portion of the {\bf anode} boundary that is in contact with the negative current collector, with a homogeneous Dirichlet BC, that is, $\phi_{s}=0$.
    \item $(\Gamma_{s}\cap \Gtop) \cup (\Gamma_{s}\cap \Gbot)$ corresponds to no inward or outward electrical current, i.e., $\bsi_s\cdot\nml=0$;
    \item $\Gpos$ lies on a part of the {\bf cathode} boundary that is in contact with the positive current collector, where a known current density is applied. The applied current density is $\iapp$, which is positive if the battery is being discharged or negative if it is being charged. We have then $\bsi_s\cdot\nml=\iapp$.
    \item $(\Gamma_{s}\cap \interface)$, as explained above, is given an interface condition based on the Butler-Volmer kinetics model: $\underbrace{\bsi_s\cdot\nml}_{=\bsi_s\cdot\nmlint}=\ibv$, computed through \eqref{eq:butler-volmer}.
\end{itemize}
It is the turn of the electrolytic medium, $\iota=e$. We consider only two subsets of $\Gamma_e$:
\begin{itemize}
    \item $(\Gamma_{e}\cap \Gtop) \cup (\Gamma_{e}\cap \Gbot)$, whose condition is flux-free (or symmetry), i.e., $\bsi_e\cdot\nml=0$;
    \item $\Gamma_{e}\cap \interface$, where, in this case, the BC imposed is 
    $\underbrace{\bsi_e\cdot\nml}_{=-\bsi_e\cdot\nmlint}=-\ibv$.
\end{itemize}

\subsubsection{Species flux and concentration}
At medium $\iota=s$, we have just two types of BC on $\Gamma_s$:
\begin{itemize}
    \item Flux-free BC, $\bsj_s\cdot\nml=0$, on $\Gneg \cup \Gpos$ (because of vanishing mass flux across the current collectors) and  on $(\Gamma_{s}\cap \Gtop) \cup (\Gamma_{s}\cap \Gbot)$ (e.g., due to symmetry);
    \item Ion flux across $(\Gamma_{s}\cap \interface)$, accounted by 
$\underbrace{\bsj_s\cdot\nml}_{=\bsj_s\cdot\nmlint}=\ibv/F$.
\end{itemize}
The electrolyte, in turn, has BC that are analogous to those for the electric potential:
\begin{itemize}
    \item $(\Gamma_{s}\cap \Gtop) \cup (\Gamma_{s}\cap \Gbot)$ with $\bsj_e\cdot\nml=0$;
    \item $(\Gamma_{e}\cap \interface)$, with 
    $\underbrace{\bsj_e\cdot\nml}_{=-\bsj_e\cdot\nmlint}=-(1-t_+)\ibv/F$.
\end{itemize}
In both of these sets of BC, we have seen that the reaction current density, $\ibv$, is divided by $F$ in order to get the respective species flux. Moreover, the value for the electrolyte is altered by factor $1-t_+$, which relates to the fact that not all the electric current is carried across the interface by the Li+ ions, but also by accompanying anions, hence their part must be subtracted.

\subsubsection{Heat flux and temperature}
For the sake of the present work, this thermal system is assumed to be adiabatic, thereby $\bsq\cdot\nml=0$ on $\Gamma$, which can be regarded as a homogeneous Neumann BC. However, the interface condition follows a more elaborate rule. Since on $\interface$ occurs the electrochemical reaction, the net interfacial Ohmic heat flux is defined by Ohm's law, yielding the interface condition
\begin{equation}
    (\bsq_s-\bsq_e)\cdot\nmlint=\eta\ibv\ \text{ on }\interface.
    \label{eq:interface_ohmic}
\end{equation}
\subsubsection{Stress and displacement}
Recall that these variables are restricted to medium $s$. We have two types of conditions:
\begin{itemize}
    \item on boundary portion $(\Gamma_{s}\cap \Gtop) \cup (\Gamma_s\cap\interface)$ there is a traction-free BC, i.e., $\bsigma\,\nml=\bzero$;
    \item on boundary portions $\Gamma_{s}\cap \Gbot$ and $\Gneg \cup \Gpos$ there holds a \emph{symmetry} or a \emph{slip} BC, respectively, which have a different interpretation in mechanics but both coincide mathematically: $\bsu\cdot\nml=0$ and $(\bsigma\,\nml)\times\nml=\bzero$ (i.e., zero normal displacement \& zero shear force).
\end{itemize}

\subsection{Initial conditions}
\label{subsec:model_ic}
Regarding the momentum balance and the energy balance, the relevant initial conditions are the strain-free state and the thermal equilibrium state, that is
\begin{itemize}
    \item $\theta(0,\bfx)=\theta_0$ for every $\bfx\in\overline{\Omega}$;
    \item $\bsu(0,\bfx)=\bzero$ for every $\bfx\in\overline{\Omega_s}$.
\end{itemize}

On the other hand, it is a more delicate matter to establish appropriate initial conditions for the electrochemical variables. Moreover, the PDEs for electric potential are quasi-static, so that an actual initial condition should not be necessary for them to be solved, but an arbitrary potential at $t=0$ could induce fluxes in an unloaded state, that is, a lack of equilibrium to start with. As explained in subsection \ref{subsec:butler-volmer}, the electrochemical equilibrium occurs if the potential difference at the interface equals the open-circuit potential of the electrode. The following choice of values allows for an initial state of electrochemical equilibrium, as the difference of voltage between each electrode and the eletrolyte equals the corresponding open-circuit potential.
\begin{itemize}
    \item $c_{sa}(0,\bfx)= c_{sa,0}$ for every $\bfx\in\overline{\Omega_{sa}}$;
    \item $c_{sc}(0,\bfx)=c_{sc,0}$ for every $\bfx\in\overline{\Omega_{sc}}$;
    \item $c_{e} (0,\bfx)=c_{e,0} $ for every $\bfx\in\overline{\Omega_{e}}$;
    \item $\phi_{sa}(0,\bfx)=0$ for every $\bfx\in\overline{\Omega_{sa}}$;
    \item $\phi_{sc}(0,\bfx)=\ocpsc(c_{sc,0}/\cscmax)-\ocpsa(c_{sa,0}/\csamax)$ for every $\bfx\in\overline{\Omega_{sc}}$.
    \item $\phi_{e} (0,\bfx)=-\ocpsa(c_{sa,0}/\csamax)$ for every $\bfx\in\overline{\Omega_{e}}$;
\end{itemize}
\subsection{Weak problem}
We have presented the strong form of the initial boundary value problem that models our Lithium-ion battery. The method chosen to numerically solve the model requires the derivation of a weak form of the problem. The relaxation procedure that delivers the weak form is complemented with the specification of suitable functional spaces for the involved variables.
\label{subsec:weak}
\subsubsection{Relaxation}
The classical weak formulation for our system of PDEs is obtained by the usual relaxation procedure. Let us begin with the energy balance equation. Let $\zeta$ be a test function, whose functional space is still unspecified. We multiply \eqref{eq:energy-balance} by $\zeta$ and integrate throughout $\Omega$:
\begin{equation*}
    \int\limits_{\Omega} \zeta\; \left[\rho_\iota C_{v,\iota} \dot{\theta} - \nabla\cdot\left(\lambda_\iota\nabla\theta\right)\right] dx
    =
    \int\limits_{\Omega} \zeta\; \bsi_\iota \cdot \nabla \phi_\iota \,dx
\end{equation*}
Now, we split the domain into the solid and electrolyte regions.
\begin{equation*}
    \int\limits_{\Omega} \zeta\; \rho_\iota C_{v,\iota} \dot{\theta} \,dx 
    - \int\limits_{\Omega_s} \zeta\; \nabla\cdot\left(\lambda_s\nabla\theta\right) \,dx
    - \int\limits_{\Omega_e} \zeta\; \nabla\cdot\left(\lambda_e\nabla\theta\right) \,dx
    =
    \int\limits_{\Omega} \zeta\; \bsi_\iota \cdot \nabla \phi_\iota \,dx
\end{equation*}
Integrating by parts and recalling the definition of heat flux, \eqref{eq:Fourierlaw},
\begin{equation*}
    \int\limits_{\Omega} \zeta\; \rho_\iota C_{v,\iota} \dot{\theta} \,dx 
    + \int\limits_{\Omega} \nabla\zeta \cdot \lambda_\iota\nabla\theta \,dx
    + \int\limits_{\bdry\Omega_s} \zeta \bsq_s \cdot \nml \,dA 
    + \int\limits_{\bdry\Omega_e} \zeta\; \bsq_e \cdot \nml \,dA 
    =
    \int\limits_{\Omega} \zeta \bsi_\iota \cdot \nabla \phi_\iota \,dx
\end{equation*}

\begin{equation*}
    \int\limits_{\Omega} \zeta\; \rho_\iota C_{v,\iota} \dot{\theta} \,dx 
    + \int\limits_{\Omega} \nabla\zeta\; \cdot \lambda_\iota\nabla\theta \,dx
    + \int\limits_{\Gamma} \zeta\; \bsq_\iota \cdot \nml \,dA 
    + \int\limits_{\interface} \zeta\; (\bsq_s-\bsq_e) \cdot \nmlint  \,dA 
    =
    \int\limits_{\Omega} \zeta\; \bsi_\iota \cdot \nabla \phi_\iota \,dx
\end{equation*}
Since the boundary conditions for the energy equation are of homogeneous Neumann type only, the third term above vanishes. Finally, by 
bringing over the interface Ohmic heat expression \eqref{eq:interface_ohmic}, we get
\begin{equation}
    \int\limits_{\Omega} \zeta\; \rho_\iota C_{v,\iota} \dot{\theta} \,dx 
    + \int\limits_{\Omega} \nabla\zeta\; \cdot \lambda_\iota\nabla\theta \,dx
    =
    \int\limits_{\Omega} \zeta\; \bsi_\iota \cdot \nabla \phi_\iota \,dx
    - \int\limits_{\interface} \zeta\;\eta \ibv \,dA 
    \label{eq:weak_energy}
\end{equation}
The weak form of the concentration PDEs follows a similar path. Let $r_\iota$ be a test function defined over $\Omega_\iota$, multiply \eqref{eq:solid-conc-pde} by $r_s$ and \eqref{eq:elyte-conc-pde} by $r_e$, and perform integration by parts:
\begin{align*}
    \int\limits_{\Omega_s}r_s\;\dot{c}_s \,dx + \int\limits_{\Omega_s}\nabla r_s\cdot (D_s(c_s,\pi)\nabla c_s)\,dx 
    &= -\int\limits_{\Gamma_s} r_s\;\bsj_s\cdot\nml \,dA,\\
    \int\limits_{\Omega_e}r_e\;\dot{c}_e \,dx + \int\limits_{\Omega_e}\nabla r_e\cdot (D_e\nabla c_e)\,dx 
    &= -\int\limits_{\Gamma_e} r_e\;\bsj_e\cdot\nml \,dA.
\end{align*}
The integral over each component boundary reduces to the interface part only, as the normal ion flux vanishes on the remaining boundary parts per the prescribed BC.
Using the interface conditions, the weak equations for $c_s$ and $c_e$ are
\begin{align}
    \int\limits_{\Omega_s}r_s\;\dot{c}_s \,dx + \int\limits_{\Omega_s}\nabla r_s\cdot (D_s(c_s,\pi)\nabla c_s)\,dx 
    &= -\frac{1}{F}\int\limits_{\Gamma_s\cap\interface} r_s\;\ibv \,dA,
    \label{eq:weak_concs}\\
    \int\limits_{\Omega_e}r_e\;\dot{c}_e \,dx + \int\limits_{\Omega_e}\nabla r_e\cdot (D_e\nabla c_e)\,dx 
    &= \frac{1-t_+}{F}\int\limits_{\Gamma_e\cap\interface} r_e\;\ibv \,dA.
    \label{eq:weak_conce}
\end{align}
Now, let us handle the relaxation of the electric potential equations. According to component $\iota$, let $\psi_\iota$ be a test function over $\Omega_\iota$, so that we can proceed with the relaxation of  \eqref{eq:solid-potential-pde} and 
\eqref{eq:elyte-potential-pde}.  
Upon breaking down the resulting surface integrals $\int_{\Gamma_s}$ and $\int_{\Gamma_e}$, we have to account for the BC and interface conditions. 
Because either the test functions or the normal current density vanish on $\Gtop$, $\Gbot$, $\Gneg$ and $\Gpos$ (only if $\iota=e$), 
there just remain the following terms
\begin{align}
\int\limits_{\Omega_s} \nabla\psi_s \cdot (\gamma_s \nabla \phi_s) \,dx
&=
-\int\limits_{\Gpos} \psi_s \iapp \,dA
-\int\limits_{\Gamma_s\cap\interface} \psi_s \, \underbrace{\ibv}_{=\bsi_s \cdot \nml} \,dA
\,, \label{eq:elpls_pde_relaxed} \\
\int\limits_{\Omega_e} \nabla\psi_e \cdot (
\kappa_e\nabla\phi_e
) \,dx
&=
-\int\limits_{\Omega_e} \nabla\psi_e \cdot (
\kappa_D\nabla\ln c_e
) \,dx
-\int\limits_{\Gamma_e\cap\interface} \psi_e\, \underbrace{(-\ibv)}_{=\bsi_e \cdot \nml} \,dA
\,. \label{eq:elple_pde_relaxed}
\end{align}
The treatment of the normal current density on $\interface$ is different from the approach given to the the previous equations. 
The idea here is to leverage on a simple third-order approximation to $\ibv$, namely,
$$
\ibv
= 2 I_c \sinh\left(\frac{F}{2 R \theta}\eta\right)
=
2 I_c \left\{\frac{F}{ 2 R \theta}\eta + \mcO[(\tfrac{F}{ 2 R \theta}\eta)^3]\right\}
\approx \frac{I_c F}{R \theta} 
\underbrace{(\phi_s - \phi_e - \ocp)}_{\eta}.
$$
By introducing this approximation into \eqref{eq:elpls_pde_relaxed} and \eqref{eq:elple_pde_relaxed}, and after reorganizing terms, we finally arrive at
\begin{align}
  \int\limits_{\Omega_s} \nabla\psi_s \cdot (\gamma_s \nabla \phi_s) \,dx
+\int\limits_{\Gamma_s\cap\interface} \psi_s \tfrac{I_c F}{R \theta} \phi_s \,dA
&=
-\int\limits_{\Gpos} \psi_s \iapp \,dA
+\int\limits_{\Gamma_s\cap\interface} \psi_s \tfrac{I_c F}{R \theta} (\phi_e + \ocp)dA
\,, \label{eq:weak_elpls} \\
\int\limits_{\Omega_e} \nabla\psi_e \cdot (
\kappa_e\nabla\phi_e
) \,dx
+\int\limits_{\Gamma_e\cap\interface} \psi_e 
\tfrac{I_c F}{R \theta} \phi_e
\,dA
&=
-\int\limits_{\Omega_e} \nabla\psi_e \cdot (
\kappa_D\nabla\ln c_e
) \,dx
+\int\limits_{\Gamma_e\cap\interface} \psi_e
\tfrac{I_c F}{R \theta} (\phi_s - \ocp)
\,dA 
\,. \label{eq:weak_elple}
\end{align}

The last PDE we are required to weaken is the elasticity equation \eqref{eq:displacement-pde}. In this case, the test function $\bsv$ is a vector field defined on $\Omega_s$ only. 
It must satisfy homogeneous Dirichlet conditions in the normal direction of $\Gpos$, $\Gneg$, $\Gamma_s\cap\Gtop$ and $\Gamma_s\cap\Gbot$. 
Moreover, since the traction is zero on the remaining boundary part, $\interface$, there is no boundary term in the weak form of this equation. We are thus left with
\begin{equation}
    \int\limits_{\Omega_s} \nabla_{_{\Sym}}\bsv\,:\,\bsfC:\nabla_{_{\Sym}}\bsu\,dx
    =
    \int\limits_{\Omega_s} \nabla\cdot\bsv\left(3 K \alpha (\theta-\theta_0) + 3 K (c_s-c_{s,0}) \right)  \,dx \,,
    \label{eq:weak_displ}
\end{equation}
where we have used the fact that $\bsfC:\nabla_{_{\Sym}}\bsu$ is a symmetric tensor field, allowing us to work with the symmetric gradient of $\bsv$ instead of the full gradient.
\begin{remark}
    Notice that in the relaxed expressions for the concentration and the electric potential, that is, in \eqref{eq:weak_concs}-\eqref{eq:weak_elple}, we have included terms that are integrated over $\Gamma_s\cap\interface$ and $\Gamma_e\cap\interface$. Those integrals, contrary to integrals over just $\interface$, preserve the orientation of $\Gamma_s$ and $\Gamma_e$, respectively, making it clearer for implementation at the element level.
\end{remark}
\subsubsection{Classical weak formulation}
Now that every PDE has been relaxed along with its corresponding BC and interface conditions, we can gather the full system again. 
Prior to the actual statement, it is imperative to define appropriate functional spaces for our trial and test functions. 

If $\mcO$ is a non-empty open set in $\R^3$, recall the $L^2(\mcO)$ space:
$$
L^2(\mcO)=\left\{ f:\mcO\to\R, \text{ a Lebesgue-measurable function such that } \int\limits_\mcO |f|^2\,dx<+\infty \right\},
$$
which is a Hilbert space with inner product and norm defined by:
$$
(f,g)_\mcO=\int\limits_\mcO fg\,dx,
\qquad
\| f \|_\mcO=\sqrt{(f,f)_\mcO}\, .
$$
Now, let $H^1(\mcO)$ denote the following space
$$
H^1(\mcO)=\left\{ f\in L^2(\mcO) \text{ such that } \nabla f \in L^2(\mcO;\R^3) \right\},
$$
where $L^2(\mcO;\R^3)$ is the space of vector-valued $L^2$ functions (each component by itself belongs to $L^2(\mcO)$). Space $H^1(\mcO)$ is a Sobolev space with norm $\| f \|_{H^1(\mcO)}=\sqrt{\| f \|_\mcO^2 + \| \nabla f\|_\mcO^2}$.
Consider then
\begin{align*}
    \scV &= \scV_{\theta}\times\scV_{c_s}\times\scV_{c_e}\times\scV_{\phi_s}\times\scV_{\phi_e}\times\scV_{\bsu}\;,\quad\text{where}\\
    \scV_{\theta} &=H^1(\Omega),\\
    \scV_{c_s} &=H^1(\Omega_s),\\
    \scV_{c_e} &=H^1(\Omega_e),\\
    \scV_{\phi_s} &=H^1_{\Gneg}(\Omega_s),\\
    \scV_{\phi_e} &=H^1(\Omega_e),\\
    \scV_{\bsu} &= H^1_{\Gamma_{s}\setminus\interface}(\Omega_s;\R^3),
\end{align*}
where the subscript of $H^1_{\Gneg}$ indicates the subspace of $H^1$ whose functions have a vanishing trace on the boundary subset $\Gneg$; similarly, $H^1_{\Gamma_{s}\setminus\interface}(\Omega_s;\R^3))$ is the subspace of vector-valued $H^1$ fields with a vanishing \emph{normal component} on $\Gamma_{s}\setminus\interface$. Product space $\scV$ is equipped with a usual product norm $\|\cdot\|_\scV$.

As $\scV$ just considers the behavior in space of our functions, time dependence must be incorporated on top of it. In that regard, we must account for functions
$\theta(t,\cdot)$, $c_s(t,\cdot)$, $c_e(t,\cdot)$, $\phi_s(t,\cdot)$, $\phi_e(t,\cdot)$, $\bsu(t,\cdot)$, that belong to the appropriate space at each $t\in[0,\tend]$. Following the notation of \cite{evans2010partial}, we must deal with a space of the type
$$
L^2(0,\tend;\scV)=\left\{ \fkv:[0,\tend]\to\scV\text{ strongly measurable, such that }\int\limits_0^{\tend} \| \fkv(t,\cdot) \|_\scV^2\,dt<+\infty \right\}.
$$

The weak problem reads:

Find $(\theta,c_s,c_e,\phi_s,\phi_e,\bsu)\in L^2(0,\tend;\scV)$ such that for every $(\zeta,r_s,r_e,\psi_s,\phi_e,\bsv)\in\scV$:
\begin{equation}
\left\{
    \begin{array}{rl}
        \left( \zeta ,\rho_\iota C_{v,\iota} \dot{\theta} \right)_{\Omega}
        + \left( \nabla\zeta , \lambda_\iota\nabla\theta \right)_{\Omega}               
            =& \left( \zeta , \bsi_\iota \cdot \nabla \phi_\iota\right)_{\Omega}
              - \langle \zeta , \eta\ibv \rangle_{\interface}
            \vspace{0.2cm}\\
        \left( r_s , \dot{c_s} \right)_{\Omega_s}
        + \left( \nabla r_s , D_s(c_s,\pi)\nabla c_s \right)_{\Omega_s}
            =& -\frac{1}{F}\langle r_s , \ibv \rangle_{\Gamma_s\cap\interface}
            \vspace{0.2cm}\\            
        \left(r_e , \dot{c_e} \right)_{\Omega_e}
        + \left( \nabla r_e , D_e\nabla c_e \right)_{\Omega_e}
            =& \frac{1-t_+}{F}\langle r_e , \ibv \rangle_{\Gamma_e\cap\interface}
            \vspace{0.2cm}\\
        \left( \nabla \psi_s , \gamma_s\nabla \phi_s \right)_{\Omega_s}
          + \langle \psi_s , \tfrac{I_c F}{R \theta}\phi_s \rangle_{\Gamma_s\cap\interface}
            =& -\langle \psi_s , \iapp \rangle_{\Gpos}
            \langle \psi_s , \tfrac{I_c F}{R \theta}(\phi_e+\ocp) \rangle_{\Gamma_s\cap\interface}
            \vspace{0.2cm}\\
        \left( \nabla \psi_e , \kappa_e\nabla\phi_e \right)_{\Omega_e}
          + \langle \psi_e , \tfrac{I_c F}{R \theta}\phi_e \rangle_{\Gamma_e\cap\interface}
            =& -\left( \nabla \psi_e , \kappa_D\nabla \ln c_e \right)_{\Omega_e}
            +\langle \psi_s , \tfrac{I_c F}{R \theta}(\phi_s-\ocp) \rangle_{\Gamma_s\cap\interface}
            \vspace{0.2cm}\\ 
        \left( \nabla_{_{\Sym}} \bsv , \bsfC : \nabla_{_{\Sym}} \bsu \right)_{\Omega_s}
            =& \left( \nabla\cdot \bsv, 3 K \alpha (\theta-\theta_0) + 3 K \omega (c_s-c_{s,0})  \right)_{\Omega_s}
            \vspace{0.2cm}
    \end{array}
\right.    
\label{eq:weak_system}
\end{equation}
subject to the initial condition
\begin{equation}
(\theta(0,\cdot),c_s(0,\cdot),c_e(0,\cdot),\phi_s(0,\cdot),\phi_e(0,\cdot),\bsu(0,\cdot))=
(\theta_0,c_{s,0},c_{e,0},\phi_{s,0},\phi_{e,0},\bsu_0)\in\scV,
    \label{eq:weak_system_init}
\end{equation}
where the suitable values of the 0 subscript quantities can be found in subsection \ref{subsec:model_ic}. The angled brackets in \eqref{eq:weak_system}, say $\langle\cdot,\cdot\rangle_{\mcF}$ for some boundary subset $\mcF$,  represent the extension of the $L^2(\mcF)$ inner product to a duality pairing for $H^{1/2}(\mcF)\times H^{-1/2}(\mcF)$, which are the appropriate spaces for the traces therein appearing.

\section{Discretization approach} 
\label{sec:method}


\subsection{Space discretization}
The space discretization of problem \eqref{eq:weak_system}, \eqref{eq:weak_system_init},
is carried out by means of the Galerkin method. So, let $\scV^{hp}\subset\scV$ be a finite-dimensional subspace defined over a shape-regular finite element mesh of $\Omega$, denoted $\mesh^{hp}$. Suppose that $\scV^{hp}$ is also a product of discrete spaces for each variable, such as $\scV^{hp}_{\theta}\subset\scV_\theta$, and so on. The basis functions  of each factor of $\scV^{hp}$ (say, $\{w_{\theta,i}\}_{i=1}^{\dim \scV^{hp}_{\theta}}$, etc.) are globally continuous (over all $\Omega$ or just the corresponding subdomain), piece-wise polynomials with variable degree $p(E)\geq1$ throughout the mesh elements $E\in\mesh^{hp}$. The mesh and the finite-element space are assumed to remain fixed along the entire time interval $[0,\tend]$.
The discrete solutions have the form
\begin{equation}
\begin{aligned}
    \theta^{hp}=\theta^{hp}(t,\bfx) &=
    \sum_{i=1}^{\dim \scV^{hp}_{\theta}}d_{\theta,i}(t) w_{\theta,i}(\bfx) \\
    c_s^{hp}=c_s^{hp}(t,\bfx) &=
    \sum_{i=1}^{\dim \scV^{hp}_{c_s}}d_{c_s,i}(t) w_{c_s,i}(\bfx) \\
    c_e^{hp}=c_e^{hp}(t,\bfx) &=
    \sum_{i=1}^{\dim \scV^{hp}_{c_e}}d_{c_e,i}(t) w_{c_e,i}(\bfx) \\
    \phi_s^{hp}=\phi_s^{hp}(t,\bfx) &=
    \sum_{i=1}^{\dim \scV^{hp}_{\phi_s}}s_{\phi_s,i}(t) w_{\phi_s,i}(\bfx) \\
    \phi_e^{hp}=\phi_e^{hp}(t,\bfx) &=
    \sum_{i=1}^{\dim \scV^{hp}_{\phi_e}}s_{\phi_e,i}(t) w_{\phi_e,i}(\bfx) \\
    \bsu^{hp}=\bsu^{hp}(t,\bfx) &=
    \sum_{i=1}^{\dim \scV^{hp}_{\bsu}}s_{\bsu,i}(t) \bsw_{\bsu,i}(\bfx)
\end{aligned}
\label{eq:discrete_solutions}
\end{equation}
Notice that we have used a $d$ for the coefficients in the first three unknowns, $\theta^{hp}$, $c_s^{hp}$ and $c_e^{hp}$, to relate to the dynamic nature of their governing equations, while we wrote an $s$ for the coefficients of $\phi_s^{hp}$, $\phi_e^{hp}$ and $\bsu^{hp}$, since these are mainly modeled by quasi-static equations. This choice of symbols is helpful below when dealing with the discretization in time.

Let $\Pi^{hp}:\scV^{\text{reg}}\to\scV^{hp}$ be the linear interpolation operator of $\scV^{hp}$, well-defined on a sufficiently regular subspace $\scV^{\text{reg}}\subset\scV$. In order to use $\Pi^{hp}$, we assume that the initial condition belongs to $\scV^{\text{reg}}$, and so do the exact solutions  $\theta(t,\cdot)$, $c_s(t,\cdot)$, $c_e(t,\cdot)$, $\phi_s(t,\cdot)$, $\phi_e(t,\cdot)$, $\bsu(t,\cdot)$, at each time $t\in(0,\tend]$.

The space-discrete problem reads:

Find $(\theta^{hp},c_s^{hp},c_e^{hp},\phi_s^{hp},\phi_e^{hp},\bsu^{hp})\in L^2(0,\tend;\scV^{hp})$ such that for every $(\zeta^{hp},r_s^{hp},r_e^{hp},\psi_s^{hp},\phi_e^{hp},\bsv^{hp})\in\scV^{hp}$:
\begin{equation}
\left\{
    \begin{array}{rl}
        \left( \zeta^{hp} ,\rho_\iota C_{v,\iota} \dot{\theta}^{hp} \right)_{\Omega}
        + \left( \nabla\zeta^{hp} , \lambda_\iota\nabla\theta^{hp} \right)_{\Omega}               
            =& \left( \zeta^{hp} , \bsi_\iota^{hp} \cdot \nabla \phi_\iota^{hp}\right)_{\Omega}
              - \langle \zeta^{hp} , \eta^{hp}\ibv^{hp} \rangle_{\interface}
            \vspace{0.2cm}\\
        \left( r_s^{hp} , \dot{c}_s^{hp} \right)_{\Omega_s}
        + \left( \nabla r_s^{hp} , D_s(c_s^{hp},\pi^{hp})\nabla c_s^{hp} \right)_{\Omega_s}
            =& -\frac{1}{F}\langle r_s^{hp} , \ibv^{hp} \rangle_{\Gamma_s\cap\interface}
            \vspace{0.2cm}\\            
        \left(r_e^{hp} , \dot{c}_e^{hp} \right)_{\Omega_e}
        + \left( \nabla r_e^{hp} , D_e\nabla c_e^{hp} \right)_{\Omega_e}
            =& \frac{1-t_+}{F}\langle r_e^{hp} , \ibv^{hp} \rangle_{\Gamma_e\cap\interface}
            \vspace{0.2cm}\\
        \left( \nabla \psi_s^{hp} , \gamma_s\nabla \phi_s^{hp} \right)_{\Omega_s}
          + \langle \psi_s^{hp} , \tfrac{I_c^{hp} F}{R \theta^{hp}}\phi_s^{hp} \rangle_{\Gamma_s\cap\interface}
            =& -\langle \psi_s^{hp} , \iapp \rangle_{\Gpos}
            \langle \psi_s^{hp} , \tfrac{I_c^{hp} F}{R \theta^{hp}}(\phi_e^{hp}+\ocp^{hp}) \rangle_{\Gamma_s\cap\interface}
            \vspace{0.2cm}\\
        \left( \nabla \psi_e^{hp} , \kappa_e\nabla\phi_e^{hp} \right)_{\Omega_e}
          + \langle \psi_e^{hp} , \tfrac{I_c^{hp} F}{R \theta^{hp}}\phi_e^{hp} \rangle_{\Gamma_e\cap\interface}
            =& -\left( \nabla \psi_e^{hp} , \kappa_D\nabla \ln c_e^{hp} \right)_{\Omega_e}
            \\ &+\langle \psi_s^{hp} , \tfrac{I_c^{hp} F}{R \theta^{hp}}(\phi_s^{hp}-\ocp^{hp}) \rangle_{\Gamma_s\cap\interface}
            \vspace{0.2cm}\\ 
        \left( \nabla_{_{\Sym}} \bsv^{hp} , \bsfC : \nabla_{_{\Sym}} \bsu^{hp} \right)_{\Omega_s}
            =& \left( \nabla\cdot \bsv^{hp}, 3 K \alpha (\theta^{hp}-\theta_0) + 3 K \omega (c_s^{hp}-c_{s,0})  \right)_{\Omega_s}
            \vspace{0.2cm}
    \end{array}
\right.    
\label{eq:space_discrete_system}
\end{equation}
subject to the initial condition
\begin{equation}
(\theta^{hp}(0,\cdot),c_s^{hp}(0,\cdot),c_e^{hp}(0,\cdot),\phi_s^{hp}(0,\cdot),\phi_e^{hp}(0,\cdot),\bsu^{hp}(0,\cdot))=
\Pi^{hp}(\theta_0,c_{s,0},c_{e,0},\phi_{s,0},\phi_{e,0},\bsu_0)\in\scV^{hp},
    \label{eq:space_discrete_init}
\end{equation}
In \eqref{eq:space_discrete_system}, the $hp$ superscript on derived quantities (fluxes, current density, hydrostatic pressure) imply their computation through the $hp$ counterparts of the solution components.

\begin{remark}
    Since the current system models a representative unit of the microstructure of a battery, at such a scale, the initial conditions can be taken as constant values (at least, piecewise), so that they should belong to every conforming $\scV^{hp}$ conceivable in practice, hence the interpolation under $\Pi^{hp}$ may produce no effect on the initial conditions \eqref{eq:space_discrete_init}.
\end{remark}

\subsection{Discrete initial-value problem}

Upon testing the Galerkin system \eqref{eq:space_discrete_system} with every basis function of $\scV^{hp}$, in addition to substituting the approximate unknowns with \eqref{eq:discrete_solutions}, the problem can be written in short as
\begin{equation}
    \begin{array}{rl}
\bfM\dot{\bfd}(t)  + \bfK_{\bfd}\bfd(t)=& \bfb_{\bfd}(t,\bfd(t),\bfs(t)) \\
\bfK_{\bfs}\bfs(t)  =& \bfb_{\bfs}(t,\bfd(t),\bfs(t))
\end{array}
\label{eq:short_semidiscrete}
\end{equation}
where we have taken advantage of the bilinearity of the LHS pairings to extract the coefficient vectors, 
$$
\bfd(t)=\begin{pmatrix}
    \bfd_{\theta}(t)\\\bfd_{c_s}(t)\\\bfd_{c_e}(t)
\end{pmatrix},
\qquad
\bfs(t)=\begin{pmatrix}
    \bfs_{\phi_s}(t)\\\bfs_{\phi_e}(t)\\\bfs_{\bsu}(t)
\end{pmatrix}.
$$
The lengths of $\bfd(t)$ and $\bfs(t)$ are
$$
m_{\bfd}=\dim(\scV_\theta^{hp}\times\scV_{c_s}^{hp}\times\scV_{c_e}^{hp}) \quad 
\text{and}\quad 
m_{\bfs}=\dim(\scV_{\phi_s}^{hp}\times\scV_{\phi_e}^{hp}\times\scV_{\bsu}^{hp}) ,
$$
respectively. For the sake of completeness, we have
$$
\bfM=\begin{pmatrix}
    \bfM_{\theta} & \bzero & \bzero \\
    \bzero & \bfM_{c_s} & \bzero \\
    \bzero & \bzero & \bfM_{c_e} \\
\end{pmatrix}\,,
$$
with its entries given by
$$
[\bfM_{\theta}]_{ij}=(w_{\theta,i},\rho C_{v,\iota} w_{\theta,j})_{\Omega},
\quad 
[\bfM_{c_s}]_{ij}=(w_{c_s,i}, w_{c_s,j})_{\Omega_s},
\quad 
[\bfM_{c_e}]_{ij}=(w_{c_e,i}, w_{c_e,j})_{\Omega_e};
$$
furthermore,
$$
\bfK_{\bfd}=\begin{pmatrix}
    \bfK_{\theta} & \bzero & \bzero \\
    \bzero & \bfK_{c_s} & \bzero \\
    \bzero & \bzero & \bfK_{c_e} \\
\end{pmatrix}\,,
$$
where the submatrices are computed through
$$
[\bfK_{\theta}]_{ij}=
(\nabla w_{\theta,i},\lambda_\iota \nabla w_{\theta,j})_{\Omega},
\quad 
[\bfK_{c_s}]_{ij}=
(\nabla w_{c_s,i}, D_s(c_s^{hp},\pi^{hp}) \nabla w_{c_s,j})_{\Omega_s},
\quad 
[\bfK_{c_e}]_{ij}=
(\nabla w_{c_e,i}, D_e \nabla w_{c_e,j})_{\Omega_e};
$$
and,
$$
\bfK_{\bfs}=\begin{pmatrix}
    \bfK_{\phi_s} & \bzero & \bzero \\
    \bzero & \bfK_{\phi_e} & \bzero \\
    \bzero & \bzero & \bfK_{\bsu} \\
\end{pmatrix}
$$
with
$$
[\bfK_{\phi_s}]_{ij}=
(\nabla w_{\phi_s,i},\gamma_s \nabla w_{\phi_s,j})_{\Omega_s} 
+ \langle w_{\phi_s,i} , \tfrac{I_c^{hp} F}{R \theta^{hp}}  w_{\phi_s,j} \rangle_{\Gamma_s\cap\interface},
$$
$$
[\bfK_{\phi_e}]_{ij}=
(\nabla w_{c_s,i}, \kappa_e \nabla w_{c_s,j})_{\Omega_s}
+ \langle w_{\phi_e,i} , \tfrac{I_c^{hp} F}{R \theta^{hp}} w_{\phi_e,j} \rangle_{\Gamma_e\cap\interface},
\quad
[\bfK_{\bsu}]_{ij}=
(\nabla_{_{\Sym}} \bsw_{\bsu,i}, \bsfC : \nabla_{_{\Sym}} \bsw_{\bsu,j})_{\Omega_s}\;.
$$

As usual, we name $\bfM$ the mass matrix, while $\bfK_{\bfd}$  and $\bfK_{\bfs}$ are known as stiffness matrices. 
These are all real, symmetric, positive-definite (SPD) matrices. With the present assumptions on material properties, $\bfM$ and $\bfK_{\bfs}$ are constant, unlike
$\bfK_{\bfd}$, which is variable because of the diffusivity of the solid electrodes.
As it can be seen, there are no off-diagonal blocks in $\bfM$ and $\bfK$ (a property that provides efficiency during assembly, and sparsity to the linear systems to be solved), thanks to how both the strong and weak systems have been organized throughout this document: any linear coupling of two $\bfd$ variables or two $\bfs$ variables, which would induce an off-diagonal block in $\bfK_{\bfd}$ or $\bfK_{\bfs}$, respectively, has been moved to the RHS of their equations. 

We do not detail the load operators in \eqref{eq:short_semidiscrete}, so we leave to the reader to fill the gaps by appropriately bringing the RHS terms of \eqref{eq:space_discrete_system} into $\bfb_{\bfd}$ or $\bfb_{\bfs}$. 
Now, let $\bff_{\bfd}$ and $\bff_{\bfs}$ be defined by
\begin{equation}
    \bff_{\bfd}(t,\bfd,\bfs)=\bfM^{-1}\left( -\bfK_{\bfd} \bfd + \bfb_{\bfd} (t,\bfd,\bfs) \right),
    \qquad 
    \bff_{\bfs}(t,\bfd,\bfs)=\bfK_{\bfs}^{-1} \bfb_{\bfs}(t,\bfd,\bfs).
\end{equation}
We note that inside $\bff_{\bfd}$ and $\bff_{\bfs}$ are all the couplings, both linear and nonlinear, that were considered in the original system of equations, and any extra coupling could be added without changing the structure.

The Galerkin formulation \eqref{eq:discrete_solutions}-\eqref{eq:space_discrete_init} in the form of a discrete initial-value problem is as follows:

Find $\bfd:[0,\tend]\to \R^{m_{\bfd}}$ and 
$\bfs:[0,\tend]\to \R^{m_{\bfs}}$ such that:
\begin{equation}
   \left\{
   \begin{array}{rll}
        \dot{\bfd}(t) =& \bff_{\bfd}(t,\bfd(t),\bfs(t)) & \text{for every } t\in(0,\tend]\\
              \bfs(t) =& \bff_{\bfs}(t,\bfd(t),\bfs(t))& \text{for every } t\in(0,\tend]  \\
              \bfd(0) =& \bfd_0 & \\
              \bfs(0) =& \bfs_0 &
   \end{array}
    \right. \label{eq:final-semidiscrete}
\end{equation}

\subsection{Discretization in time}
We at last proceed with the integration in time of \eqref{eq:final-semidiscrete}. 
For some time step length, $0 < \delta t < 1 $, we discretize the time interval $[0,\tend]$ in $N\geq 1$ equal time steps, where $\tend = N \delta t = t_N$. We denote the time at each time step by
\begin{equation*}
    t_n = n \delta t ,\qquad n=0,1,\dots,N.
\end{equation*}
Then, for each $n$ we should solve the problem:
\begin{equation}
\begin{aligned}
(\text{STAGE 1})& & \bfd(t_n) =&
   \bfd (t_{n-1}) +
   \int\limits_{t_{n-1}}^{t_n}
    \bff_{\bfd}(\tau,\bfd(\tau),\bfs(\tau))
   d\tau\,; \\
(\text{STAGE 2})& & \bfs(t_{n}) =&
    \bff_{\bfs}(t_n,\bfd(t_n),\bfs(t_{n}))\,.
\end{aligned}
\label{eq:two-stage-exact}
\end{equation}
In \eqref{eq:two-stage-exact}, STAGE 1 corresponds to the integral form of the ODE for $\bfd$, and STAGE 2 is an implicit equation for $\bfs(t_n)$ that can be subsequently solved. The approximation in time is primarily motivated by the impossibility to evaluate exactly the integral on the RHS of the first line of \eqref{eq:two-stage-exact}.

Assuming that the evolution in time of $\bfd$ is sufficiently smooth, through Taylor expansion we can get
\begin{equation}
   \bfd (t_{n}) 
   =
   \bfd (t_{n-1}) + 
   \delta t\;
   \bff_{\bfd}\left[ 
   \frac{t_{n-1}+t_{n}}{2}, 
   \bfd\left( \frac{t_{n-1}+t_{n}}{2}\right), 
   \bfs\left( \frac{t_{n-1}+t_{n}}{2}\right) \right] 
    + \mcO(\delta t^3)\,.
    \label{eq:taylor_expansion}
\end{equation}
It is well known that, upon introducing the approximations
$$
\bfd\left( \frac{t_{n-1}+t_{n}}{2}\right) \approx 
\frac{\bfd(t_{n-1})+\bfd(t_{n})}{2}
\qquad
\text{and}
\qquad
\bfs\left( \frac{t_{n-1}+t_{n}}{2}\right) \approx 
\frac{\bfs(t_{n-1})+\bfs(t_{n})}{2} \,,
$$
and by dropping the $\mcO(\delta t^3)$ in \eqref{eq:taylor_expansion}, 
the \emph{implicit midpoint rule}
is obtained, which delivers an A-stable second order approximation \cite{gautschi2012numerical}, so that we can find a time-discrete version of $\bfd(t)$ and $\bfs(t)$. We denote the elements of this discretization by
\begin{equation*}
    \bfd_{n} \approx \bfd(t_n)\,,\qquad
    \bfs_{n} \approx \bfs(t_n)\,,
\end{equation*}
for every $n=0,1,\dots,N$ (notice that the case $n=0$ is given, since it is the initial condition of the problem).
For our case, the caveat of a direct implementation of such an approach is that, at every time step, a system of nonlinear equations would need to be solved, which can turn too costly for the amount of variables in play.

In this work, we apply a semi-implicit 
midpoint
rule as in \cite{boffi2019higher}: the values of the variables at $t_n$ on the RHS of the equations \eqref{eq:taylor_expansion} are first approximated with an explicit method. The integration scheme therefore reads as follows:
\begin{equation}
\begin{aligned}
(\text{STAGE 1})& & \bfd_{n} =&
   \bfd_{n-1} + 
   \delta t\;
   \bff_{\bfd}
   \left( 
   \frac{t_{n-1}+t_{n}}{2}, 
   \frac{\bfd_{n-1}+\bfd_{n}}{2}, 
   \frac{\bfs_{n-1}+\bfs_{n}}{2}
   \right) 
    \,; \\
(\text{STAGE 2})& & \bfs_{n} =&
    \bff_{\bfs}(t_{n},\bfd_{n}, \tilde{\bfs}_{n})\,.
\end{aligned}
\label{eq:two-stage-numerical}
\end{equation}
where the vectors $\tilde{\bfd}_{n}$ and $\tilde{\bfs}_{n}$ of the RHS correspond to
\begin{equation}
    \begin{aligned}
    \tilde{\bfd}_{n} =& 
    \left\{
    \begin{array}{ll}
       \bfd_0 + \delta t\, \bff_{\bfd}(0    , \bfd_0 , \bfs_0 )  & n=1 \,;\\
       2\bfd_{n-1} - \bfd_{n-2}  & n\geq 2 \,;
    \end{array}
    \right. \\
    \tilde{\bfs}_{n} =& 
    \left\{
    \begin{array}{ll}
       \bff_{\bfs}(t_{n},\tilde{\bfd}_1 , \bfs_0)  & n=1 \,; \\
       2\bfs_{n-1} - \bfs_{n-2}  & n\geq 2 \,.
    \end{array}
    \right.  
    \end{aligned}
    \label{eq:predictors}
\end{equation}

If restricted to STAGE 1 of the scheme \eqref{eq:two-stage-numerical}-\eqref{eq:predictors}, the integration method for $\bfd$ may be regarded as a predictor-corrector method. 
In \S 6.2.3 of \cite{gautschi2012numerical}, it is explained that if an explicit method of order $k$ is first used (predictor method), and it is followed by an implicit method of the same order (corrector method), then the full method is of order $k$ as well. 
It is easy to check too, that the full scheme keeps order $k$ when the predictor is of order $k-1$. 
In \eqref{eq:two-stage-numerical}-\eqref{eq:predictors}, we have $k=2$ since the corrector is the 
midpoint
rule, while there is one predictor for the first time step and another for the subsequent steps. For the first time step, the predictor is Euler method (which is of first order). Following \cite{boffi2019higher}, starting at $n=2$ we use as a predictor a simple formula of $\mcO({\delta t^2})$ error (which amounts to linear local truncation error). The formula (i.e., the $n\geq2$ cases in \eqref{eq:predictors}) comes from the centered finite difference approximation of the second derivative, namely,
$$
\Ddot{\bfd}(t_{n-1})=\frac{1}{\delta t^2}\left(\bfd(t_{n-2})- 2 \bfd(t_{n-1}) + \bfd(t_{n})\right) + \mcO({\delta t^2}) 
\quad \Rightarrow\quad
\bfd({t_{n}}) = 2 \bfd(t_{n-1}) - \bfd(t_{n-2}) + \mcO({\delta t^2}).
$$
Obviously, the price to pay for this approach is the storage in memory of the solutions at $t_{n-2}$ as well.
\begin{remark}
    Even though in \cite{boffi2019higher} this scheme performs the second stage a single time, we can also regard the predictor as a convenient seed for a fixed-point method where the corrector is the iteration operator. This view is favored in the numerical experiments reported below.
\end{remark}

\subsection{Convergence of the full scheme}

Regarding the fully discrete scheme, take as a sample the complete equation for $c_s$ at time $t=t_n$. In its expanded form, it can be written like this:
\begin{multline}
    \left( r_s^{hp} , c_{s;n}^{hp} \right)_{\Omega_s}
    + \delta t\;\left( \nabla r_s^{hp} , D_s(\tilde{c}_{s;\text{mid}}^{hp},\tilde{\pi}^{hp}_{\text{mid}})\nabla c_{s;n}^{hp} \right)_{\Omega_s}
        =
    \left( r_s^{hp} , c_{s;n-1}^{hp} \right)_{\Omega_s}
    -\delta t\;\left( \nabla r_s^{hp} , D_s(c_{s;\text{mid}}^{hp},\pi^{hp}_{\text{mid}})\nabla c_{s;n-1}^{hp} \right)_{\Omega_s} \\
    -\delta t\;\frac{1}{F}\langle r_s^{hp} ,\tilde{\ibv}_{;\text{mid}}^{hp}
    \rangle_{\Gamma_s\cap\interface}\,,
    \label{eq:fully_discrete_concs}
\end{multline}
where $r_s^{hp}$ is the test function, and $c_{s;n}^{hp}$ is the unknown, both belonging to $\scV^{hp}_{c_s}$, every $n$ or $n-1$ subscript represents to which time step the quantity corresponds, the subscript $_\text{mid}$ means that the quantity is evaluated at the average of solutions of step $n$ and $n-1$, and the tilde indicates that such values are computed with the predictors \eqref{eq:predictors}. At this point, problem \eqref{eq:fully_discrete_concs} can be seen as a variational problem which is linear with respect to $c_{s;n}^{hp}$ that has been discretized in space with a Galerkin method. Even in the infinite-dimensional case (i.e., before the spatial $hp$-discretization), it is easy to check both the continuity and coercivity of the bilinear form at the LHS, as well as the continuity of the RHS; that is,
$$
C^{\text{coer}}_{c_s}\| c_{s;n} \|^2_{H^1(\Omega_s)}
\leq
\left( c_{s;n} , c_{s;n} \right)_{\Omega_s}
    + \delta t\left( \nabla c_{s;n} , D_s(\tilde{c}_{s;\text{mid}},\tilde{\pi}_{\text{mid}})\nabla c_{s;n} \right)_{\Omega_s}\;,
$$
$$        
    \left( r_{s} , c_{s;n} \right)_{\Omega_s}
        + \delta t \left( \nabla r_{s} , D_s(\tilde{c}_{s;\text{mid}},\tilde{\pi}_{\text{mid}})\nabla c_{s;n} \right)_{\Omega_s}
    \leq 
    C^{\text{LHS}}_{c_s}
    \| r_{s} \|_{H^1(\Omega_s)} 
    \| c_{s;n} \|_{H^1(\Omega_s)}\;,
$$
$$
    \left( r_s , c_{s;n-1} \right)_{\Omega_s}
    -\delta t \left( \nabla r_s , D_s(c_{s;\text{mid}},\pi_{\text{mid}})\nabla c_{s;n-1} \right)_{\Omega_s}
    -\tfrac{\delta t}{2}\frac{1}{F}\langle r_s , \tilde{\ibv}_{;\text{mid}}
    \rangle_{\Gamma_s\cap\interface}
    \leq 
    C^{\text{RHS}}_{c_s}
    \| r_{s} \|_{H^1(\Omega_s)}\;,
$$
where $C^{\text{coer}}_{c_s}$, $C^{\text{LHS}}_{c_s}$ and $C^{\text{RHS}}_{c_s}$ are positive constants, independent of $r_{s;n} , c_{s;n}$, but may depend on previous values of $c_s$ and other unknowns and material properties.
According to Lax-Milgram's theorem, the discrete-in-time weak problem at $t=t_n$ is well-posed (with solution $c_{s;n}\in\scV_{c_s}$), and given the conformity of $\scV^{hp}_{c_s}$, Cea's lemma guarantees that $c_{s;n}^{hp}$, that is, the discrete solution to \eqref{eq:fully_discrete_concs}, converges in a stable fashion to $c_{s;n}$ \cite{demkowicz2023mathematical}. A similar analysis extends to all the remaining equations of our formulation.

Based on the numerical experiments and the theoretical results on the semi-implicit 
midpoint
rule available in \cite{boffi2019higher} (therein applied to a fluid-structure interaction problem), we expect stability in the entire run of our scheme \eqref{eq:two-stage-numerical}. If the solution is sufficiently regular, both in space and time, we can gather the order of accuracy of the 
midpoint
rule and the usual approximation properties of finite element spaces, to postulate an expected \emph{a priori} error estimate of the form
\begin{equation}
    \left\| \left(\theta(t_N,\cdot), c_s(t_N,\cdot), c_e(t_N,\cdot), \phi_s(t_N,\cdot), \phi_e(t_N,\cdot), \bsu(t_N,\cdot)\right) 
    -
    \left(\theta_N^{hp}, c_{s;N}^{hp}, c_{e;N}^{hp}, \phi_{s;N}^{hp}, \phi_{e;N}^{hp}, \bsu_{N}^{hp}\right)
    \right\|_{\scV}
    \leq C
    h_{\max}^{p} \delta t ^2\,
    \label{eq:full-error-conjec}
\end{equation}
where $h_{\max}$ is the maximum diameter of the elements in $\mesh^{hp}$, $p$ is assumed to be uniform, and the constant $C$ may depend on $p$. Proving this estimate goes beyond the scope of the present paper.

\begin{remark}
    We can mention two technical details that can shed light on the validity of the above estimate. First, regarding the coercivity and continuity constants, $C^{\text{coer}}_{c_s}$ and $C^{\text{LHS}}_{c_s}$, respectively. The existence of these constants, as positive real numbers, is due to the fact that the variable mass diffusivity of the electrodes is bounded below and above (see \eqref{eq:stress-diffusivity}), namely, the value of $D_s(\tilde{c}_{s;n},\tilde{\pi}_{n})$ ranges from $e^{-\beta_D}D_{s,0}$ through $e^{\alpha_D}D_{s,0}$. The equations for $c_s$ are the only ones that involve (non-linearly) varying diffusion coefficients, thereby the satisfaction of coercivity and continuity in the remaining equations is certainly more immediate. 
    
    On the other hand, it is well-known that the well-posedness of initial-value problems is hardly attainable without a Lipschitz condition on the rate of change functions, say, function $\bff_{\bfd}$ in \eqref{eq:final-semidiscrete}. In our formulation, we have a function within $\bff_{\bfd}$ that fails to satisfy that condition but within some bounds. We are speaking of $I_c$, defined in \eqref{eq:butler-volmer}. This issue is caused by the presence of the square root (evaluated at $c_e$, at $c_{s,\max}-c_s$ and at $c_s$), which is not Lipschitz at zero. For that reason, we must restrain ourselves from allowing too small values of $c_e$, as well as values of $c_s$ too close both to zero and to $c_{s,\max}$.
    If this restriction is properly applied and preserved along all time steps, there should be no other limitation to have a second order approximation in time.
    \label{rk:about-estimate}
\end{remark}

%
\section{Numerical examples} 
\label{sec:results}
\subsection{Materials and parameters}
We have chosen the following materials for the electrodes and the electrolytic medium.
\begin{itemize}
    \item Solid anode material: Lithiated graphite (Li$_x$C$_6$).
    \item Electrolyte: Lithium hexafluorophosphate (LiPF$_6$) in ethylene carbonate-diethyl carbonate (EC-DEC).
    \item Solid cathode material: Lithium manganese oxide (Li$_x$Mn$_2$O$_4$).
\end{itemize}
There is plenty of data for this set of substances, so that we can carry out the simulation of our multiphyisics model.
All the parameters that we require to run the simulation are given in Table \ref{tab:parameters}.
\begin{table}[H]
\centering
\begin{tabular}{|c|c|c|c|l|}
\hline
\textbf{Parameter}  & \textbf{SI units}                     &\textbf{Value}&\textbf{Source}                & \textbf{Notes}           \\ \hline
$\rho_{sa}C_{v,sa}$ & kg m$^{-1}$ s$^{-2}$ K$^{-1}$         & 3.8235$\times 10^{6}$  & \cite{miranda_effect_2019}        & Li$_x$C$_6$                          \\
$E_{sa}$            & kg m$^{-1}$ s$^{-2}$                  & 3.64$\times 10^{9}$  & \cite{wang2018mechanical}         & Similar anode material               \\
$\nu_{sa}$          & -                                     & 0.3  & \cite{bai2019two}                 & Extended from cathode              \\
$\alpha_{sa}$       & K$^{-1}$                              & 1$\times 10^{-5}$  & \cite{carlstedt_coupled_2022}     & An estimate based on theirs        \\
$\omega_{sa}$       & m$^3$ mol$^{-1}$                      & 3.499$\times 10^{-6}$  & \cite{bai2019two}                 & Extended from cathode              \\
$\gamma_{sa}$       & kg$^{-1}$ m$^{-3}$ s$^3$ A$^2$        & 100  & \cite{miranda_effect_2019}        & Li$_x$C$_6$                          \\
$\lambda_{sa}$      & kg m s$^{-3}$ K$^{-1}$                & 1.04  & \cite{miranda_effect_2019}        & Li$_x$C$_6$                          \\
$D_{0;sa}$          & m$^2$ s$^-1$                          & 3.9$\times 10^{-14}$  & \cite{miranda_effect_2019}        & Li$_x$C$_6$                          \\
$c_{sa,\max}$       & mol m$^{-3}$                          & 3.1507$\times 10^{4}$  & \cite{miranda_effect_2019}        & Li$_x$C$_6$                          \\
$\rho_{sc}C_{v,sc}$ & kg m$^{-1}$ s$^{-2}$ K$^{-1}$         & 9.0371$\times 10^{5}$  & \cite{miranda_effect_2019}        & Li$_x$Mn$_2$O$_4$                    \\
$E_{sc}$            & kg m$^{-1}$ s$^{-2}$                  & 2.5$\times 10^{9}$  & \cite{wang2018mechanical}         & Similar cathode material \\
$\nu_{sc}$          & -                                     & 0.3 & \cite{bai2019two}                 & Li$_x$Mn$_2$O$_4$                    \\
$\alpha_{sc}$       & K$^{-1}$                              & 1$\times 10^{-5}$  & \cite{carlstedt_coupled_2022}     & An estimate based on theirs        \\
$\omega_{sc}$       & m$^3$ mol$^{-1}$                      & 3.499$\times 10^{-6}$  & \cite{bai2019two}                 & Li$_x$Mn$_2$O$_4$                  \\
$\gamma_{sc}$       & kg$^{-1}$ m$^{-3}$ s$^3$ A$^2$        & 3.8  & \cite{miranda_effect_2019}        & Li$_x$Mn$_2$O$_4$                  \\
$\lambda_{sc}$      & kg m s$^{-3}$ K$^{-1}$                & 6.2  & \cite{miranda_effect_2019}        & Li$_x$Mn$_2$O$_4$                  \\
$D_{0;sc}$          & m$^2$ s$^-1$                          & 1$\times 10^{-13}$  & \cite{miranda_effect_2019}        & Li$_x$Mn$_2$O$_4$                  \\
$c_{sc,\max}$       & mol m$^{-3}$                          & 2.286$\times 10^{4}$  & \cite{miranda_effect_2019}        & Li$_x$Mn$_2$O$_4$                  \\
$\rho_{e}C_{v,e}$   & kg m$^{-1}$ s$^{-2}$ K$^{-1}$         & 1.9979$\times 10^{6}$  & \cite{miranda_effect_2019}        & LiPF$_6$ in EC-DEC                 \\
$\kappa_{e}$        & kg$^{-1}$ m$^{-3}$ s$^3$ A$^2$        & 0.2  & \cite{trembacki_fully_2015}       & LiPF$_6$ in EC-DEC                 \\
$\lambda_{e}$       & kg m s$^{-3}$ K$^{-1}$                & 0.344  & \cite{miranda_effect_2019}        & LiPF$_6$ in EC-DEC                 \\
$D_{e}$             & m$^2$ s$^{-1}$                        & 7.5$\times 10^{-11}$  & \cite{trembacki_fully_2015}       & LiPF$_6$ in EC-DEC                 \\
$t_+$               & -                                     & 0.363  & \cite{trembacki_fully_2015}       & LiPF$_6$ in EC-DEC                 \\
$R$                 & kg m$^2$ s$^{-2}$ mol$^{-1}$ K$^{-1}$ & 8.314462618 & \cite{linden2002handbook}         & Universal constant                 \\
$F$                 & A s mol$^{-1}$                        & 96485.33212 & \cite{linden2002handbook}         & Universal constant                 \\
$k_{\text{BV}}$     & m$^{2.5}$ mol$^{-0.5}$ s$^{-1}$       & 1.1$\times 10^{-11}$     & \cite{trembacki_fully_2015}       & See Table 2 of reference           \\
$\alpha_D$                 & -                                     & 6         & \cite{gritton_using_2017}         & An estimate based on theirs        \\
$\beta_D$                 & -                                     & 1.5         & \cite{gritton_using_2017}         & An estimate based on theirs        \\
$\pi_{\max}$        & kg m$^{-1}$ s$^{-2}$                  & 1$\times 10^{9}$  & \cite{gritton_using_2017}         & An estimate based on theirs        \\
$\theta_0$          & K                                     & 298.15  & -                     & Frequent value \\
$c_{e,0}$           & mol m$^{-3}$                          & 2$\times 10^{3}$  & \cite{trembacki_fully_2015}       & See Table 1 of reference          \\
\hline
\end{tabular}
\caption{Parameter values used in the simulations}
\label{tab:parameters}
\end{table}
The experimental models for the open-circuit potential are taken from \cite[Appendix A]{fuller_simulation_1994}. The anode follows the function
\begin{equation*}
    \ocpsa(\soca)=-0.16+1.32e^{-3\soca}+10 e^{-2000\soca},
\end{equation*}
while the cathode open-circuit potential is described by
\begin{equation*}
\begin{aligned}
    \ocpsc(\socc) =
    & 4.06279 + 0.0677504\tanh{(-21.8502\socc + 12.8262)} \\
    & -0.105734\left(\frac{1}{(1.00167-\socc)^{0.379571}}-1.576\right) \\
    & -0.045 e^{-71.69\socc^8} +0.01e^{-200(\socc-0.19)}
\end{aligned}
\end{equation*}

\subsection{Simulation scenarios}

\begin{figure}
    \includegraphics[width=0.2\textwidth]{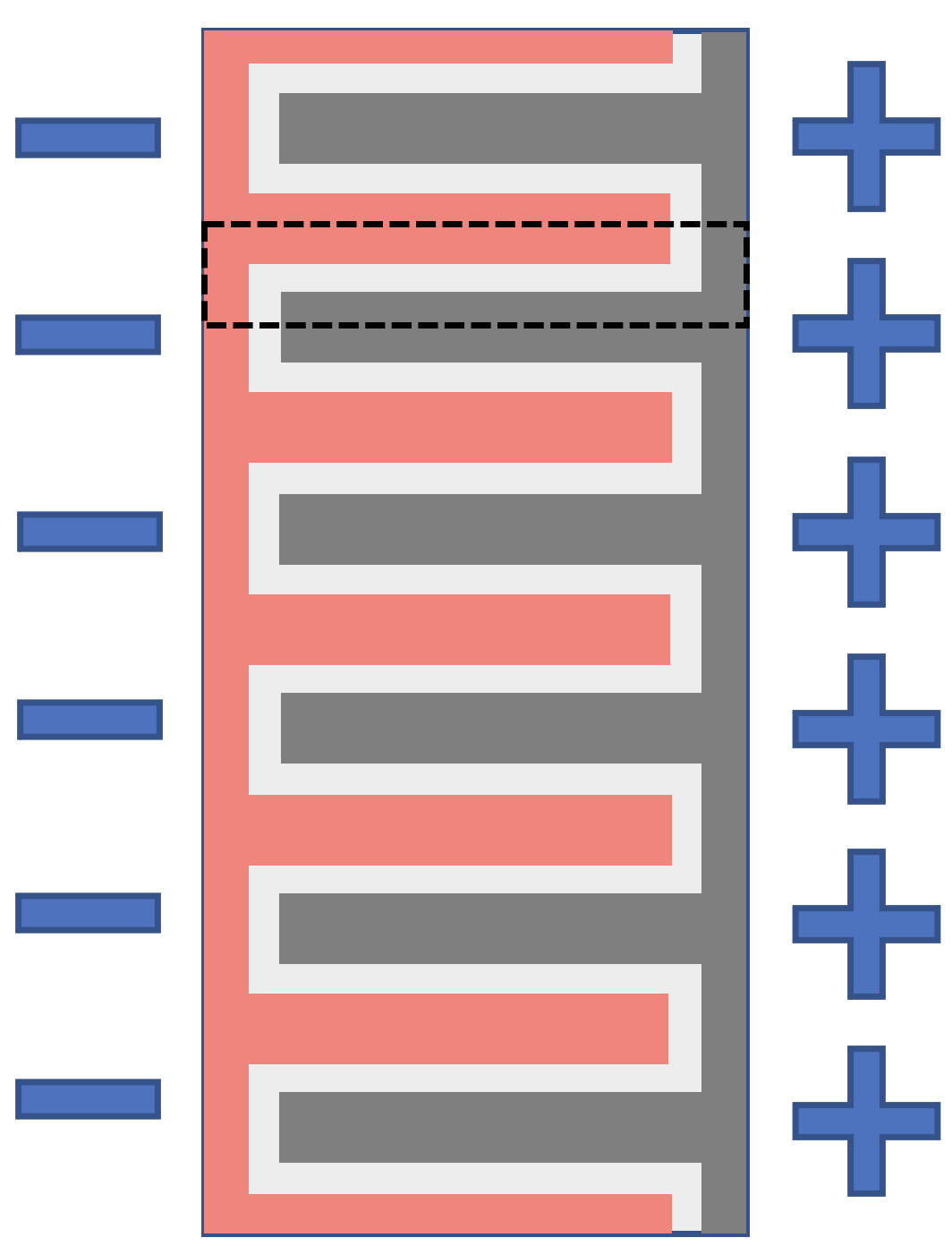}
    \hfill
    \includegraphics[width=0.7\textwidth]{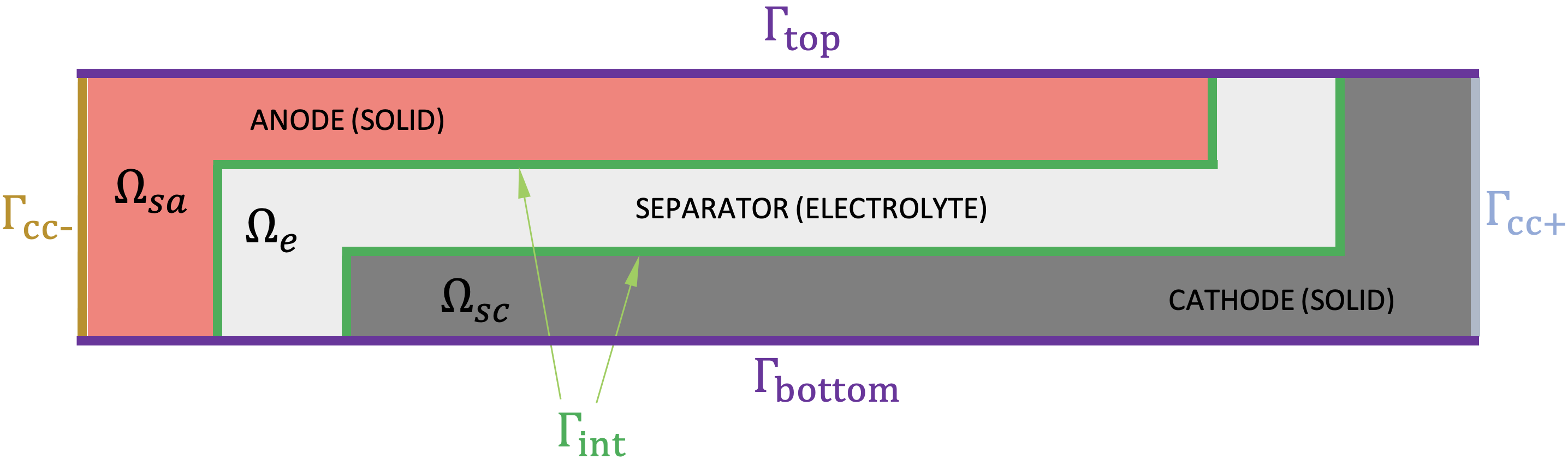}
\caption{Diagram of an interdigitated battery cell structure: (left) multiple digits, showing a representative geometry with a dashed line; (right) detail of the representative geometry, specifying subdomains and boundaries.}
\label{fig:interdigitated-geometry}
\end{figure}

In order to see the performance of our model in a near realistic battery geometry, we have identified a certain design that has been studied by previous researchers \cite{trembacki_fully_2015,miranda_effect_2019}. It can be referred to as the \emph{interdigitated plates} design. Figure \ref{fig:interdigitated-geometry} shows how the structure of this battery cell design looks like, along with a repeating detailed feature that has 
been chosen as the computational domain for our simulations. We remark that the location of the subdomains, the interface and the boundary partition are all specified in the figure. The boundary conditions are applied to the corresponding boundary parts as explained above. 

\begin{figure}[!ht]
    \centering
    \includegraphics[width=0.65\textwidth]{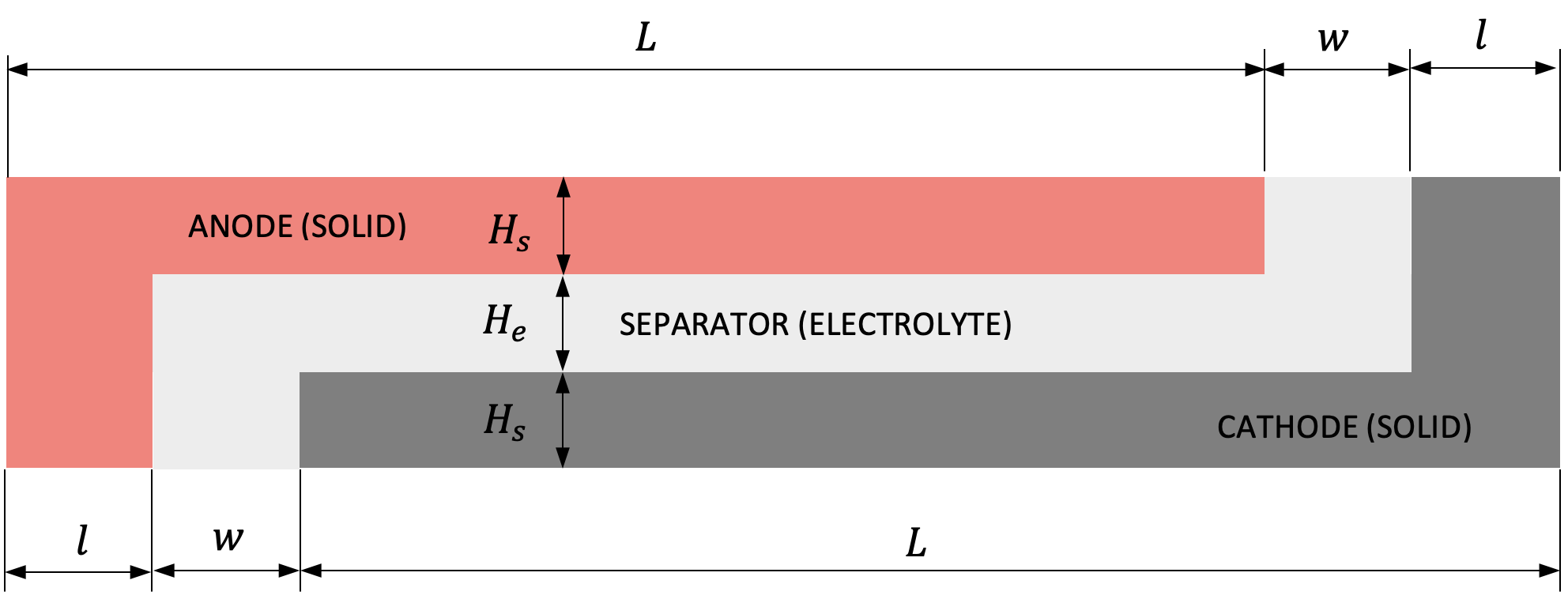}
    \caption{Main dimensions of the representative geometry of an interdigitated battery.}
    \label{fig:dimensions_drawing}
\end{figure}

The representative geometry herein chosen has the dimensions given in Figure \ref{fig:dimensions_drawing} and Table \ref{tab:dimensions_table}.
\begin{table}[!ht]
\centering
\begin{tabular}{|l|l|l|l|l|l|}
\hline
\textbf{Parameter} & $H_s$ & $H_e$ & $L$ & $w$ & $l$ \\ \hline
\textbf{Value}     & 30    & 40    & 900 & 40  & 60  \\ \hline
\end{tabular}
\caption{Dimensions of the interdigitated plates design. All values in {\textmu}m}
\label{tab:dimensions_table}
\end{table}

The applied current density on $\Gpos$ will vary according to each of the following scenarios:
\begin{itemize}
    \item low current discharge: applied current density $\iapp=$5A/m$^2$, end time $\tend=4$ hours;
    \item high current discharge: applied current density $\iapp=$20A/m$^2$, end time $\tend=1$ hour;
    \item low current charge: applied current density $\iapp=$-5A/m$^2$, end time $\tend=4$ hours;
    \item high current charge: applied current density $\iapp=$-20A/m$^2$, end time $\tend=1$ hour.
\end{itemize}
All scenarios start at a 50\% state of charge for both electrodes.

Moreover, the four scenarios are also run at a strain-free, isothermal condition, which means that the model effectively covers just the concentration and electric potential evolution. Below, we refer to this as the \emph{electrochemical model}. By adding these with the four runs of the full multi-physics model, we get a total of eight simulations.

\subsection{Quantities of interest}
\label{sec:qoi}
A transient multidimensional simulation of a multiphysics model like the present one outputs a very big amount of data. The presentation of such results, besides field visualizations, can be aided by simpler observables or quantities of interest that summarize much of the relevant information generated by the computer simulations. Next, we show some quantities of interest that must be post-processed.
\begin{enumerate}
    \item Cell output voltage:
        \begin{equation}
            V_{\text{out}}(t)= \iapp \int\limits_{\Gpos} \phi_{sc}(\cdot,t)   dA.
            \label{eq:v-out}
        \end{equation}
    \item Average electric potential in electrolyte:
        \begin{equation}
            \overline{\phi_e}(t)= \int\limits_{\Omega_{e}} \phi_{e}(\cdot,t)   dx.
            \label{eq:phie-avg}
        \end{equation}
    \item Average state of charge at each electrode:
        \begin{equation}
            \overline{\soca}(t)= \int\limits_{\Omega_{sa}} \soca(\cdot,t)   dx\,,
            \qquad\qquad 
            \overline{\socc}(t)= \int\limits_{\Omega_{sc}} \socc(\cdot,t)   dx.
            \label{eq:soc-avg}
        \end{equation}
    \item Average temperature:
        \begin{equation}
            \overline{\theta}(t)= \int\limits_{\Omega} \theta(\cdot,t)   dx.
            \label{eq:tmprt-avg}
        \end{equation}
\end{enumerate}
An important measure in battery engineering is the \emph{power density} $\overline{P}$, which is usually computed for the discharge scenarios.
\begin{equation}
    \overline{P}=\int\limits_{0}^{\tend} V_{out}(t) \iapp  dt.
    \label{eq:power-energy-density}
\end{equation}
Even though we do not intend to assess differences in performance as the design changes, we plan to compute this quantity of interest to know how close it is when using the electrochemical model with respect to the the full multiphysical model.

In addition to all of the above, the mechanical behavior is observed by means of displacements and stresses. To simplify things, we will look at the maxima of these variables at end time. Since the stress is a second order tensor, we turn it into a scalar quantity by means of the Von Mises equivalent tensile stress \cite{anand_continuum_2020}:
\begin{equation*}
    \sigma_{_{VM}}=
    \sqrt{\frac{
    (\sigma_{11}-\sigma_{33})^2+(\sigma_{22}-\sigma_{11})^2+(\sigma_{33}-\sigma_{22})^2
    }{2}
    +3(\sigma_{12}^2+\sigma_{23}^2+\sigma_{31}^2)}
\end{equation*}
Recalling that in the present 2D simulations we hold a plane-strain assumption, shear components $\sigma_{23}$ and $\sigma_{31}$ are nill, while the third axial stress component $\sigma_{33}$ can be calculated through \eqref{eq:me-strain}-\eqref{eq:generalized_Hooke} with $\epsilon_{33}=0$.

\subsection{Implementation details}
The implementation and the solution of the present discrete model is built on top of the $hp2D$ finite element library written in Fortran, which allows for the specification of multiple physical attributes and subdomains where they are active or inactive, along with finite element spaces that are conforming to all the energy spaces (those of the exact Hilbert complex $H^1\to H(\curl)\to H(\div)\to L^2$) with variable anisotropic order throughout the mesh, local adaptive-refinement capabilities, and transfinite interpolation among other features \cite{hpbook,GattoDemkowicz10,Fuentes2015}. The linear systems at every step are solved with the MUMPS library \cite{amestoy2000mumps}, and the assembly is accelerated thanks to the usage of Open-MP shared memory parallelization.

The presence of functions with a restricted domain, such as the exchange current density $I_c$ in \eqref{eq:butler-volmer}, may induce runtime errors if the new time step solution falls out of bounds (see remark \ref{rk:about-estimate}). Apart from restricting the limits where our simulation can start, this has forced us to run a few extra fixed point iterations after the corrector step is made, at each time step. This increase of iterations keeps correcting the solution, preventing it from approaching too much to the bounds. The lack of extra iterations may cause some NaNs to appear after trying the computation of the square root of negative numbers, which directly kills the simulation.

To start the time integration process, we have run two unloaded time steps to fill solution degrees of freedom (DOF) at equilibrium. This is not a mandatory step, but we have made this choice as it diminishes voltage overshoots and undershoots at the first few time steps when the electrical load is applied.

\begin{figure}[!ht]
    \centering
    \includegraphics[width=\textwidth]{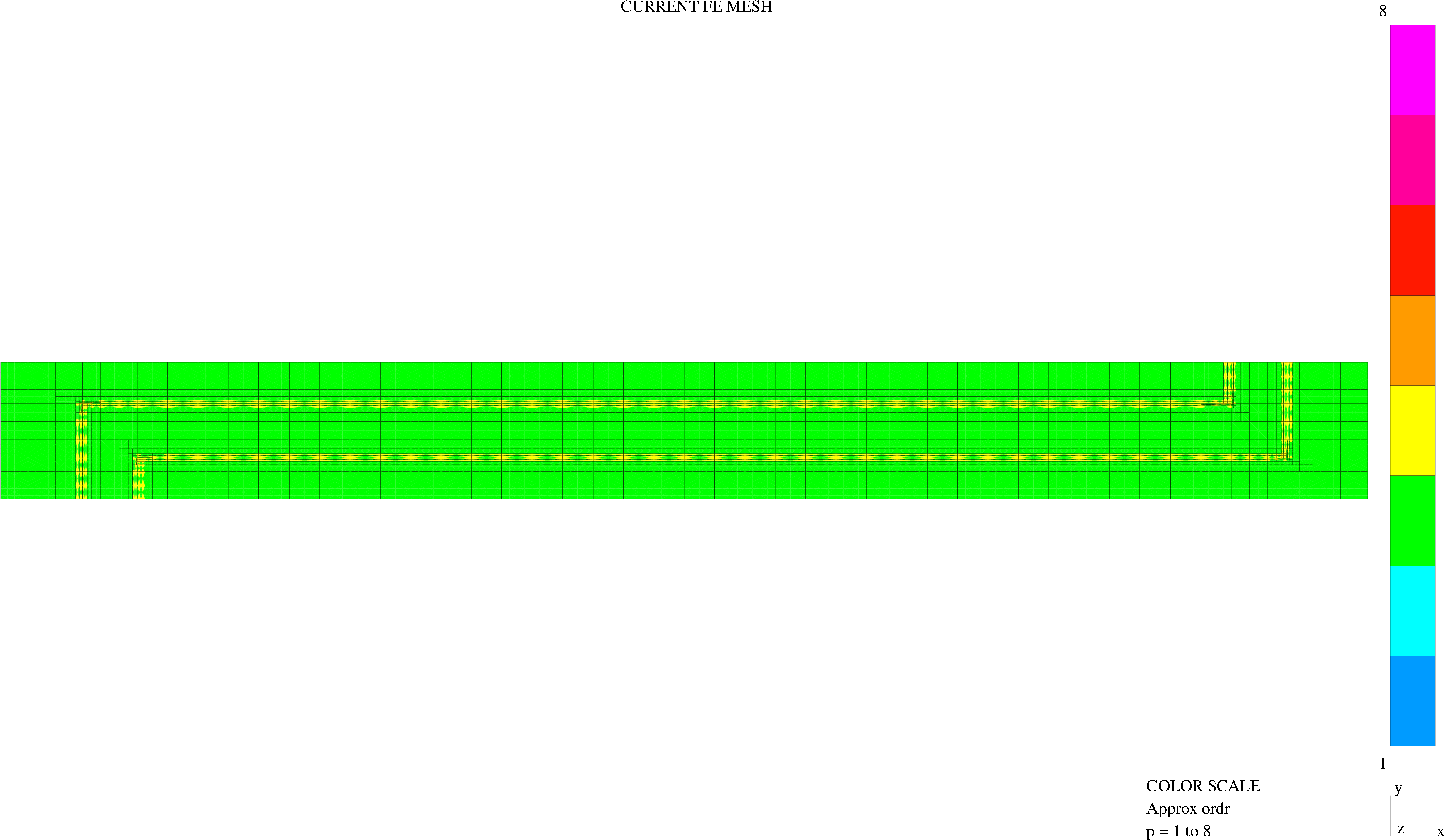}
    \caption{Finite element mesh used for all the numerical experiments herein reported. }
    \label{fig:interdig-mesh}
\end{figure}

Finally, we have imposed a layered $hp$ finite element mesh to improve accuracy near the interface and the boundaries, which is the part of the domain where the most complex behavior arises, both physically and numerically. These layers were generated with thinner elements of gradually higher polynomial degree closer to the interface. The result is an anisotropic mesh, since the refinements are made only in the direction normal to the boundary. 
For the computations reported below, the mesh comprises 1148 quadrilateral elements, and a total of 56232 degrees of freedom. We show said mesh in Figure \ref{fig:interdig-mesh}. Notice how the element size gets finer near the interface. The color scale indicates the order, which is cubic throughout most of the domain, but is quartic in the normal direction on the most refined region.

The utilization of these layered $hp$ meshes has diminished overshoots and undershoots, and decreases the number of extra fixed point iterations required to avoid cases of NaN appearances in the computation of $I_c$. The number of extra fixed point iterations of all the results below is four (besides the first corrector step). A convergence study focused on this type of meshes for complex problems, such as the present application, may be the subject of a future paper.

Finally, given the vast diversity of magnitudes of the variables and parameters at play (see Table \ref{tab:parameters}), an appropriate dimensional analysis was carried out, so that the magnitude range of the numbers within the computations got significantly diminished. This caused, for instance, that the practical length unit was $10^{-4}$ m, and that the practical time unit became $1$ minute. With that in mind, the time step length herein chosen ranged from 0.05 to 0.1 units, e.g., between 3 and 6 seconds. The longest simulation required 4800 time steps.

\subsection{Results}

\begin{figure}[!ht]
    \includegraphics[width=0.49\textwidth]{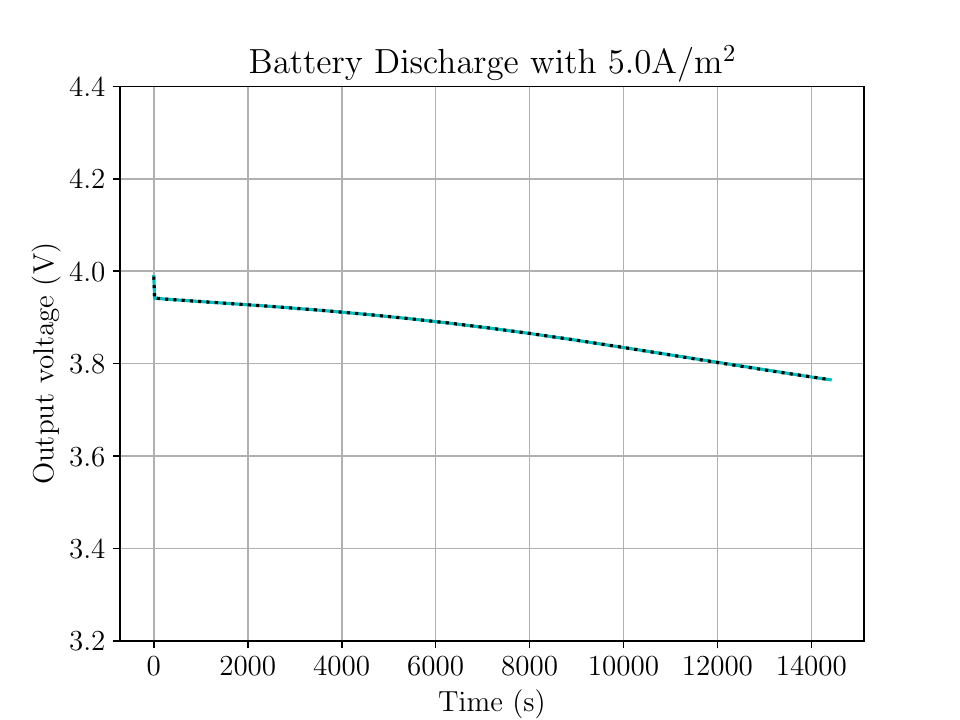} \hfill
    \includegraphics[width=0.49\textwidth]{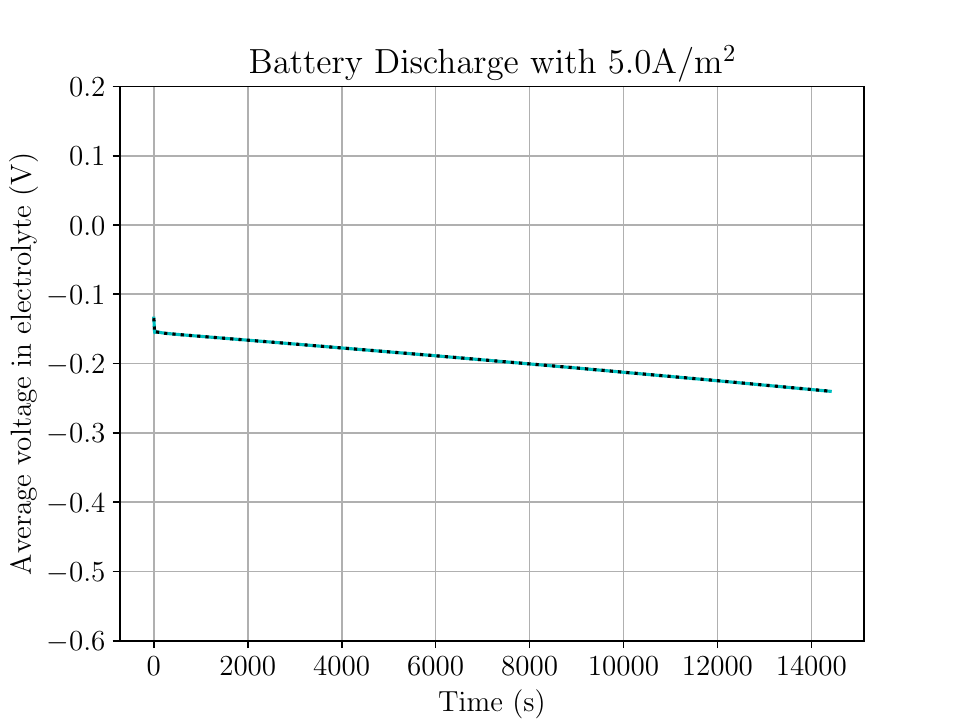} \\
    \includegraphics[width=0.49\textwidth]{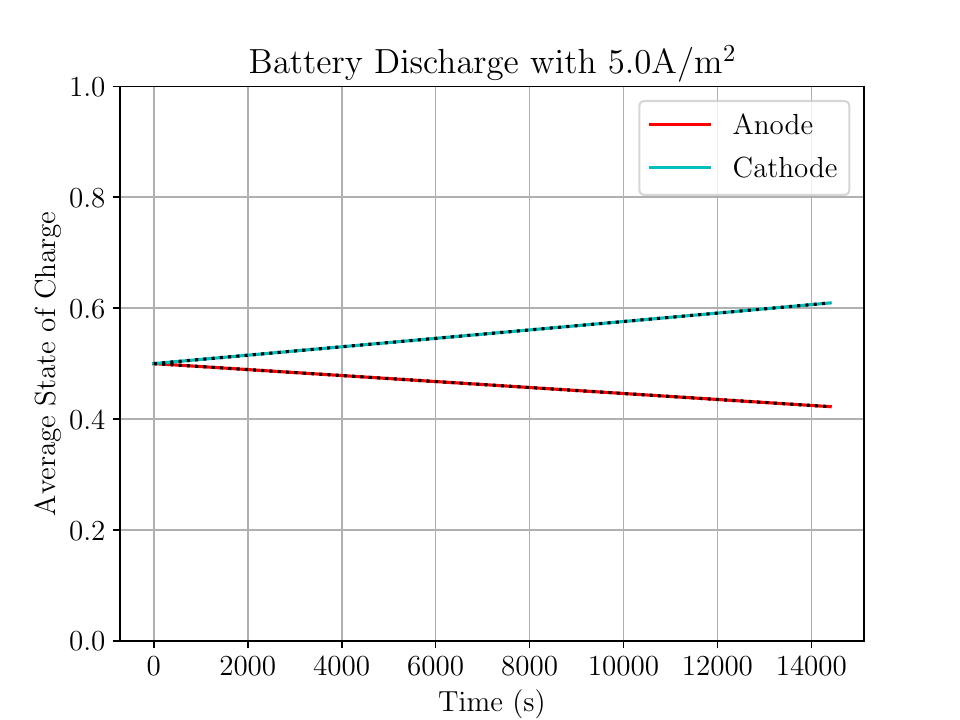} \hfill
    \includegraphics[width=0.49\textwidth]{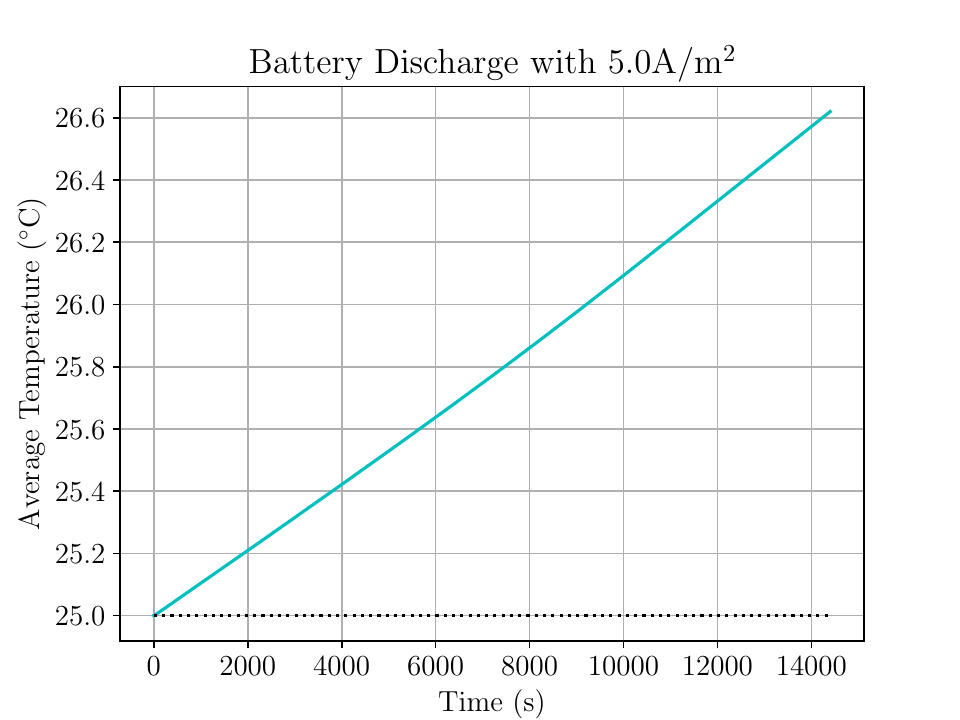}
    \caption{Low current discharge scenario. 
    Results comparison between full multiphysical model (solid lines)
    and electrochemical model (dotted lines): 
    (top left) output voltage vs. time; 
    (top right) average potential in electrolyte vs. time;
    (bottom left) average electrode state of charge vs. time;
    (bottom right) average temperature (in Celsius degrees) vs. time}
    \label{fig:b2dl_comparison_plots}
\end{figure}

\begin{figure}[!ht]
    \includegraphics[width=0.49\textwidth]{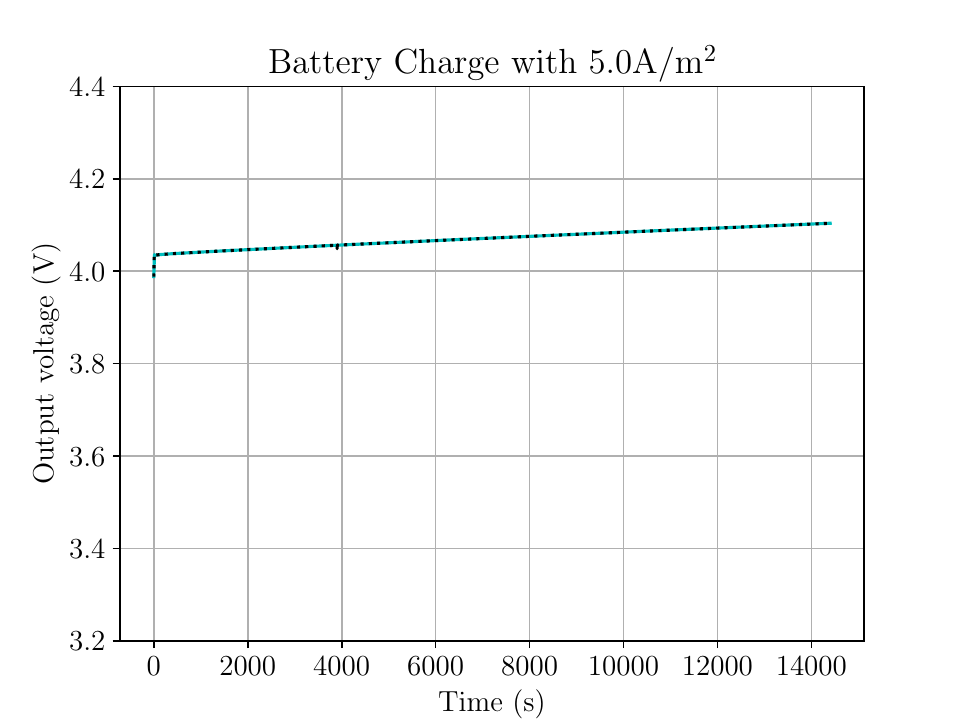} \hfill
    \includegraphics[width=0.49\textwidth]{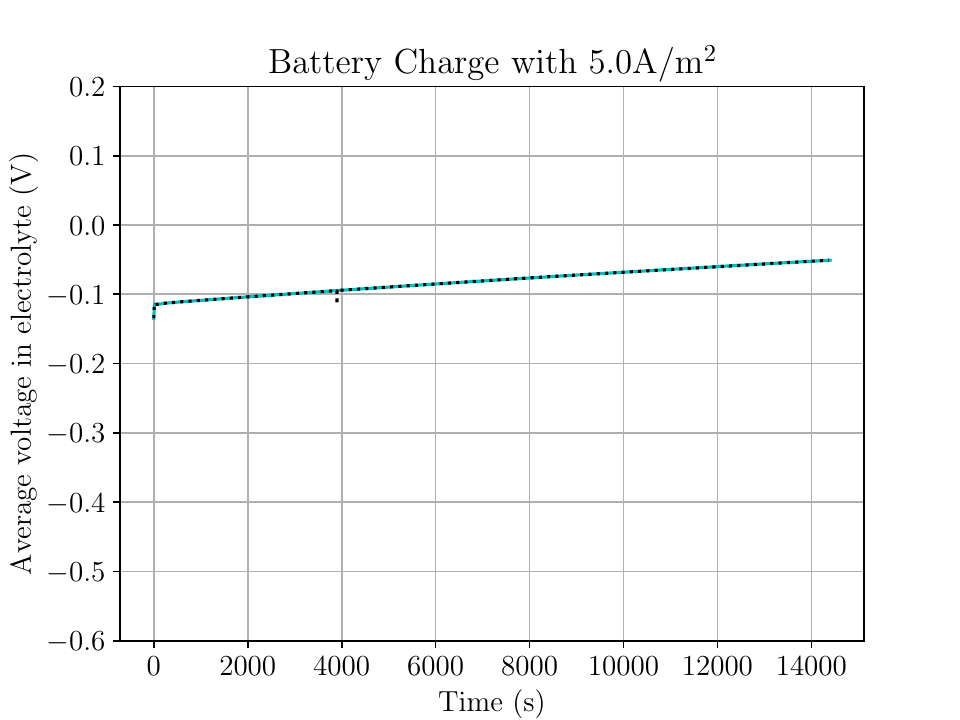} \\
    \includegraphics[width=0.49\textwidth]{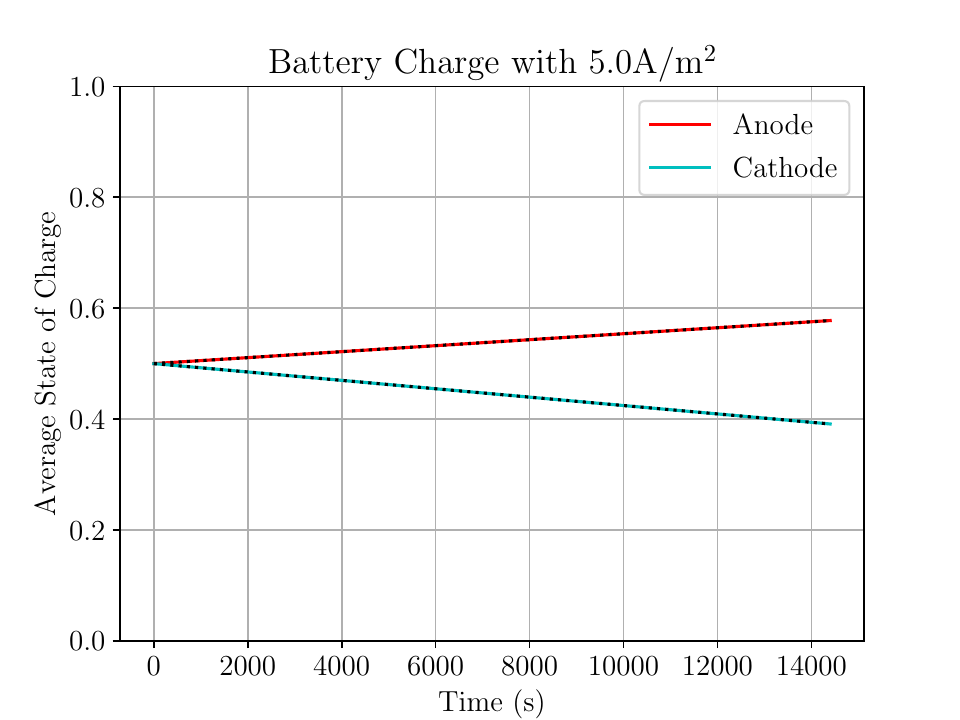} \hfill
    \includegraphics[width=0.49\textwidth]{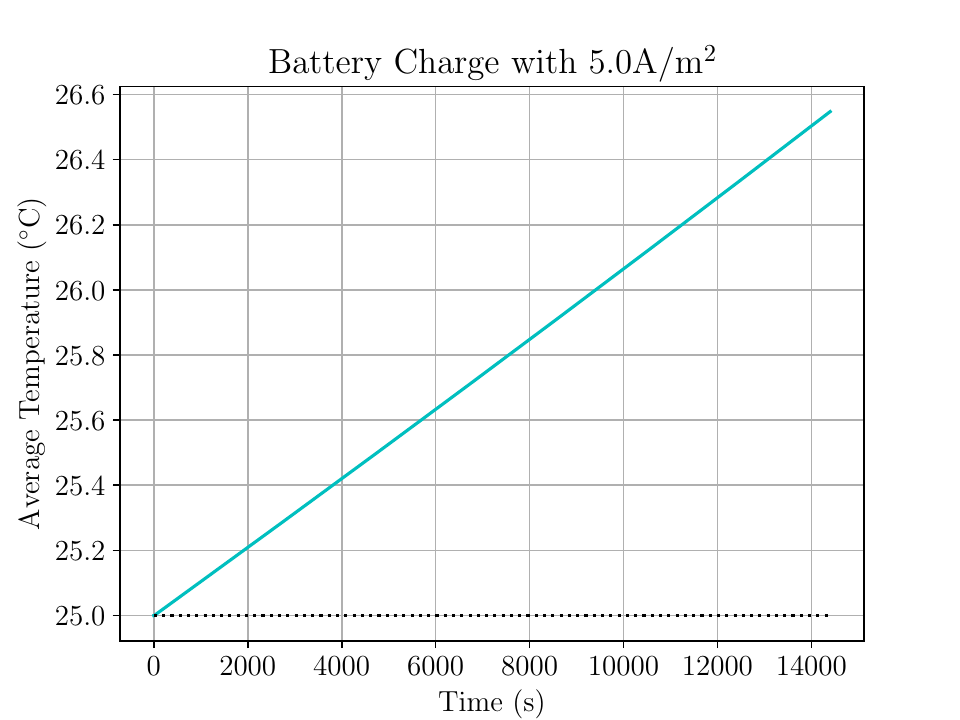}
    \caption{Low current charge scenario. 
    Results comparison between full multiphysical model (solid lines)
    and electrochemical model (dotted lines): 
    (top left) output voltage vs. time; 
    (top right) average potential in electrolyte vs. time;
    (bottom left) average electrode state of charge vs. time;
    (bottom right) average temperature (in Celsius degrees) vs. time}
    \label{fig:b2cl_comparison_plots}
\end{figure}

\begin{figure}[!ht]
    \includegraphics[width=0.49\textwidth]{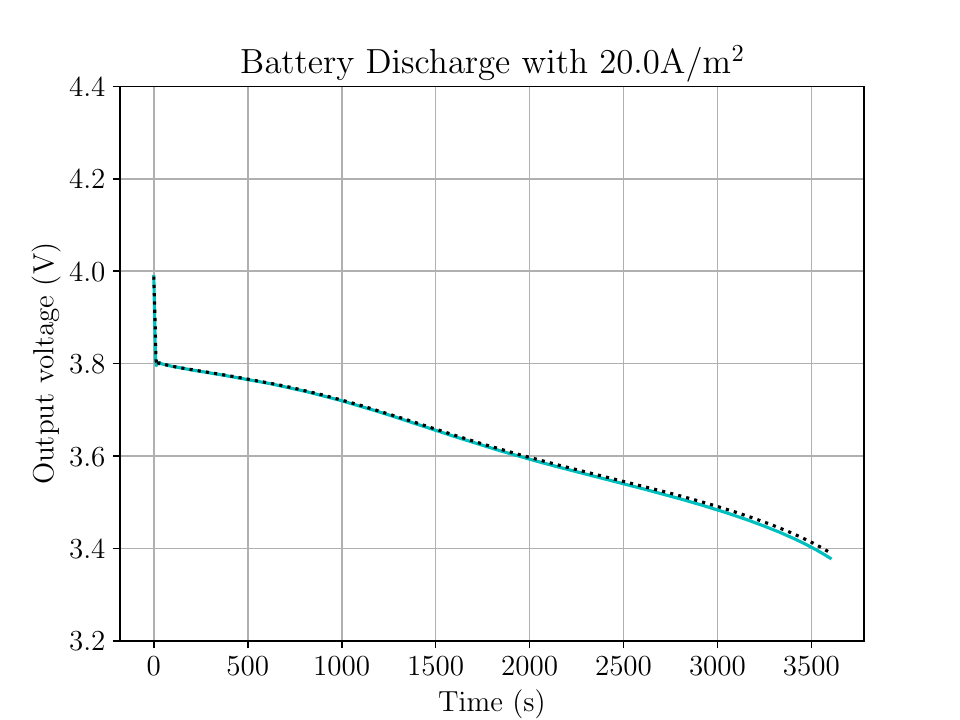} \hfill
    \includegraphics[width=0.49\textwidth]{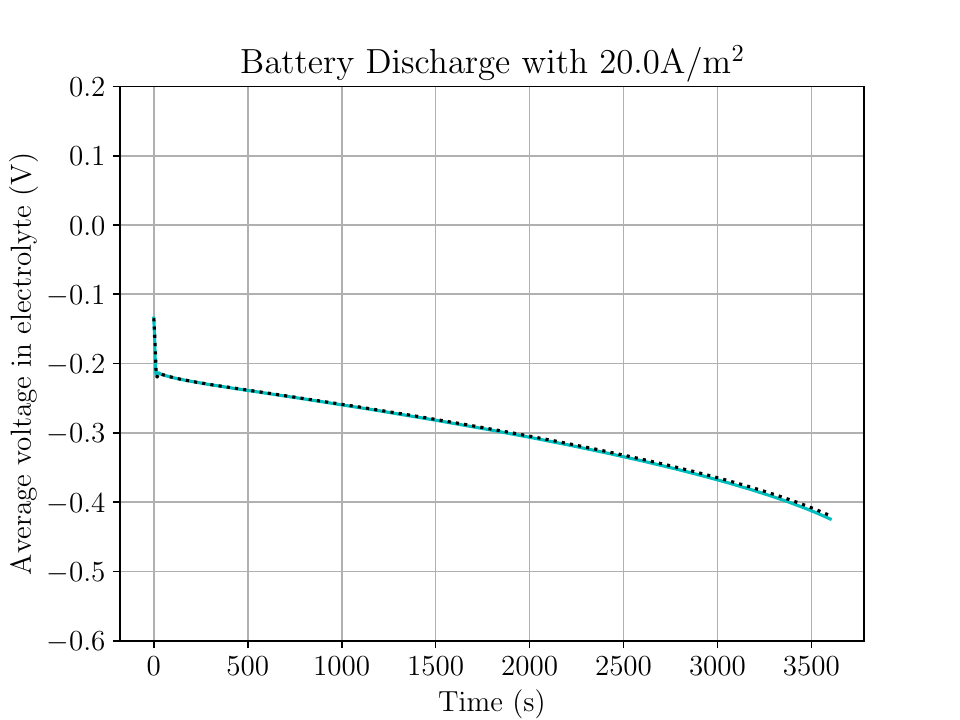} \\
    \includegraphics[width=0.49\textwidth]{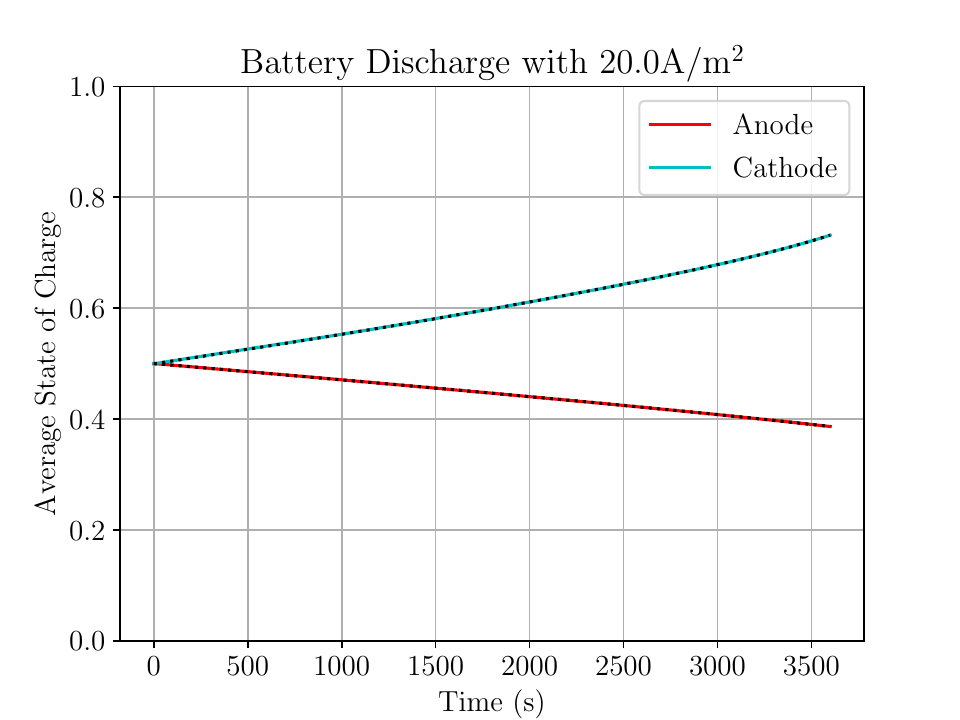} \hfill
    \includegraphics[width=0.49\textwidth]{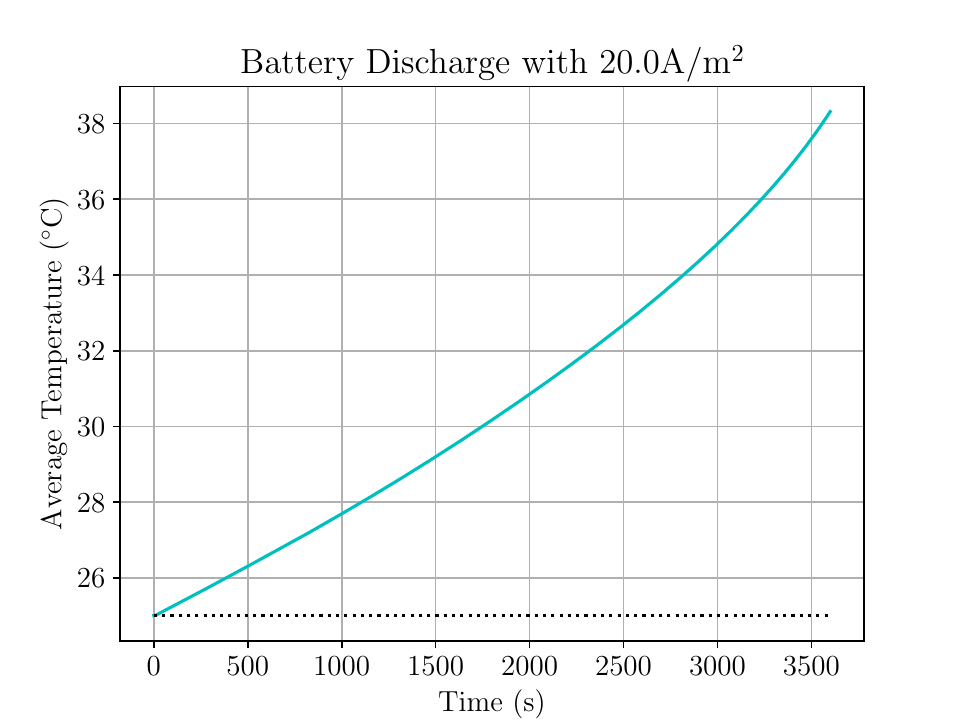}
    \caption{High current discharge scenario. 
    Results comparison between full multiphysical model (solid lines)
    and electrochemical model (dotted lines): 
    (top left) output voltage vs. time; 
    (top right) average potential in electrolyte vs. time;
    (bottom left) average electrode state of charge vs. time;
    (bottom right) average temperature (in Celsius degrees) vs. time}
    \label{fig:b2dh_comparison_plots}
\end{figure}

\begin{figure}[!ht]
    \includegraphics[width=0.49\textwidth]{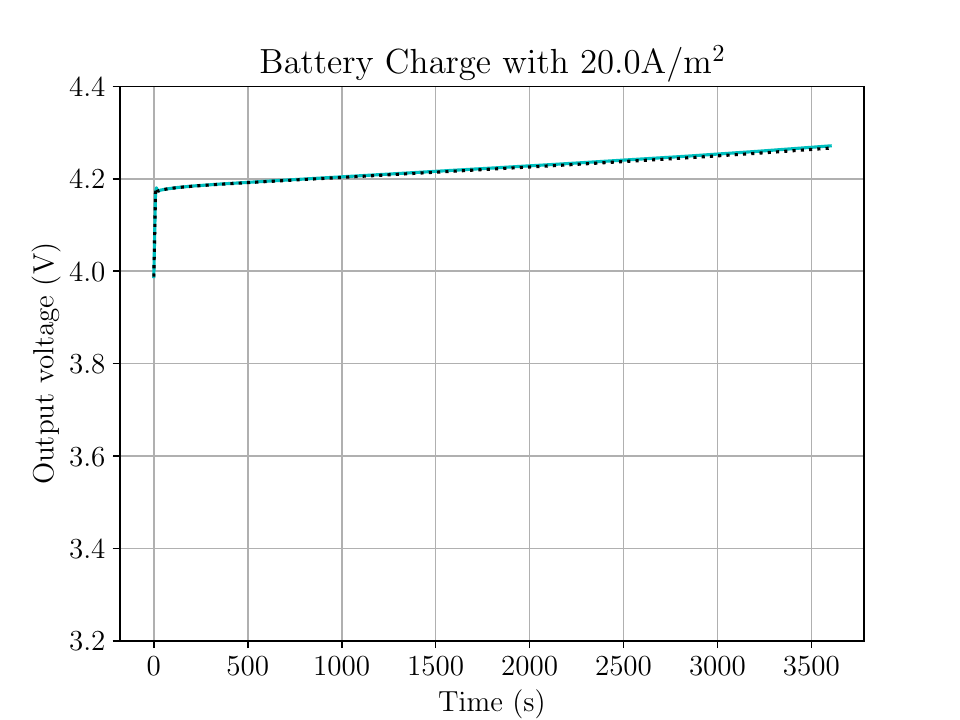} \hfill
    \includegraphics[width=0.49\textwidth]{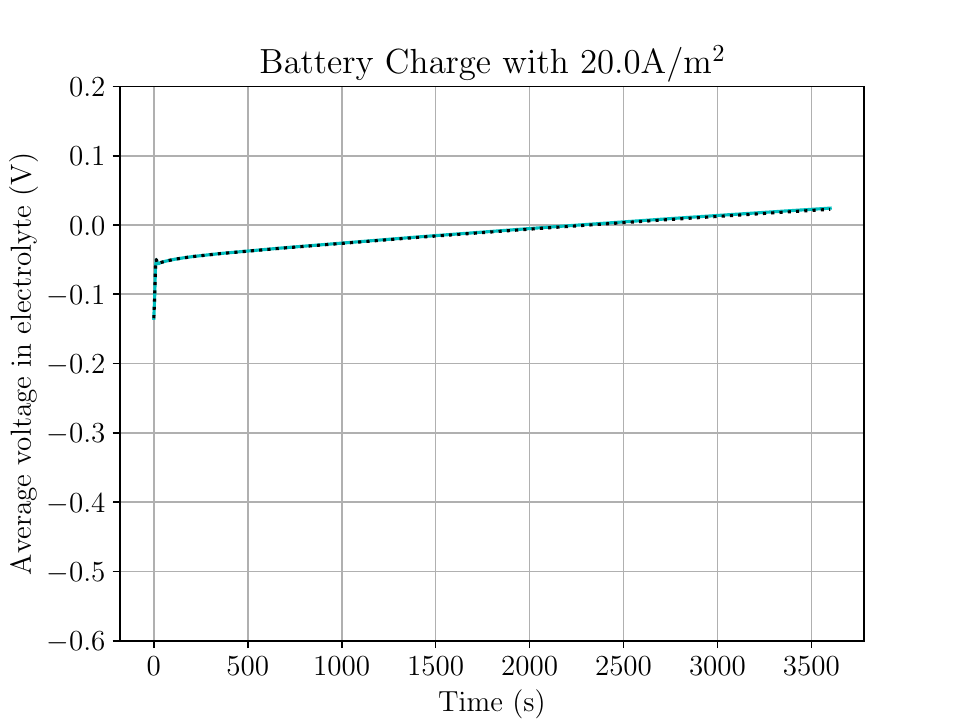} \\
    \includegraphics[width=0.49\textwidth]{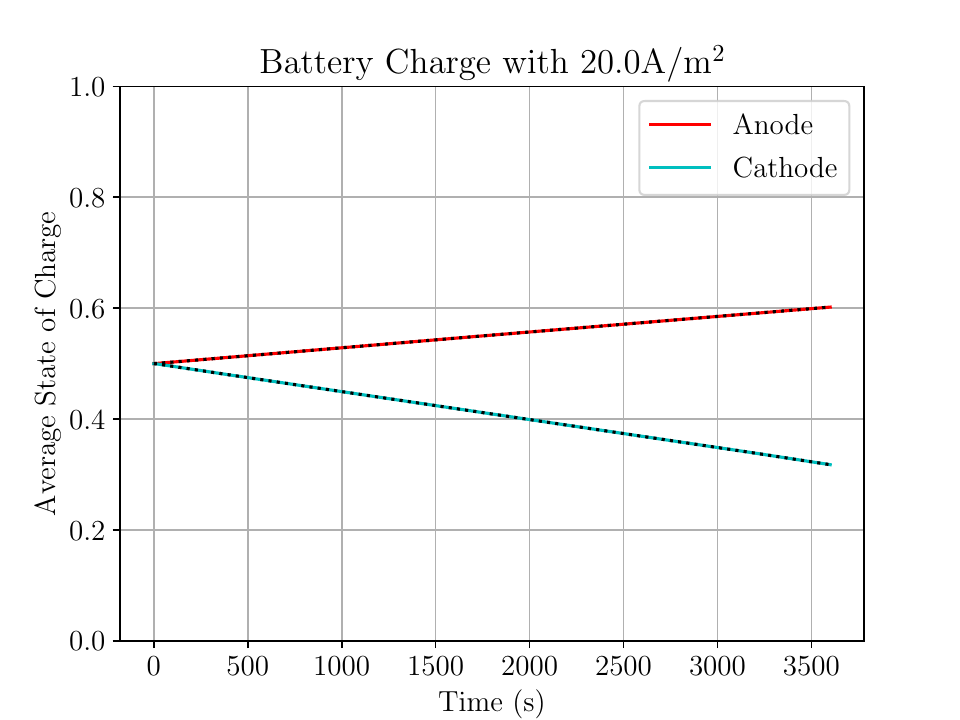} \hfill
    \includegraphics[width=0.49\textwidth]{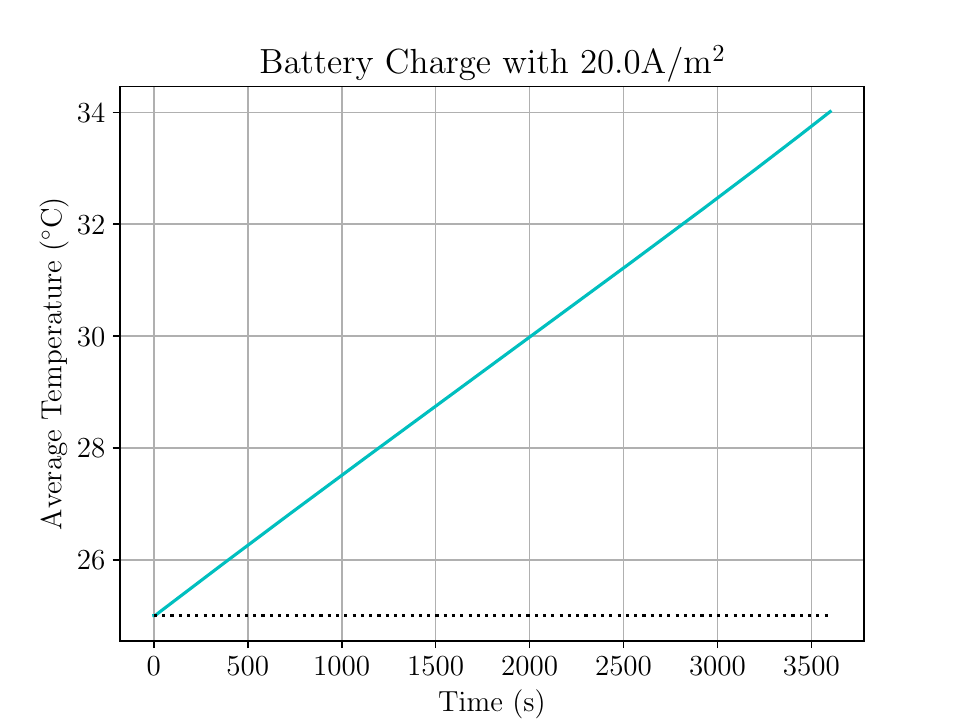}
    \caption{High current charge scenario. 
    Results comparison between full multiphysical model (solid lines)
    and electrochemical model (dotted lines): 
    (top left) output voltage vs. time; 
    (top right) average potential in electrolyte vs. time;
    (bottom left) average electrode state of charge vs. time;
    (bottom right) average temperature (in Celsius degrees) vs. time}
    \label{fig:b2ch_comparison_plots}
\end{figure}

After computing with the model and post-processing the results, we are able to present time plots and field plots that depict the behavior of the
battery discharge and charge processes. Taking advantage of formulas \eqref{eq:v-out} through \eqref{eq:tmprt-avg}, in Figures \ref{fig:b2dl_comparison_plots} and \ref{fig:b2cl_comparison_plots} we graph the evolution in time for the cell output voltage $V_{\text{out}}(t)$, the average electrolyte voltage $\bar{\phi_e}(t)$, the states of charge $\bar{\check{c}_{sa}}$ and $\bar{\check{c}_{sa}}$, and the average cell temperature $\bar{\theta}(t)$ under the low current discharge and charge, respectively. The plots for each model are distinguished by the line type: solid lines for multiphysics model, dotted lines for electrochemical model. In a similar fashion, Figures \ref{fig:b2dh_comparison_plots}
and \ref{fig:b2ch_comparison_plots} correspond to the high current scenarios.
%
%

What we can observe here is that the behavior changes quite slightly if just the electrochemical model is used. Particularly, the state of charge of the electrodes seems to be indistinguishable one from another. Where a very small difference shows up is in the last portion of the output voltage graphs, but only in the high current scenarios. We can say that, when accounting for all the physics, a slightly shorter time is predicted for the discharge (charge) process to reach a certain low (high) voltage. To be more precise, in the hour of simulated discharge at 20 A/m$^2$, the electrochemical model reached its lowest output potential about 45 seconds later than the full model, while in the charge scenario, the highest voltage of the electrochemical model was met nearly two and a half minutes earlier by the multiphysical model. We expect this to be more significant and relevant as the applied current density increases, so that thermo-mechanical effects should not be neglected.

Another result to point out is the sharp electric potential drop/surge produced right at the start of the discharge/charge process. Next, we attempt to explain the reasons for that to happen. We have chosen the initialization of the electric potential fields for maintaining equilibrium if we have zero current applied to the cathode's current collector (see subsection \ref{subsec:model_ic}). But, when $\iapp\neq0$, then by Gauss' law \eqref{eq:gauss}, in order to make the divergence of $\bsi_{sc}$ vanish, another flux of the same magnitude but opposite sign must be present at other portion of the cathode's boundary. The only such place allowed by our problem's BC is the interface, and that is why a nonzero overpotential is immediately induced, thus making both $\phi_{sc}$ and $\phi_e$ vary such that the corresponding Butler-Volmer current $\ibv$ can balance out $\iapp$.

It is of interest as well to capture the average temperature rise as the discharge and charge goes on, displayed in the bottom left plots of Figures \ref{fig:b2dl_comparison_plots}-\ref{fig:b2ch_comparison_plots}. The dotted lines are horizontal because in the electrochemical model there is no 
evolution of temperature (it is assumed to be isothermal). We note that both discharge scenarios are associated with a greater temperature variation than their charging counterparts. Particularly, the high current discharge presents an increase of over 13 degrees after an hour.

\begin{figure}[!ht]
    \begin{subfigure}{0.49\textwidth}
        \centering
        \includegraphics[width=\textwidth]{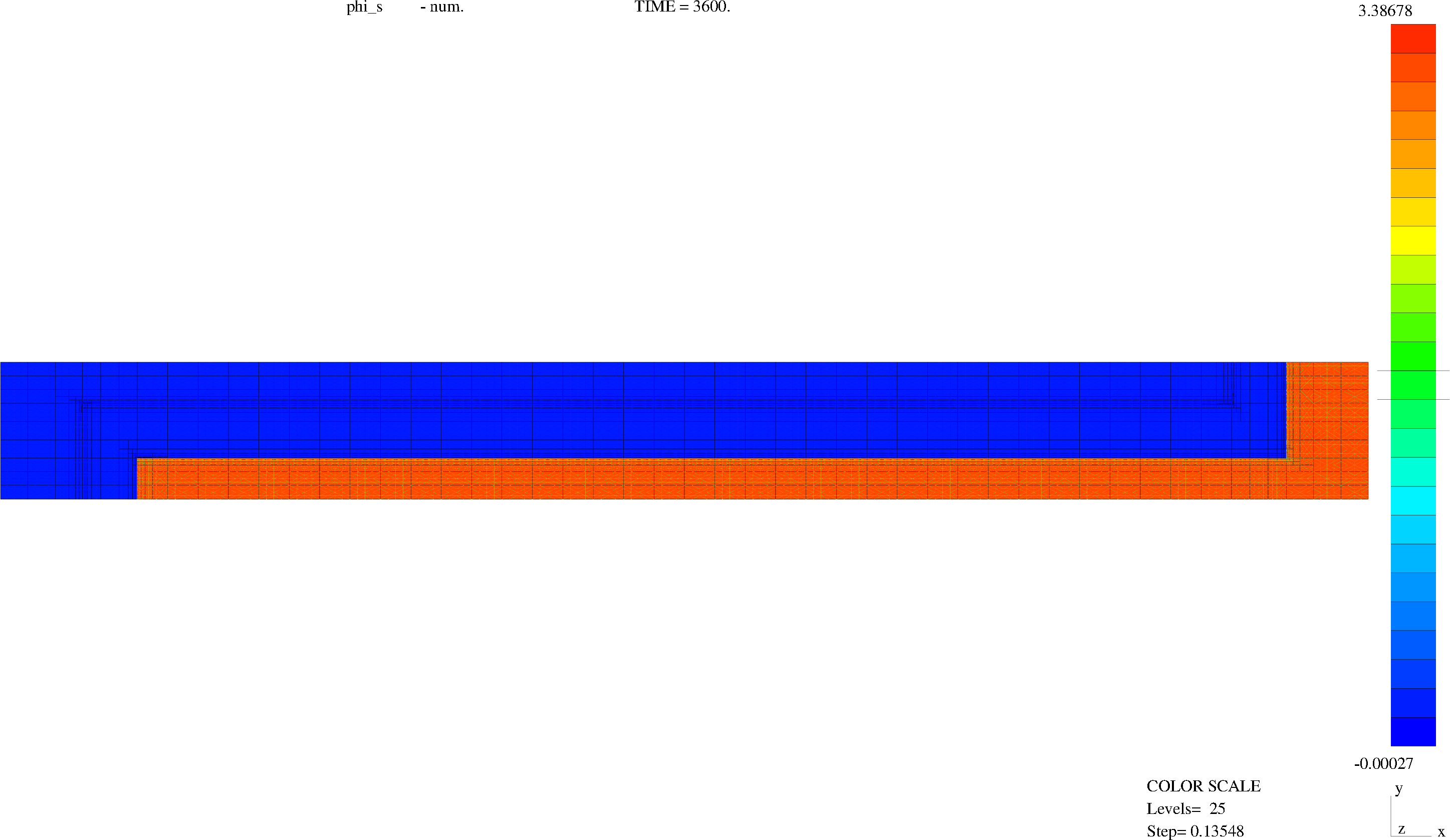}
        \caption{}
    \end{subfigure}
    \begin{subfigure}{0.49\textwidth}
        \centering
        \includegraphics[width=\textwidth]{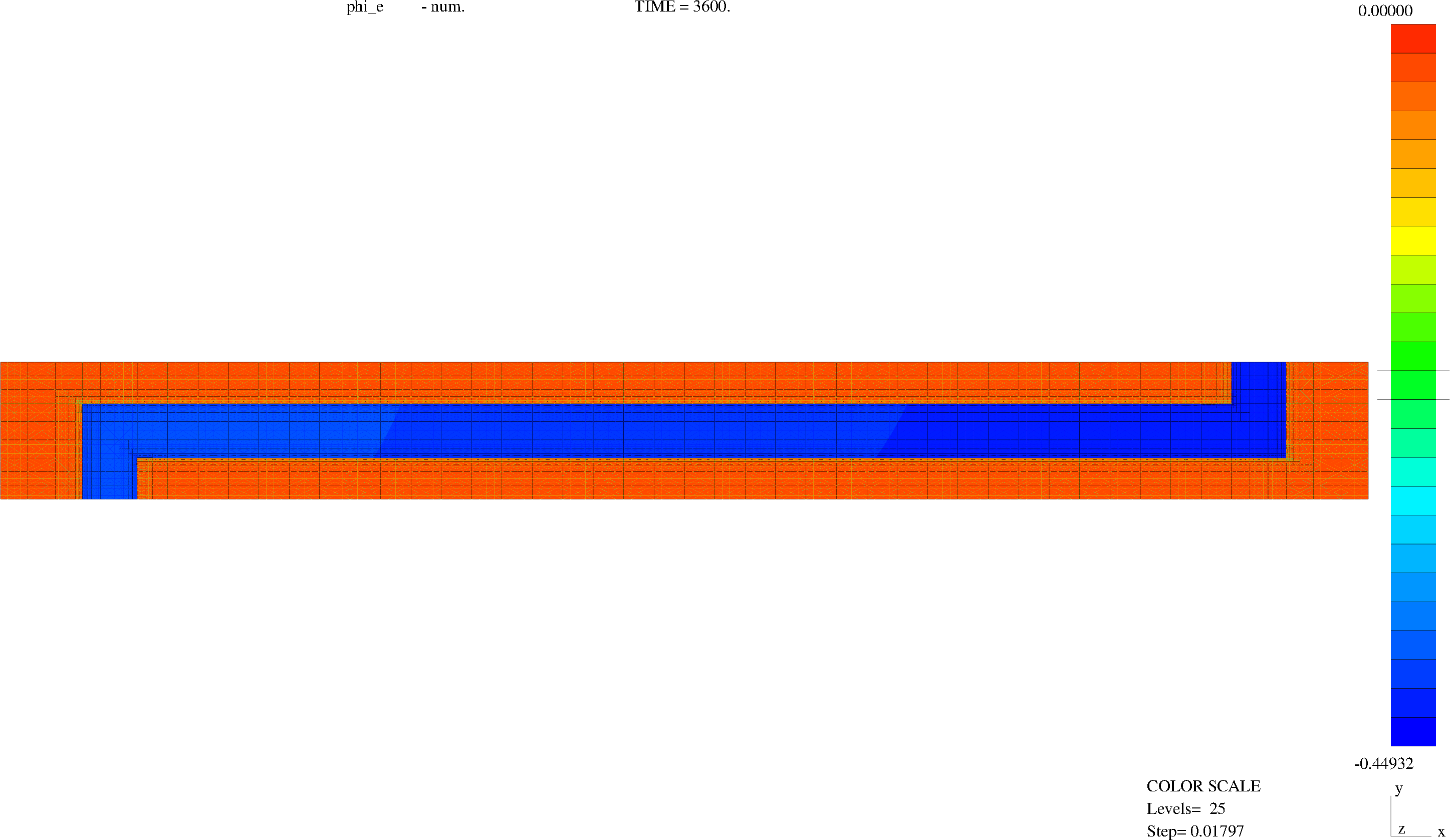}
        \caption{}
    \end{subfigure} \\[12pt]
    \begin{subfigure}{0.49\textwidth}
        \centering
        \includegraphics[width=\textwidth]{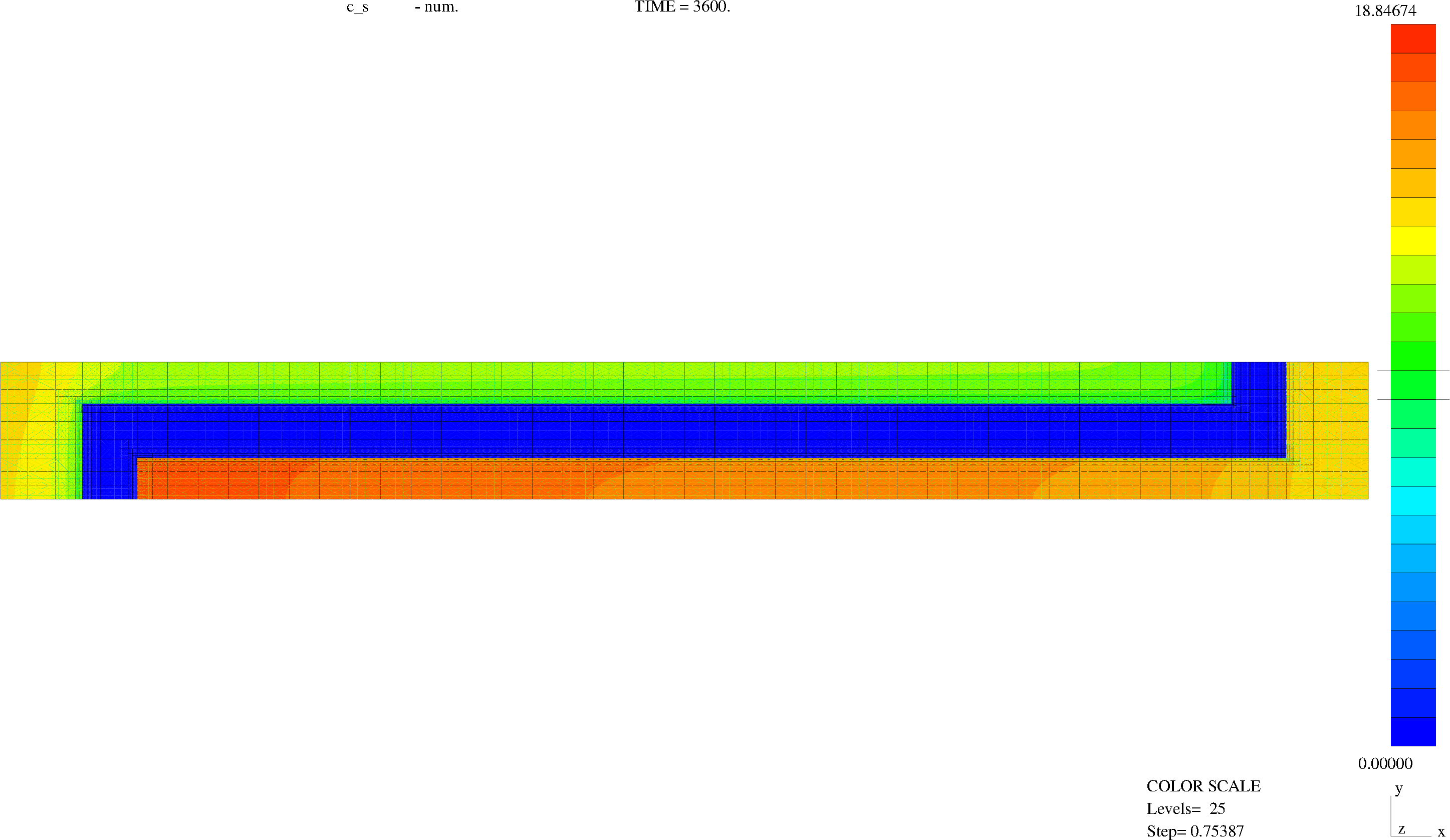}
        \caption{}
    \end{subfigure}
    \begin{subfigure}{0.49\textwidth}
        \centering
        \includegraphics[width=\textwidth]{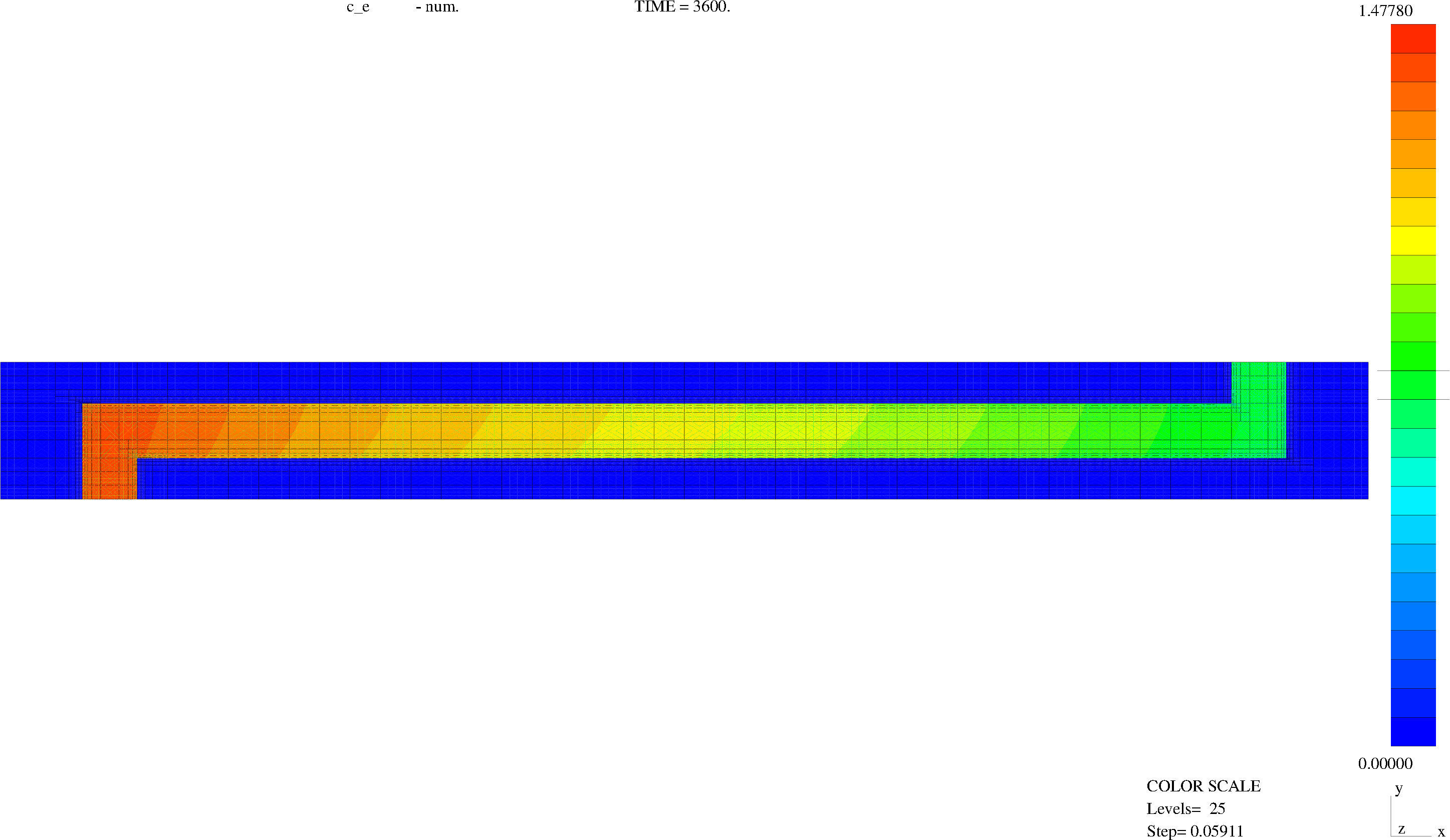}
        \caption{}
    \end{subfigure} \\[12pt]
    \begin{subfigure}{0.49\textwidth}
        \centering
        \includegraphics[width=\textwidth]{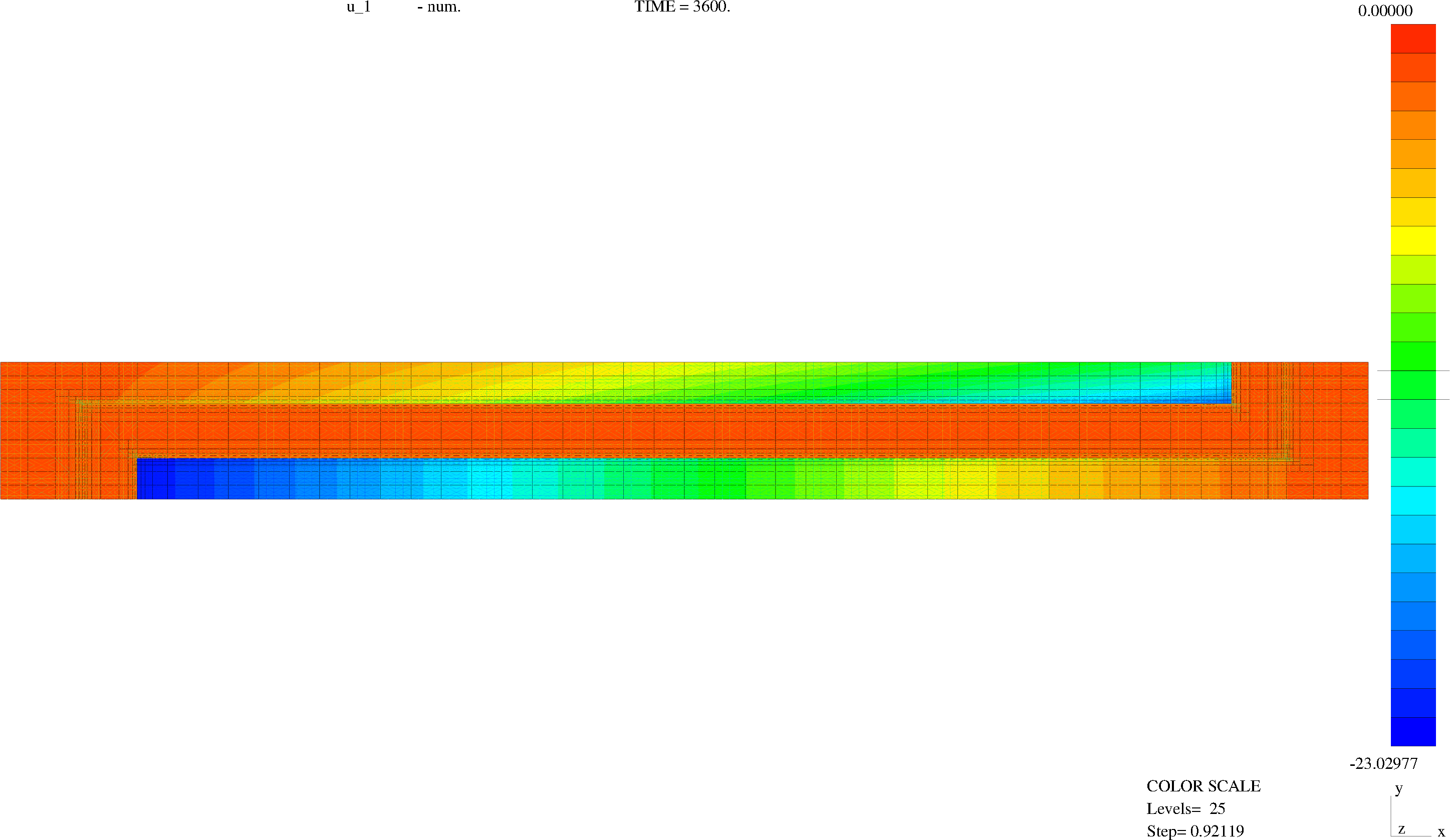}
        \caption{}
    \end{subfigure}
    \begin{subfigure}{0.49\textwidth}
        \centering
        \includegraphics[width=\textwidth]{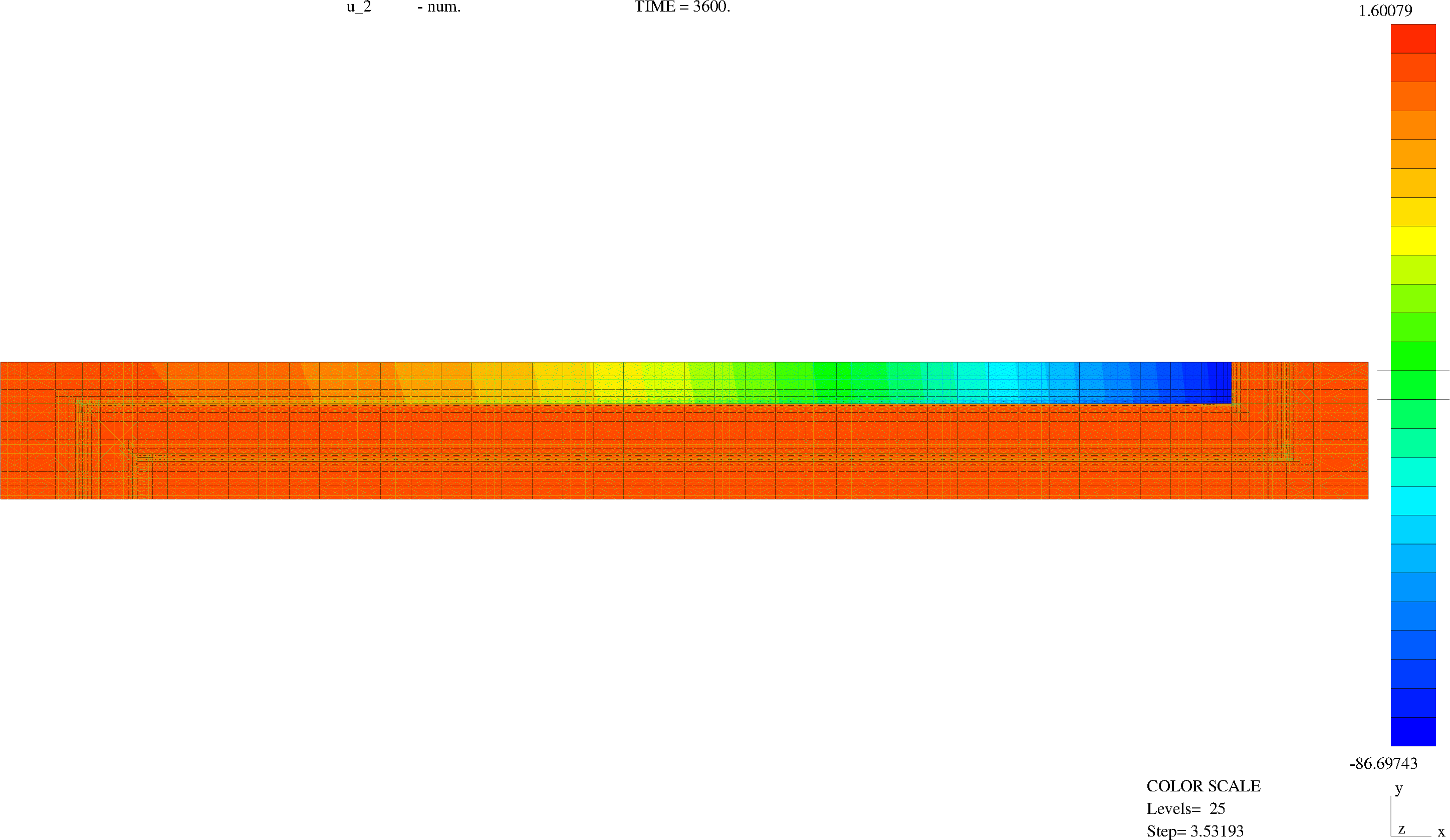}
        \caption{}
    \end{subfigure} \\[12pt]
    \begin{subfigure}{0.49\textwidth}
        \centering
        \includegraphics[width=\textwidth]{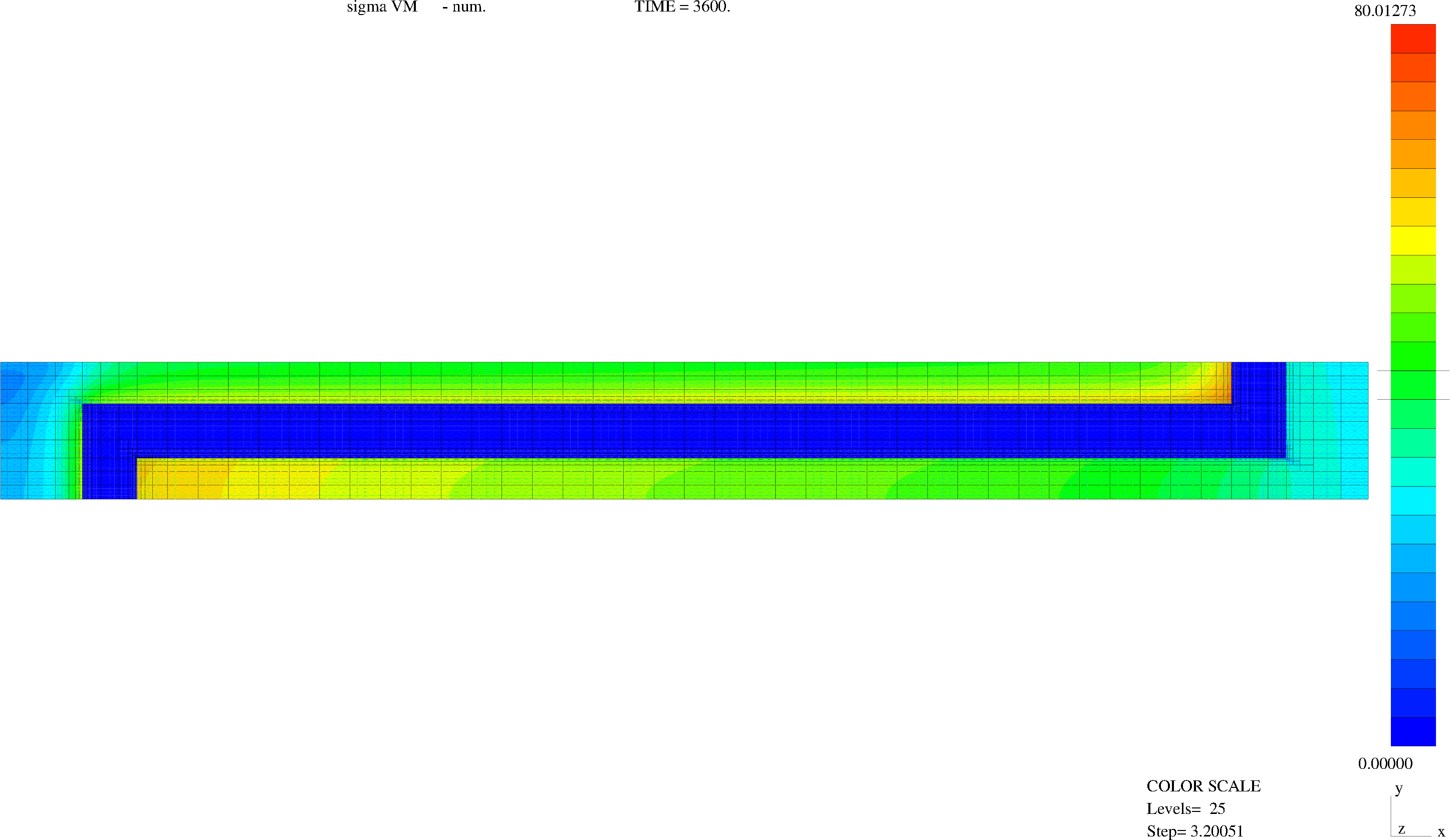}
        \caption{}
    \end{subfigure}
    \begin{subfigure}{0.49\textwidth}
        \centering
        \includegraphics[width=\textwidth]{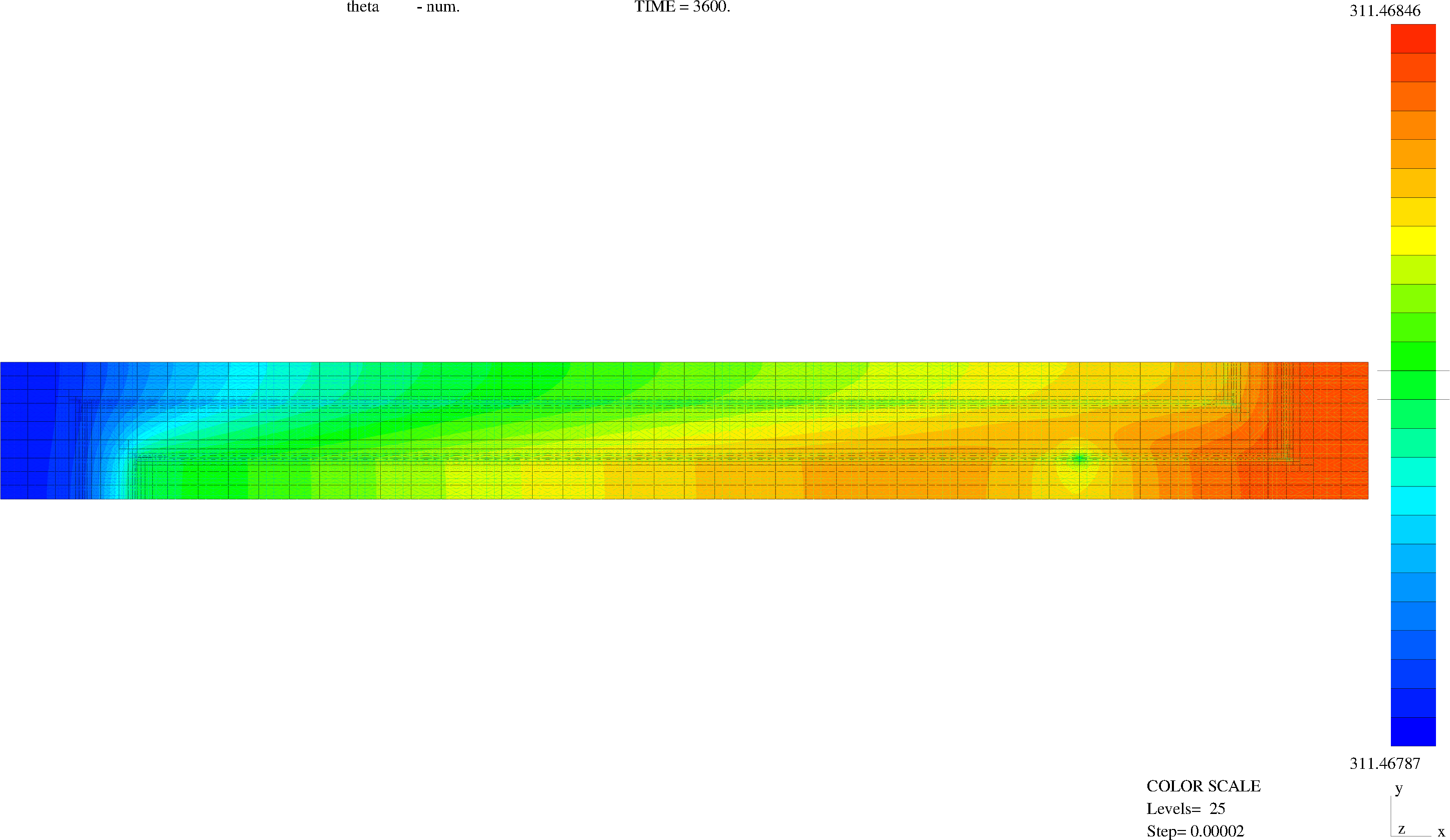}
        \caption{}
    \end{subfigure}
    \caption{High current discharge scenario - full multiphysical model. 
    Distribution of field variables in the representative domain at end time, $t$=3600s = 1 hour:
    (a) $\phi_{sc}$ in Volts; 
    (b) $\phi_e$ in Volts; 
    (c) $c_s$ in mol/dm$^3$; 
    (d) $c_e$ in mol/dm$^3$; 
    (e) $u_1$ in {\textmu}m;  
    (f) $u_2$ in {\textmu}m; 
    (g) $\sigma_{_{VM}}$ in MPa;
    (g) $\theta$ in Kelvin.}
    \label{fig:b2dh4_field_plots}
\end{figure}

\begin{figure}[!ht]
    \begin{subfigure}{0.49\textwidth}
        \centering
        \includegraphics[width=\textwidth]{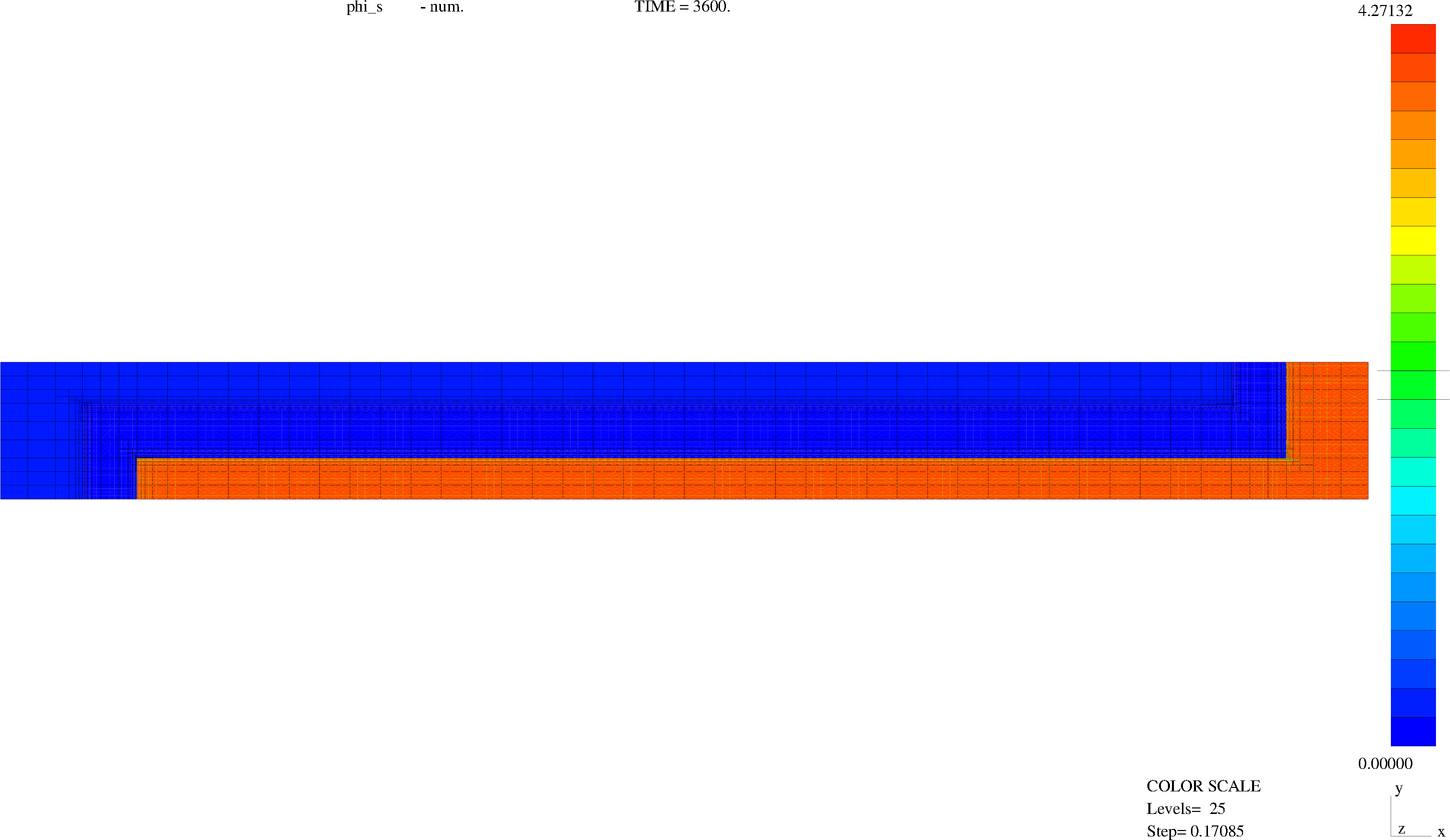}
        \caption{}
    \end{subfigure}
    \begin{subfigure}{0.49\textwidth}
        \centering
        \includegraphics[width=\textwidth]{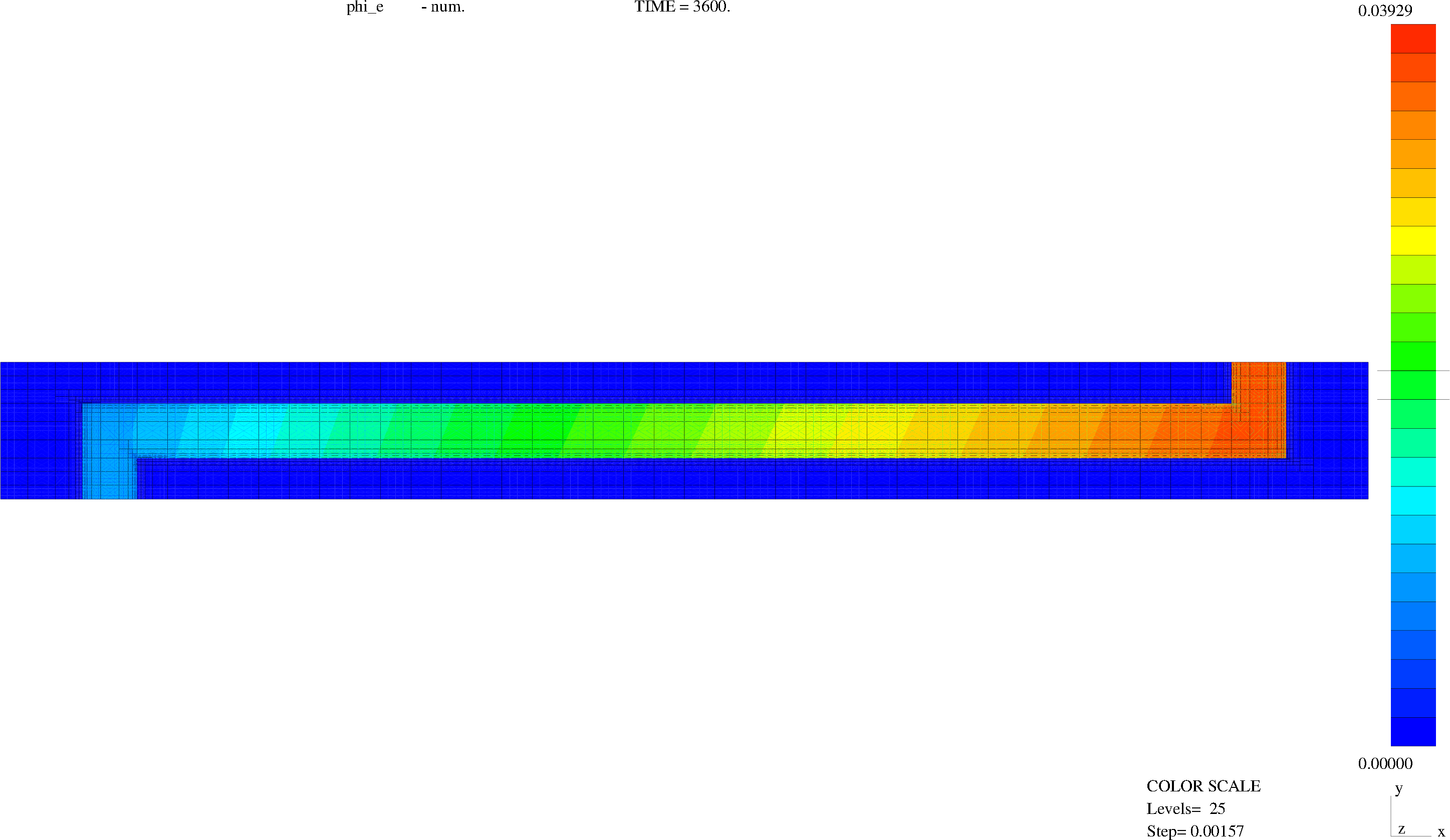}
        \caption{}
    \end{subfigure} \\[12pt]
    \begin{subfigure}{0.49\textwidth}
        \centering
        \includegraphics[width=\textwidth]{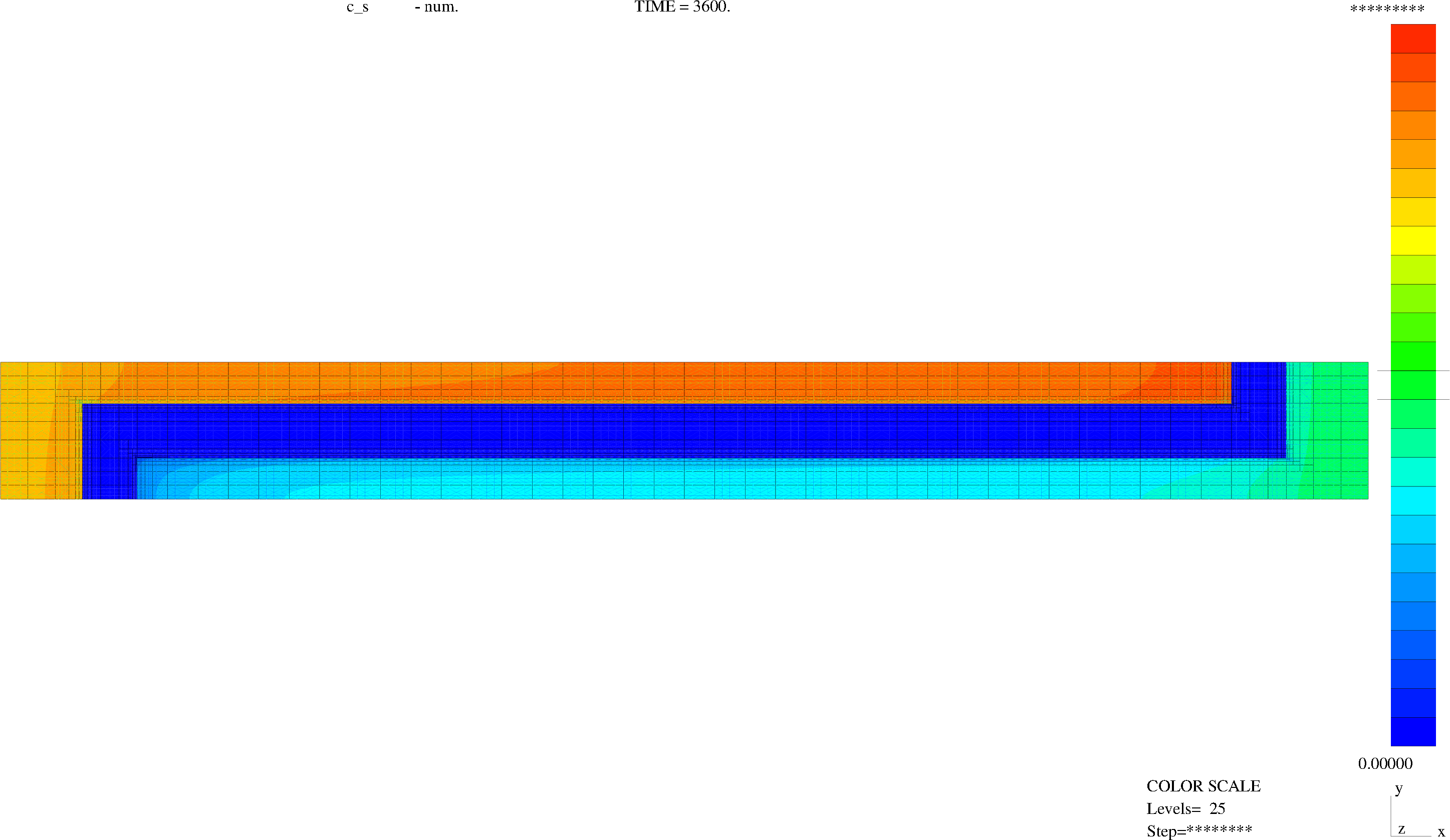}
        \caption{}
    \end{subfigure}
    \begin{subfigure}{0.49\textwidth}
        \centering
        \includegraphics[width=\textwidth]{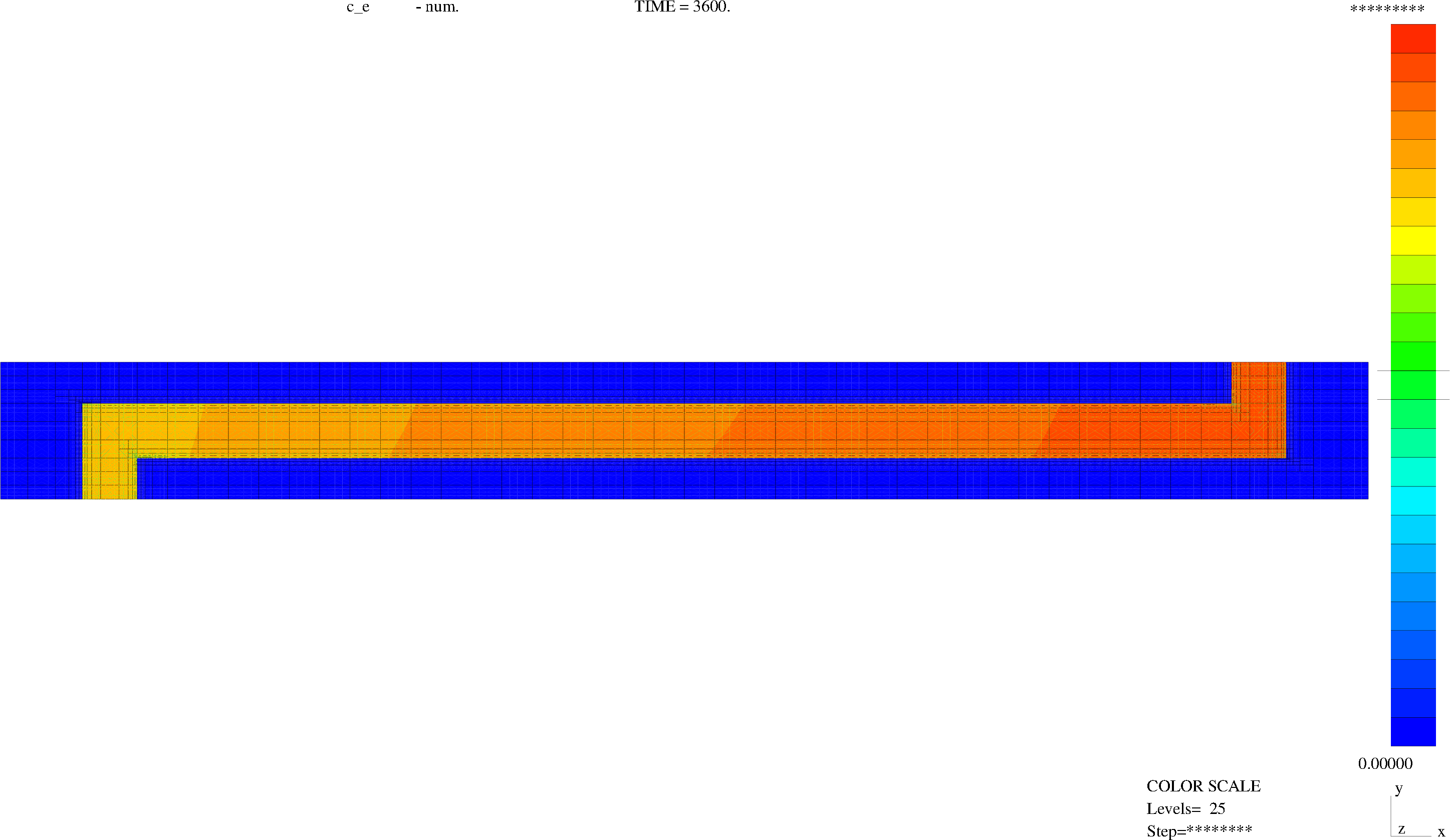}
        \caption{}
    \end{subfigure} \\[12pt]
    \begin{subfigure}{0.49\textwidth}
        \centering
        \includegraphics[width=\textwidth]{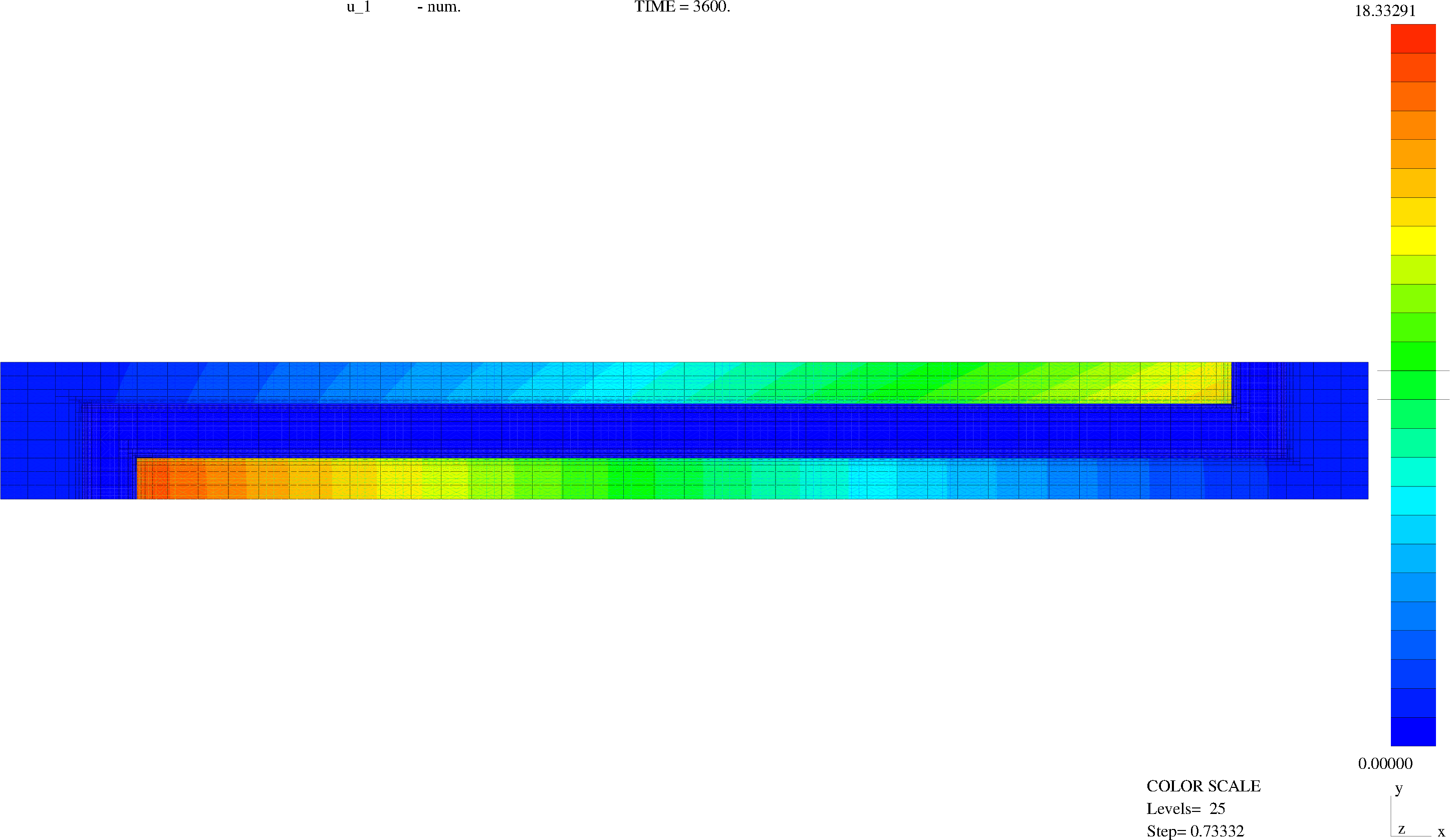}
        \caption{}
    \end{subfigure}
    \begin{subfigure}{0.49\textwidth}
        \centering
        \includegraphics[width=\textwidth]{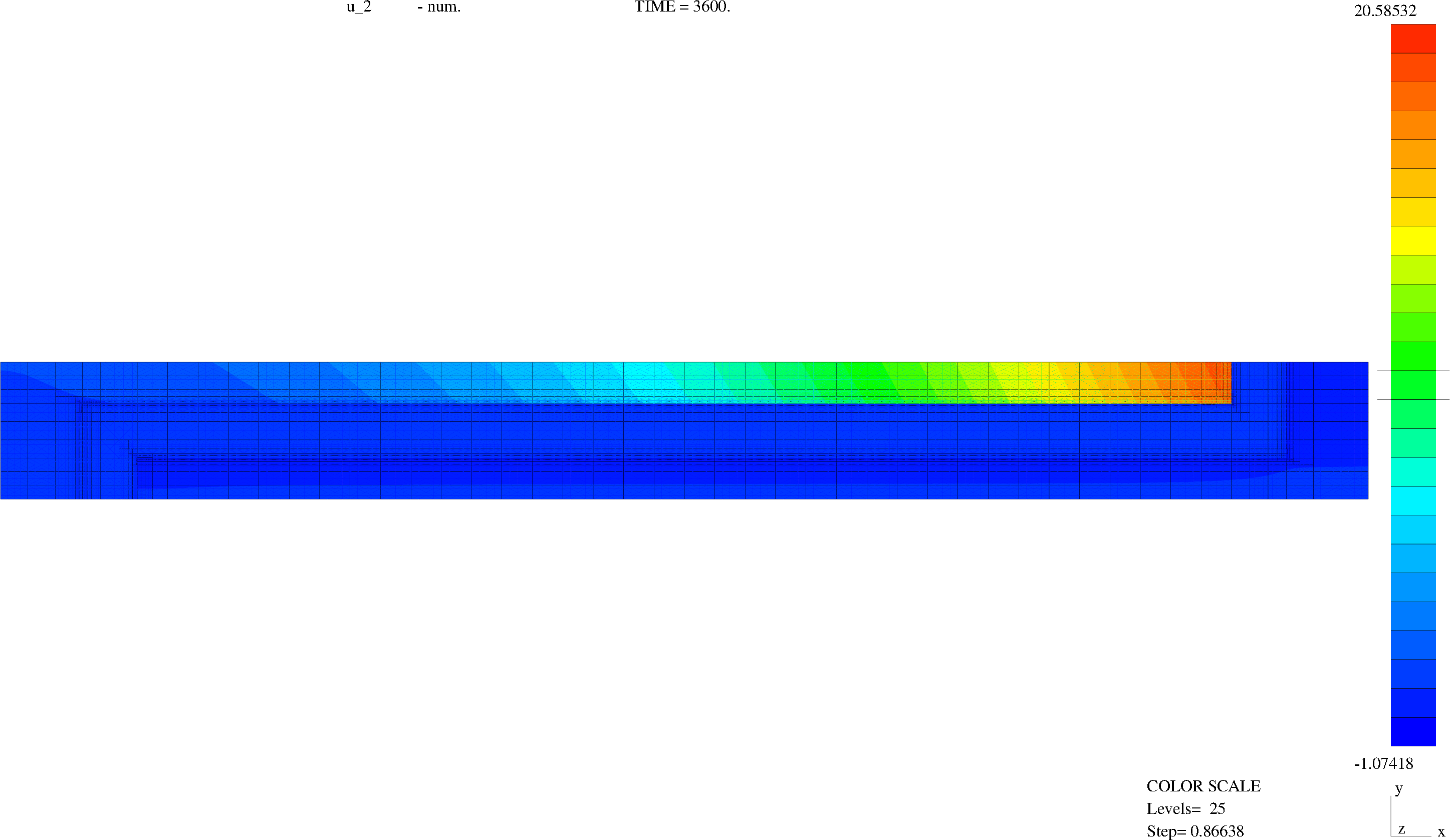}
        \caption{}
    \end{subfigure} \\[12pt]
    \begin{subfigure}{0.49\textwidth}
        \centering
        \includegraphics[width=\textwidth]{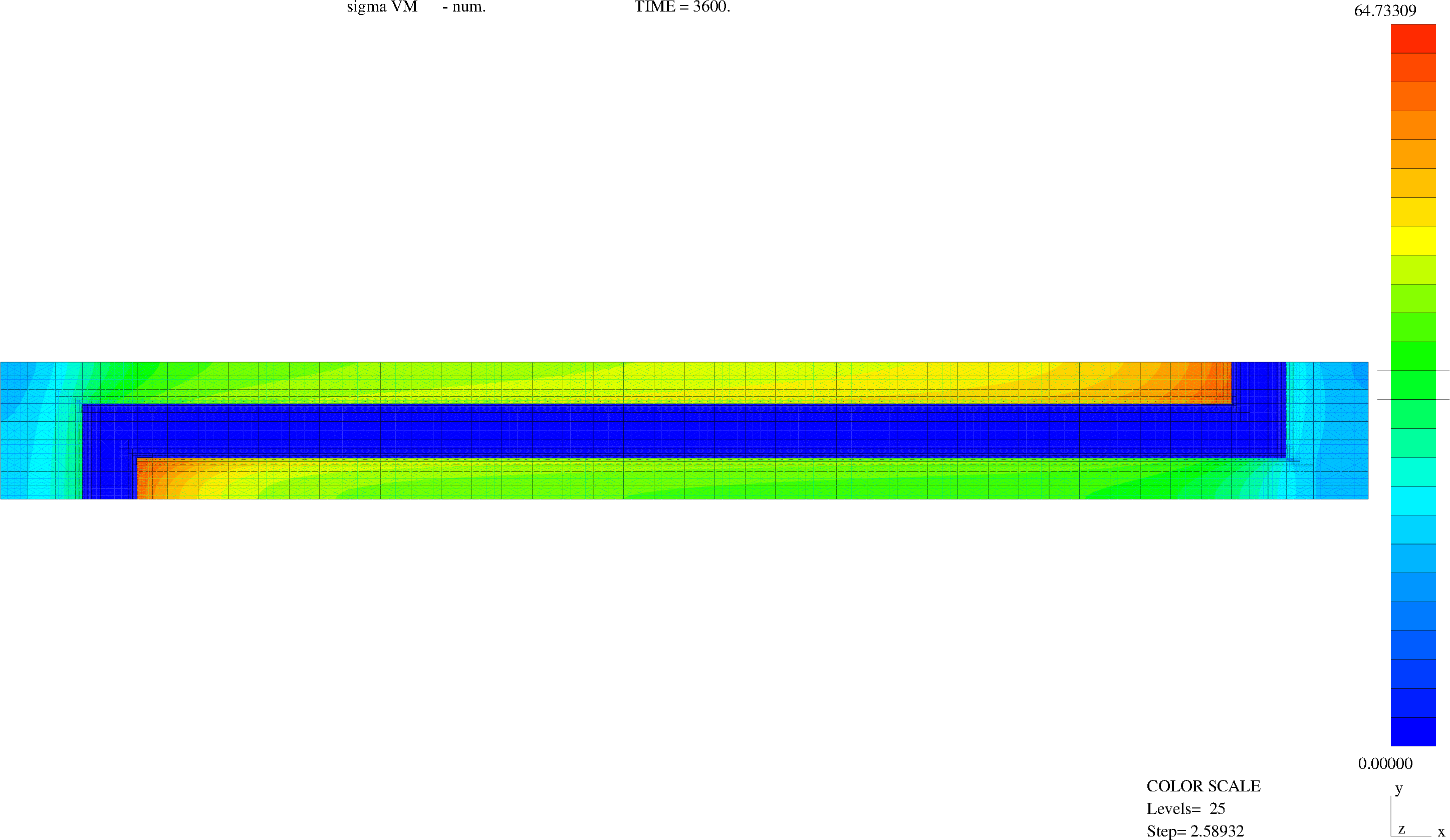}
        \caption{}
    \end{subfigure}
    \begin{subfigure}{0.49\textwidth}
        \centering
        \includegraphics[width=\textwidth]{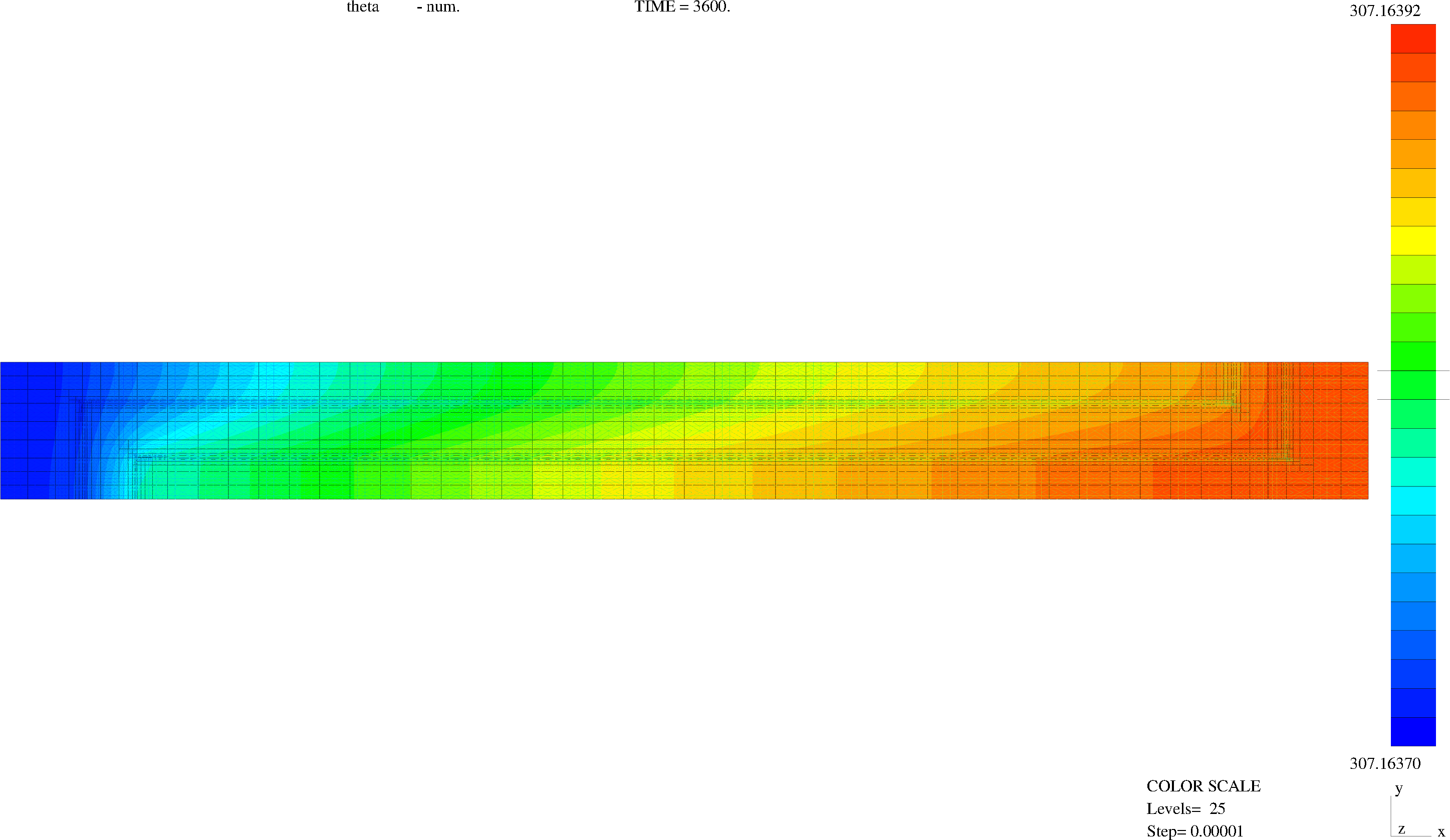}
        \caption{}
    \end{subfigure}
    \caption{High current charge scenario - full multiphysical model. 
    Distribution of field variables in the representative domain at end time, $t$=3600s = 1 hour:
    (a) $\phi_{sc}$ in Volts; 
    (b) $\phi_e$ in Volts; 
    (c) $c_s$ in mol/dm$^3$; 
    (d) $c_e$ in mol/dm$^3$; 
    (e) $u_1$ in {\textmu}m;  
    (f) $u_2$ in {\textmu}m; 
    (g) $\sigma_{_{VM}}$ in MPa;
    (g) $\theta$ in Kelvin.}
    \label{fig:b2ch4_field_plots}
\end{figure}

%
Regarding the spatial distribution of the field variables, 
we restrict the presentation of results to the full multiphysical model applied to both high current scenarios. 
Figures \ref{fig:b2dh4_field_plots} and \ref{fig:b2ch4_field_plots} 
consist of eight field plots each, that summarize the state of the representative cell domain at $t=\tend$,
for the discharge and charge cases, respectively. It is important to remark, that whenever the variable is not
defined over the entire domain $\Omega$, we deliver the obtained value at the variable support only, and zero elsewhere.
With that in mind, in the figures we find the distribution of:
electric potential in the solid electrodes, $\phi_s$ (in Volts);
electric potential in the electrolyte, $\phi_e$ (in Volts);
Lithium ion concentration in the solid electrodes, $c_s$ (in mol/dm$^3$);
Lithium ion concentration in the electrolyte, $c_e$ (in mol/dm$^3$);
displacement in coordinate 1, $u_1$ (in {\textmu}m);
displacement in coordinate 2, $u_2$ (in {\textmu}m);
Von Mises equivalent tensile stress $\sigma_{_{VM}}$ (in MPa);
and temperature, $\theta$ (in Kelvin).
%

Although the first four pictures of each figure just supplement the information of the time plots at $t=\tend$, namely, the spatial distribution of the electric potential and concentration variables therein averaged, the field plots provide new information, especially regarding the mechanical and thermal aspects of the system. For instance, we see that the displacements, which are induced by thermal dilation and Lithium intercalation, are far from negligible, since we see values of up to 23{\textmu}m horizontally, and 86{\textmu}m vertically, considering that the domain is 1000 {\textmu}m long and 100 {\textmu}m wide. On the same line, the state of stresses, summarized by the Von Mises stress reach very significant values, especially around corners, with a maximum of 80Mpa in the discharge, and 65MPa in the charge scenario. These mechanical insights can be relevant to a designer in order to prevent material failure.

The last plot in Figures \ref{fig:b2dh4_field_plots} and \ref{fig:b2ch4_field_plots} demonstrates that the temperature distribution is really smooth over all the subdomains of the cell, but with an almost negligible gradient. We can assert that, although the temperature increase in time is significant as a whole, the entire domain practically gets heated uniformly. This makes it plausible, for a future work, to consider a model where an averaged temperature is solved for at the macro level instead of seeking a detailed solution for this variable.

At last, we use the notion of power density to make a final comparison between the two models considered in the numerical experiments. Table \ref{tab:power_comparison} includes the result of applying \eqref{eq:power-energy-density} to the results of every scenario and both models. There is no doubt that the difference between both models is minute for all scenarios, but the relative difference in the high current scenarios has increased importantly with respect to the low current scenarios. We anticipate that, when considering much higher values of $\iapp$, the difference in power density may prove to be large enough so as to prefer the multiphysical model's higher fidelity over the simpler electrochemical model, whenever a decision has to be made concerning design, manufacturing or selection of batteries.
\begin{table}[!ht]
\centering
\begin{tabular}{|c|c|c|c|}
\hline
\textbf{Scenario}                                & \textbf{Model}           & $\overline{P}$ [W/dm$^3$] & \textbf{Relative difference} \\ \hline
\multirow{2}{*}{\textbf{Low current discharge}}  & \textbf{Multiphysical}   & 0.89517                     & \multirow{2}{*}{0.002\%}     \\ \cline{2-3}
                                                 & \textbf{Electrochemical} & 0.89519                     &                              \\ \hline
\multirow{2}{*}{\textbf{Low current charge}}     & \textbf{Multiphysical}   & -0.94241                    & \multirow{2}{*}{-0.003\%}    \\ \cline{2-3}
                                                 & \textbf{Electrochemical} & -0.94239                    &                              \\ \hline
\multirow{2}{*}{\textbf{High current discharge}} & \textbf{Multiphysical}   & 3.34967                     & \multirow{2}{*}{0.106\%}     \\ \cline{2-3}
                                                 & \textbf{Electrochemical} & 3.35322                     &                              \\ \hline
\multirow{2}{*}{\textbf{High current charge}}    & \textbf{Multiphysical}   & -3.91057                    & \multirow{2}{*}{-0.054\%}    \\ \cline{2-3}
                                                 & \textbf{Electrochemical} & -3.90844                    &                              \\ \hline
\end{tabular}
\caption{Power density and Energy density results comparison.}
\label{tab:power_comparison}
\end{table}

\section{Conclusions} 
\label{sec:conclusions}
\subsection{Summary and contributions}
We have presented a thorough derivation of a multiphysical model for Lithium-ion batteries, 
based on continuum theories, under few assumptions on the battery cell's materials and configuration. 
The resulting nonlinear PDE system along with suitable initial conditions, 
boundary conditions and interface conditions is turned into a discrete model,
by applying a Galerkin discretization to the weak form of the problem, 
along with a stable second order semi-implicit time stepping technique.

The numerical experiments have delivered a whole set of results that may help enlightening 
the understanding of the complex physics going on inside a battery cell, 
which is aided by the contrast of the full multiphysical model with a simpler electrochemical model. 
As the load increases, the consideration of thermal and mechanical effects 
(only available in the full multiphysical model) becomes more relevant in the entire system. 
Moreover, the temperatures, deformations and stresses revealed by the simulations
can provide important insights for battery scientists and engineers and, 
in particular, may open up opportunities for design optimization.

Given the major complexity of the multi-way coupling present in the system of equations, 
the full thermo-electro-chemo-mechanical model and its corresponding discrete scheme constitute a solid development
for advanced and highly accurate multiphysics simulations, 
so becoming significant mathematical and computational contributions to the science and engineering of lithium-ion batteries.

\subsection{Future extensions}
While our work manages to capture the intricate coupling between physics, a noticeable limitation
is the relatively small values of the applied current density, $\iapp$,
if comparing to other sources which have handled up to several hundreds of A/m$^2$. 
Future work will expand the applicable range of $\iapp$ by applying load stepping, 
so that the overshoots or undershoots generated at the beginning become better controlled. 
Furthermore, we have identified viable code optimizations that will allow for faster solution of each time step, 
thus making it feasible to use much smaller $\delta t$, again improving the accuracy in time for larger loads.

We have introduced special layered meshes that improve the solution accuracy around the interface. 
The mesh, constructed right before the simulation begins, remains fixed along all time steps. 
While this approach has been effective, we plan to improve it by including automatic adaptive meshes 
driven by residual-based error indicators.
Two promising formulations to achieve this goal are discontinuous Petrov-Galerkin (DPG) methods or 
First Order System Least-squares (FOSLS) methods \cite{DPGOverview,goalorientedDPG,bochev2004least,bochev2009least}.

\paragraph{Acknowledgments.}
The author thanks the institutional research funding by \emph{Fundación Universitaria Konrad Lorenz}, in its 2022 and 2023 calls.
The author acknowledges the significant help by Jacob Salazar (UT Austin) with regards to the visualization of results and proofreading.

\bibliographystyle{apalike}
\bibliography{main}
\end{document}